\documentstyle[
12pt,
epsf,stmaryrd,amscd,amstex,amssymb]{amsart}
\newcommand{\nc}{\newcommand}
\newcommand{\rc}{\renewcommand}

\nc{\ad}{{\mbox{\bf{ad}}}}
\nc{\AJ}{{\operatorname{aj}}}
\nc{\Aut}{{\operatorname{Aut}}}
\nc{\Bls}{{{\cal B}ls}}
\nc{\Boxtimes}{{\fbox{$\times$}}}
\nc{\blt}{{\bullet}}
\nc{\bSt}{{\mbox{\bf{St}}}}
\nc{\card}{{\operatorname{card}}}
\nc{\Cch}{{\check{C}}}
\nc{\Ch}{{\operatorname{Ch}}}
\nc{\chara}{{\operatorname{char}}}
\nc{\CHom}{{\cal{H}om}}
\nc{\Coker}{{\operatorname{Coker}}}
\nc{\codim}{{\operatorname{codim}}}
\nc{\Cone}{{\operatorname{Cone}}}
\nc{\cSgn}{{\cal{S}gn}}
\nc{\depth}{{\operatorname{depth}}}
\nc{\dirlim}{{\underset{\rightarrow}{\operatorname{lim}}}}
\nc{\dotbox}{{\overset{\bullet}{\boxtimes}}}
\nc{\dotimes}{{\overset{\bullet}{\otimes}}}
\nc{\Ed}{{\operatorname{Edge}}}
\nc{\Ext}{{\operatorname{Ext}}}
\nc{\Fac}{{\cal{F}ac}}
\nc{\Fun}{{\operatorname{F}}}
\nc{\FS}{{\cal{FS}}}
\nc{\Hom}{{\operatorname{Hom}}}
\nc{\had}{{{\hat{\mbox{\bf{ad}}}}}}
\nc{\hgt}{{\operatorname{ht}}}
\nc{\Id}{{\operatorname{Id}}}
\nc{\id}{{\operatorname{id}}}
\nc{\Ima}{{\operatorname{Im}}}
\nc{\ind}{{\operatorname{ind}}}
\nc{\Ind}{{\operatorname{Ind}}}
\nc{\infi}{{\operatorname{inf}}}
\nc{\infh}{{\frac{\infty}{2}}}
\nc{\invlim}{{\underset{\leftarrow}{\operatorname{lim}}}}
\nc{\Jac}{{{\cal J}ac}}
\nc{\Ker}{{\operatorname{Ker}}}
\nc{\lcm}{{\operatorname{lcm}}}
\nc{\Locsys}{{{\cal L}ocsys}}
\nc{\Map}{{{\cal M}ap}}
\nc{\Mor}{{\operatorname{Mor}}}
\nc{\MS}{{\cal{MS}}}
\nc{\Ob}{{\operatorname{Ob}}}
\nc{\opp}{{\operatorname{opp}}}
\nc{\Or}{{{\cal O}r}}
\nc{\Ord}{{{\cal O}rd}}
\nc{\Part}{{{\cal P}art}}
\nc{\PGL}{{\operatorname{PGL}}}
\nc{\Pic}{{\operatorname{Pic}}}
\nc{\Rep}{{{\cal{R}}ep}}
\nc{\rk}{{\operatorname{rk}}}
\nc{\Sets}{{{\cal{S}}ets}}
\nc{\Sew}{{{\cal{S}}ew}}
\nc{\sgn}{{\operatorname{sgn}}}
\nc{\Sh}{{{\cal S}h}}
\nc{\Sign}{{{\cal S}ign}}
\nc{\Spe}{{\mbox{\bf{Sp}}}}
\nc{\supr}{{\operatorname{sup}}}
\nc{\Supp}{{\operatorname{Supp}}}
\nc{\supp}{{\operatorname{supp}}}
\nc{\Teich}{{{\cal{T}}eich}}
\nc{\tFS}{{\widetilde{\cal{FS}}}}
\nc{\Tor}{{\operatorname{Tor}}}
\nc{\totimes}{{\tilde{\otimes}}}
\nc{\tr}{{\operatorname{tr}}}
\nc{\tRep}{{\widetilde{{\cal R}ep}}}
\nc{\tTeich}{{\widetilde{{\cal T}eich}}}
\nc{\Vect}{{{\cal V}ect}}
\nc{\Ve}{{\operatorname{Vert}}}
\nc{\wt}{{\widetilde}}

\nc{\bo}{{\mbox{\bf{0}}}}
\nc{\One}{{\mbox{\bf{1}}}}
\nc{\one}{{\mbox{\bf{1}}}}

\nc{\BA}{{\Bbb A}}
\nc{\bA}{{\overline{A}}}
\nc{\ba}{{\mbox{\bf{a}}}}
\nc{\baB}{{\overline{B}}}
\nc{\baeta}{{\bar{\eta}}}
\nc{\baJ}{{\bar{J}}}
\nc{\BB}{{\Bbb B}}
\nc{\bB}{{\mbox{\bf{B}}}}
\nc{\bc}{{\mbox{\bf{c}}}}
\nc{\bC}{{\overline{C}}}
\nc{\BC}{{\Bbb{C}}}
\nc{\bCC}{{\overline{\cal{C}}}}
\nc{\bCM}{{\overline{\cal{M}}}}
\nc{\bD}{{\bar{D}}}
\nc{\BD}{{\overline{D}}}
\nc{\bd}{{\mbox{\bf{d}}}}
\nc{\bp}{{\mbox{\bf{p}}}}
\nc{\BE}{{\overline{E}}}
\nc{\BF}{{\overline{F}}}
\nc{\bF}{{\mbox{\bf{F}}}}
\nc{\bg}{{\mbox{\bf{g}}}}
\nc{\bG}{{\mbox{\bf{G}}}}
\nc{\BG}{{\Bbb G}}
\nc{\bGamma}{{\overline{\Gamma}}}
\nc{\bbH}{{\bar{\mbox{\bf{H}}}}}
\nc{\bH}{{\mbox{\bf{H}}}}
\nc{\bI}{{\mbox{\bf{I}}}}
\nc{\bL}{{\mbox{\bf{L}}}}
\nc{\BL}{{\Bbb{L}}}
\nc{\blambda}{{\bar{\lambda}}}
\nc{\bM}{{\mbox{\bf{M}}}}
\nc{\bmu}{{\vec{\mu}}}
\nc{\bN}{{\mbox{\bf{N}}}}
\nc{\BN}{{\Bbb{N}}}
\nc{\bnu}{{\mbox{\boldmath{${\nu}$}}}}
\nc{\bof}{{\mbox{\bf{f}}}}
\nc{\BP}{{\Bbb P}}
\nc{\bP}{{\mbox{\bf{P}}}}
\nc{\BPO}{{\overset{\circ}{\BP}}}
\nc{\BQ}{{\Bbb Q}}
\nc{\bq}{{\mbox{\bf{q}}}}
\nc{\BR}{{\Bbb{R}}}
\nc{\bR}{{\mbox{\bf{R}}}}
\nc{\br}{{\mbox{\bf{r}}}}
\nc{\breta}{{\bar{\eta}}}
\nc{\bs}{{\mbox{\bf{s}}}}
\nc{\bS}{{\mbox{\bf{S}}}}
\nc{\bt}{{\mbox{\bf{t}}}}
\nc{\bU}{{\mbox{\bf{U}}}}
\nc{\bV}{{\mbox{\bf{V}}}}
\nc{\bu}{{\mbox{\bf{u}}}}
\nc{\BUpsilon}{{\bar{\Upsilon}}}
\nc{\bw}{{\mbox{\bf{w}}}}
\nc{\bx}{{\mbox{\bf{x}}}}
\nc{\bX}{{\mbox{\bf{X}}}}
\nc{\BZ}{{\Bbb{Z}}}
\nc{\bz}{{\mbox{\bf{z}}}}
\nc{\bZ}{{\mbox{\bf{Z}}}}
\nc{\bzero}{\mbox{\boldmath{$0$}}}

\nc{\CA}{{\cal A}}
\nc{\CAD}{{\overset{\bullet}{\cal{A}}}}
\nc{\CAO}{{\overset{\circ}{\cal{A}}}}
\nc{\CB}{{\cal B}}
\nc{\CalD}{{\cal D}}
\nc{\CE}{{\cal E}}
\nc{\CF}{{\cal F}}
\nc{\CG}{{\cal G}}
\nc{\CH}{{\cal H}}
\nc{\CI}{{\cal I}}
\nc{\CID}{{\overset{\bullet}{\cal{I}}}}
\nc{\CJ}{{\cal J}}
\nc{\CK}{{\cal K}}
\nc{\CL}{{\cal L}}
\nc{\CM}{{\cal M}}
\nc{\CN}{{\cal N}}
\nc{\CO}{{\cal O}}
\nc{\CP}{{\cal P}}
\nc{\CPO}{{\overset{\circ}{\cal{P}}}}
\nc{\CQ}{{\cal Q}}
\nc{\CR}{{\cal R}}
\nc{\CS}{{\cal S}}
\nc{\CT}{{\cal T}}
\nc{\CTD}{{\overset{\bullet}{\cal{T}}}}
\nc{\CTPO}{{\overset{\circ}{\cal{T}\cal{P}}}}
\nc{\CU}{{\cal{U}}}
\nc{\CV}{{\cal V}}
\nc{\CW}{{\cal W}}
\nc{\CX}{{\cal X}}
\nc{\CY}{{\cal Y}}
\nc{\CZ}{{\cal Z}}

\nc{\dCL}{{\overset{\bullet}{\cal{L}}}}
\nc{\dd}{{\operatorname{d}}}
\nc{\ddelta}{{\overset{\bullet}{\delta}}}
\nc{\dfu}{{\overset{\bullet}{\frak{u}}}}
\nc{\dlambda}{{\overset{\bullet}{\lambda}}}
\nc{\DO}{{\overset{\circ}{D}}}
\nc{\dpar}{{\partial}}
\nc{\dS}{{\overset{\bullet}{S}}}
\nc{\dT}{{\overset{\bullet}{T}}}

\nc{\etabar}{{\bar{\eta}}}
\nc{\sbar}{{\bar{s}}}
\nc{\Kbar}{{\bar{\KK}}}

\nc{\hCH}{{\hat{\cal{H}}}}
\nc{\hCI}{{\hat{\cal{I}}}}
\nc{\hfC}{{\hat{\frak{C}}}}
\nc{\hfg}{{\hat{\frak{g}}}}
\nc{\hL}{{\hat{L}}}
\nc{\OH}{{\overset{\circ}{H}}}
\nc{\hpsi}{{\hat{\psi}}}
\nc{\hx}{{\hat{x}}}

\nc{\jo}{{\overset{\circ}{j}}}

\nc{\phid}{{\overset{\bullet}{\phi}}}

\nc{\tA}{{\tilde{A}}}
\nc{\ta}{{\tilde{a}}}
\nc{\tB}{{\tilde{B}}}
\nc{\tb}{{\tilde{b}}}
\nc{\tBP}{{\tilde{\BP}}}
\nc{\tC}{{\tilde{C}}}
\nc{\tc}{{\tilde{c}}}
\nc{\tCA}{{\tilde{\cal{A}}}}
\nc{\tCC}{{\tilde{\cal{C}}}}
\nc{\tCH}{{\tilde{\cal{H}}}}
\nc{\tCI}{{\tilde{\cal{I}}}}
\nc{\tCO}{{\tilde{\cal{O}}}}
\nc{\tCP}{{\tilde{\cal{P}}}}
\nc{\tCT}{{\tilde{\cal{T}}}}
\nc{\tD}{{\tilde{D}}}
\nc{\tDelta}{{\tilde{\Delta}}}
\nc{\tE}{{\tilde E}}
\nc{\tF}{{\tilde F}}
\nc{\tfD}{{\tilde{\frak{D}}}}
\nc{\tfF}{{\tilde{\frak{F}}}}
\nc{\tff}{{\tilde{\frak{f}}}}
\nc{\tfu}{{\tilde{\frak{u}}}}
\nc{\tJ}{{\tilde{J}}}
\nc{\tj}{{\tilde{j}}}
\nc{\tK}{{\tilde K}}
\nc{\tL}{{\tilde{L}}}
\nc{\tM}{{\tilde{M}}}
\nc{\tP}{{\tilde{P}}}
\nc{\tPhi}{{\tilde{\Phi}}}
\nc{\tpi}{\tilde{\pi}}
\nc{\TPO}{{\overset{\circ}{T\BP}}}
\nc{\tR}{{\tilde{R}}}
\nc{\tS}{{\tilde S}}
\nc{\tT}{{\tilde{T}}}
\nc{\ttau}{{\tilde{\tau}}}
\nc{\ttheta}{{\tilde{\theta}}}
\nc{\tU}{{\tilde{U}}}
\nc{\tUpsilon}{{\tilde{\Upsilon}}}
\nc{\tW}{{\tilde W}}
\nc{\ty}{{\tilde y}}
\nc{\tY}{{\tilde Y}}
\nc{\txi}{{\tilde{\xi}}}

\nc{\UD}{{\overset{\bullet}{U}}}
\nc{\UO}{{\overset{\circ}{U}}}

\nc{\vA}{{\vec{A}}}
\nc{\valpha}{{\vec{\alpha}}}
\nc{\vbeta}{{\vec{\beta}}}
\nc{\vc}{{\vec{c}}}
\nc{\vD}{{\vec{D}}}
\nc{\vd}{{\vec{d}}}
\nc{\vgamma}{{\vec{\gamma}}}
\nc{\vK}{{\vec{K}}}
\nc{\vlambda}{{\vec{\lambda}}}
\nc{\vmu}{{\vec{\mu}}}
\nc{\vnu}{{\vec{\nu}}}
\nc{\vo}{{\vec{0}}}
\nc{\vu}{{\vec{u}}}
\nc{\vx}{{\vec{x}}}
\nc{\vy}{\vec{y}}
\nc{\vzero}{\vec{0}}

\nc{\XO}{{\overset{\circ}{X}}}

\nc{\ya}{{\operatorname{aj}}}

\nc{\nen}{\newenvironment}
\nc{\ol}{\overline}
\nc{\ul}{\underline}

\nc{\Lra}{\Longrightarrow}

\nc{\Llra}{\Longleftrightarrow}
\nc{\hra}{\hookrightarrow}
\nc{\iso}{\overset{\sim}{\lra}}
\nc{\rlh}{\rightleftharpoons}

\nc{\IC}{{\cal{IC}}}
\nc{\PS}{{\cal{PS}}}
\nc{\oCG}{{\overline{\cal G}}}
\nc{\oCQ}{{\overline{\cal Q}}}
\nc{\oCZ}{{\overline{\cal Z}}}
\nc{\dZ}{{\overset{\bullet}{\cal Z}}{}}
\nc{\ddZ}{{\ddot{\cal Z}}{}}
\nc{\oZ}{{\overset{\circ}{\cal Z}}{}}

\nc{\dP}{{\overset{\bullet}{\cal P}}{}}
\nc{\oP}{{\overset{\circ}{\cal P}}{}}
\nc{\oQ}{{\overset{\circ}{\cal Q}}{}}
\nc{\obp}{{\overset{\circ}{{\bf p}}}{}}
\nc{\tbj}{{\tilde{\bf j}}{}}
\nc{\tbp}{{\tilde{\bf p}}{}}
\nc{\tfC}{{\widetilde{\frak C}}{}}
\nc{\tfE}{{\widetilde{\frak E}}{}}
\nc{\tfj}{{\widetilde{\frak j}}{}}
\nc{\tfQ}{{\widetilde{\frak Q}}{}}
\nc{\tfp}{{\widetilde{\frak p}}{}}
\nc{\ofQ}{{\overset{\circ}{{\frak Q}}}{}}
\nc{\tGQ}{{\widetilde{\cal{GQ}}}{}}
\nc{\oGQ}{{\overset{\circ}{\cal{GQ}}}{}}
\nc{\ooGQ}{{\overset{\circ\circ}{\cal{GQ}}}{}}
\nc{\oGZ}{{\overset{\circ}{\cal{GZ}}}{}}
\nc{\tGZ}{{\widetilde{\cal{GZ}}}{}}

\nc{\Ue}{{U_\varepsilon}}
\nc{\Upe}{{\Upsilon_\varepsilon}}
\nc{\crho}{{\check{\rho}}}
\nc{\ctheta}{{\check{\theta}}}

\nc{\pr}{\protect}

\nc{\nn}{{\newline}}
\nc{\np}{{\newpage}}    

\nc{\lab}{      \label}
\nc{\npp}{{     \newpage\setcounter{page}{0}    }}
\nc{\setpart}{{         \setcounter{part}       }}
\nc{\setpage}{{         \setcounter{page}       }}
\nc{\setsection}{{      \setcounter{section}    }}
\nc{\nd}{ $$\text{ This version is preliminary and approximate,
it is not for distribution. }$$ }

\nc{\noi}{{\noindent}}

\nc{\cont}{\tableofcontents}

\nc{\sbr}{{\smallpagebreak}}
\nc{\mbr}{{\medpagebreak}}
\nc{\bbr}{{\bigpagebreak}}
\nc{\bbb}{ \boldsymbol         }
\nc{\bul}{ \bullet         }    

\nc{\bem}{{     \begin{em}      }}
\nc{\eem}{{     \end{em}        }}

\nc{\bbox}{{    \blackbox       }}      

\nc{\ra}{{      \rightarrow     }}
\nc{\laa}{{     \leftarrow      }}      
\nc{\lra}{{\longrightarrow}}
\nc{\lr}{{\leftrightarrow}}             
\nc{\lrs}{{\rightleftarrows}}           

\nc{\imp}{{\Rightarrow}}                
\nc{\eq}{{\Leftrightarrow}}             
        \nc{\Ra}{{\Rightarrow}}                 
        \nc{\LRa}{{\Leftrightarrow}}            

\nc{\inj}{{\pr  \hookrightarrow }}              
\nc{\injj}{{\pr \hookleftarrow  }}              

\nc{\sur}{{     \twoheadrightarrow      }}      
\nc{\surr}{{    \twoheadleftarrow       }}      

\nc{\mm}{{\mapsto}}                             
\nc{\mmm}{{\mapsfrom}}                          

\nc{\va}{{\uparrow}}                            

\nc{\bb}{\pr\underset}           
\rc{\aa}{\pr\overset}            

\nc{\indd}{{ ${} \ \ \ \ \  \ \        {} $     }}      
\nc{\inddd}{{   \indd\indd                      }}      
\nc{\nnd}{{     \nn  \indd                      }}      
\nc{\nndb}{{    \nn  \indd $\bul$               }}      

\nc{\bss}{{\backslash}}                         
\nc{\barr}{     \overline       }               
\nc{\ud}{       \underline      }               

\nc{\ti}{\tilde}              
\nc{\tii}{\widetilde}         

\nc{\hatt}{\widehat}                            
\nc{\hata}{{    \bbb{ \hat{} }          }}      

\nc{\ch}{\check}                                
\nc{\cha}{{     \bbb{ \check{} }        }}      

\nc{\sub}{{     \subseteq       }}         
\nc{\subb}{{    \supseteq       }}         
\nc{\nsub}{{    \nsubseteq      }}         
\nc{\nsubb}{{   \nsupseteq      }}         %

\nc{\nin}{{     \notin  }}

\nc{\lb}{\langle}                                       
\nc{\rb}{\rangle}

\nc{\lB}{       \left(  }                               
\nc{\rB}{       \right) }
\nc{\BBl}{{     \bbb{ \left( \right.}   }}              
\nc{\BBr}{{     \bbb{ \left. \right)}   }}

\nc{\mat} {             \left(          \matrix }
\nc{\emat}{             \endmatrix      \right) }
\nc{\sm} {              \left(          \smallmatrix    }
\nc{\esm}{              \endsmallmatrix \right) }
\nc{\smat} {            \left(          \smallmatrix    }
\nc{\esmat}{            \endsmallmatrix \right) }
\nc{\imat} {            \left.          \matrix }
\nc{\eimat}{            \endmatrix      \right. }
\nc{\ism} {             \left.          \smallmatrix    }
\nc{\eism}{             \endsmallmatrix \right. }

\nc{\ca}{               \left\{         \smallmatrix    }       
\nc{\eca}{              \endsmallmatrix \right\}        }
\nc{\Ca}{               \left\{         \matrix         }
\nc{\eCa}{              \endmatrix      \right\}        }
\nc{\Eca}{              \endmatrix      \right.         }

\nc{\com}{      \begin{diagram} }
\nc{\ecom}{       \end{diagram} }

\nc{\tab}{      \begin{tabular}         }
\nc{\etab}{     \end{tabular}           }       
\nc{\hl}{{      \hline                  }}

\nc{\Eq}{       \begin{equation}        }
\nc{\Eeq}{      \end{equation}  }
\nc{\aln}{      \begin{align}   }
\nc{\ealn}{     \end{align}     }

\nc{\se}{       \section                }
\nc{\sus}{      \subsection             }
\nc{\sss}{      \subsubsection          }

\nc{\Lemm}{     \subsection{Lemma}              }
\nc{\lemm}{     \subsubsection{\bf Lemma}               }
\nc{\slemm}{    \subsubsection*{\bf Lemma}              }

\nc{\sublemm}{  \subsubsection{\bf Sublemma}            }
\nc{\ssublemm}{         \subsubsection*{\bf Sublemma}           }

\nc{\Pro}{      \subsection{Proposition}        }
\nc{\pro}{      \subsubsection{\bf Proposition} }
\nc{\spro}{     \subsubsection*{\bf Proposition}        }

\nc{\Corr}{     \subsection{Corollary}          }
\nc{\cor}{     \subsubsection{\bf Corollary}   }
\nc{\scor}{    \subsubsection*{\bf Corollary}  }

\nc{\Theo}{     \subsection{Theorem}            }
\nc{\theo}{     \subsubsection{\bf Theorem}             }
\nc{\stheo}{    \subsubsection*{\bf Theorem}    }

\nc{\rem}{      \subsubsection{Remark}          }
\nc{\srem}{     \subsubsection*{Remark} }
\nc{\rems}{     \subsubsection{Remarks}         }
\nc{\srems}{    \subsubsection*{Remarks}        }

\nc{\conj}{     \subsubsection{Conjecture}      }
\nc{\sconj}{    \subsubsection*{Conjecture}     }

\nc{\ex}{       \subsubsection{Example}         }
\nc{\sex}{      \subsubsection*{Example}        }
\nc{\exs}{      \subsubsection{Examples}        }
\nc{\sexs}{     \subsubsection*{Examples}       }

\nc{\h}{{       \hslash }}      
\nc{\All}{{     \forall }}
\nc{\yy}{\infty}
\nc{\ys}{{  \frac{\infty}{2}  }}

\nc{\pl}{{\oplus}}                              
\nc{\tim}{{\times}}
\nc{\btim}{{\boxtimes}}
\nc{\ltim}{\ltimes}                     %
\nc{\rtim}{\rtimes}                     %
\nc{\ltr}{\triangleleft}        %
\nc{\rtr}{\triangleright}       %

\nc{\ten}{{     \otimes         }}
\nc{\Lten}{{    \aa{L}\otimes   }}            
\nc{\Ltim}{{    \aa{L}\times    }}            
\nc{\tenA}{     \bb{A}\ten      }
\nc{\tenB}{     \bb{B}\ten      }
\nc{\tenZ}{     \bb{\Z}\ten     }
\nc{\tenR}{     \bb{\R}\ten     }
\nc{\tenC}{     \bb{\C}\ten     }
\nc{\tenk}{     \bb{\k}\ten     }
\nc{\bten}{{\boxtimes}}                         

\nc{\con}{{ @>{\protect\cong}>> }}      
\nc{\conl}{{    @>{\cong}>>     }}      
\nc{\conn}{{    @<{\cong}<<     }}      
\nc{\Con}{{     \equiv          }}      
\nc{\appr}{{    \sim            }}      
\nc{\eqr}{{     \sim            }}      

\nc{\fra}{      \frac   }       

\nc{\ha}{{ \frac{1}{2} }}               
        \nc{\half}{{ \frac{1}{2} }}

\nc{\ci}{{\circ}}               
\nc{\cd }{{\cdot}}              
\nc{\cdd}{{\cdot}}              
                \nc{\cdx}{{\cdot}}                
\nc{\cddd}{{\cdot\cdot\cdot}}   

\nc{\ox}{{      \OO_X           }}               

\nc{\omx}{{     \om_X           }}               
\nc{\Omx}{{     \Om_X^1         }}               

\nc{\cupp}{\bigcup}             
\nc{\capp}{\bigcap}
\nc{\pll}{\bigoplus}

\nc{\pii}{\prod}                
\nc{\ppii}{\bigprod}            
\nc{\cci}{\sqcup}              
\nc{\ccii}{\bigsqcup}

\nc{\wwe}{\bigwedge}            
\nc{\cce}{\bigcoprod}           

\nc{\aaa}{      \stackrel       }

\nc{\edd}{{ \end{document}      }}

\nc{\tx}{       \text           }               

\nc{\df}{{      \overset{\text{def}}=\  }}       
\nc{\dff}{{\    \overset{\text{def}}=\  }}       
\nc{\inv}{{     {}^{-1}                 }}          
\nc{\thh}{      ^{\text{th}}            }           

\nc{\emp}{{   \emptyset}}                       

\nc{\we}{{\wedge}}                              
\nc{\wee}{{     \aa{\bul}\wedge }}              
\nc{\wetwo}{{     \overset{2}\wedge       }}    

\nc{\limp}{{    \underset {\leftarrow}\lim      }}
\nc{\limi}{{    \underset {\rightarrow}\lim     }}
\nc{\plimp}{{\pro\underset {\leftarrow}\lim     }}
\nc{\plimi}{{\pro\underset {\rightarrow}\lim    }}

\nc{\ppp}{{ {\Bbb P}^1 }}                       
\nc{\ppn}{{ {\Bbb P}^n }}                       
\nc{\pt}{       { \text{pt} }   }               

\nc{\qlb}{{ \barr{{\Bbb Q}_l} }}                
\nc{\ffq}{{  {\Bbb F}_q  }}                     
\nc{\ffp}{{  {\Bbb F}_p  }}                     
\nc{\tw}{   {}^{(1)}    }               

\nc{\Spec}{{ \text{Spec}                }}
\nc{\Specf}{{ \text{Specf}                      }}
\nc{\aand}{{\ \ \ \text{and}\ \ \       }}
\nc{\oor}{{\ \  \text{or}\ \    }}
\nc{\hk}{{       \mathrm{ hyperk\ddot{a}hler }    }}

\rc{\Im}{{      \text{Im}       }}
\nc{\rank}{{    \ \text{rank}\  }}
\nc{\Res}{{     \  \text{Res}   }}
\nc{\End}{{     \text{End}      }}
\nc{\RHom}{{    \text{RHom}     }}
\nc{\HHom}{{    \text{$\HH$om}  }}
\nc{\EEnd}{{    \text{$\EE nd$} }}
\nc{\AAut}{{    \text{$\AA ut$} }}
\nc{\RHHom}{{   \text{R$\HH$om} }}
\nc{\Der}{{     \text{Der}      }}

\nc{\ord        }{{ \text{ord} }}                       
\nc{\divv       }{{ \text{div} }}                       
\nc{\Lie        }{{ \text{Lie} }}

\nc{\timA} {{   \underset{A}\tim             }}
\nc{\timB} {{   \underset{B}\tim             }}
\nc{\timC} {{   \underset{C}\tim             }}
\nc{\timG} {{   \underset{G}\tim             }}
\nc{\timH} {{   \underset{H}\tim             }}
\nc{\timN} {{   \underset{N}\tim             }}
\nc{\timP}{{    \underset{P}\tim             }}
\nc{\timQ}{{    \underset{Q}\tim             }}
\nc{\timS} {{   \underset{S}\tim             }}
\nc{\timT} {{   \underset{T}\tim             }}
\nc{\timU} {{   \underset{U}\tim             }}
\nc{\timV} {{   \underset{V}\tim             }}
\nc{\timX} {{   \underset{X}\tim             }}
\nc{\timY} {{   \underset{Y}\tim             }}
\nc{\timZ} {{   \underset{Z}\tim             }}

\nc{\ab}{{       ^{\text{ab}}                   }}
\nc{\af}{{       ^{\text{aff}}                  }}

\nc{\cod}{\text{codim}} 

\rc{\AA}{{\cal A}}
\rc{\BB}{{\cal B}}
\nc{\CC}{{\cal C}}
\nc{\DD}{{\cal D}}
\nc{\EE}{{\cal E}}
\nc{\FF}{{\cal F}}
\nc{\GG}{{\cal G}}
\nc{\HH}{{\cal H}}
\nc{\II}{{\cal I}}
\nc{\JJ}{{\cal J}}
\nc{\KK}{{\cal K}}
\nc{\LL}{{\cal L}}
\nc{\MM}{{\cal M}}
\nc{\NN}{{\cal N}}
\nc{\OO}{{\cal O}}
\nc{\PP}{{\cal P}}
\nc{\QQ}{{\cal Q}}
\nc{\RR}{{\cal R}}
\rc{\SS}{{\cal S}}
\nc{\TT}{{\cal T}}
\nc{\UU}{{\cal U}}
\nc{\VV}{{\cal V}}
\nc{\WW}{{\cal W}}
\nc{\ZZ}{{\cal Z}}
\nc{\XX}{{\cal X}}
\nc{\YY}{{\cal Y}}

\nc{\A}{{\Bbb A }}
\nc{\B}{{\Bbb B}}
\nc{\C}{{\Bbb C}}
                \nc{\cc}{{\Bbb C}}
\nc{\Cs}{{\Bbb C^*}}
                \nc{\cs}{{\Bbb C^*}}
                \nc{\ccs}{{\Bbb C^*}}
\nc{\D}{{\Bbb D}}
\nc{\DDD}{{\Bbb D}}
\nc{\E}{{\Bbb E}}
\nc{\F}{{\Bbb F}}
\nc{\G}{{\Bbb G}}
        \nc{\hH}{{\Bbb H}}
\nc{\I}{{\Bbb I}}
\nc{\J}{{\Bbb J}}
\nc{\K}{{\Bbb K}}
        \nc{\lL}{{\Bbb L}}
\nc{\M}{{\Bbb M}}
\nc{\N}{{\Bbb N}}
        \nc{\oO}{{\Bbb O}}
        \nc{\pP}{{\Bbb P}}
\nc{\Q}{{\Bbb Q}}
\nc{\R}{{\Bbb R}}
        \nc{\sS}{{\Bbb S}}
\nc{\T}{{\Bbb T}}
\nc{\U}{{\Bbb U}}
\nc{\V}{{\Bbb V}}
\nc{\W}{{\Bbb W}}
\nc{\Z}{{\Bbb Z}}
\nc{\X}{{\Bbb X}}
\nc{\Y}{{\Bbb Y}}

\nc{\k}{{\Bbbk}}

\let\L\lL
\let\O\oO

\nc{\fA}{{\frak A}}
\nc{\fB}{{\frak B}}
\nc{\fC}{{\frak C}}
\nc{\fD}{{\frak D}}
\nc{\fE}{{\frak E}}
\nc{\fF}{{\frak F}}
\nc{\fG}{{\frak G}}
\nc{\fH}{{\frak H}}
\nc{\fI}{{\frak I}}
\nc{\fJ}{{\frak J}}
\nc{\fK}{{\frak K}}
\nc{\fL}{{\frak L}}
\nc{\fM}{{\frak M}}
\nc{\fN}{{\frak N}}
\nc{\fO}{{\frak O}}
\nc{\fP}{{\frak P}}
\nc{\fQ}{{\frak Q}}
\nc{\fR}{{\frak R}}
\nc{\fS}{{\frak S}}
\nc{\fT}{{\frak T}}
\nc{\fU}{{\frak U}}
\nc{\fV}{{\frak V}}
\nc{\fW}{{\frak W}}
\nc{\fZ}{{\frak Z}}
\nc{\fX}{{\frak X}}
\nc{\fY}{{\frak Y}}
\nc{\fa}{{\frak a}}
\nc{\fb}{{\frak b}}
\nc{\fc}{{\frak c}}
\nc{\fd}{{\frak d}}
\nc{\fe}{{\frak e}}
\nc{\ff}{{\frak f}}
\nc{\fg}{{\frak g}}
\nc{\fh}{{\frak h}}
\nc{\fiI}{{\frak i}}  
        \nc{\ffi}{{\frak i}}  
\nc{\fj}{{\frak j}}
\nc{\fk}{{\frak k}}
\nc{\fl}{{\frak{l}}}
\nc{\fm}{{\frak m}}
\nc{\fn}{{\frak n}}
\nc{\fo}{{\frak o}}
\nc{\fp}{{\frak p}}
\nc{\fq}{{\frak q}}
\nc{\fr}{{\frak r}}
\nc{\fs}{{\frak s}}
\nc{\ft}{{\frak t}}
\nc{\fu}{{\frak u}}
\nc{\fv}{{\frak v}}
\nc{\fw}{{\frak w}}
\nc{\fz}{{\frak z}}
\nc{\fx}{{\frak x}}
\nc{\fy}{{\frak y}}

\nc{\Aa}{{  \text{A}    }}
\nc{\Bb}{{  \text{B}    }}
\nc{\Cc}{{  \text{C}    }}
\nc{\Dd}{{  \text{D}    }}
\nc{\Ee}{{  \text{E}    }}
\nc{\Ff}{{  \text{F}    }}
\nc{\Gg}{{  \text{G}    }}
\nc{\Hh}{{  \text{H}    }}
\nc{\Ii}{{  \text{I}    }}
\nc{\Jj}{{  \text{J}    }}
\nc{\Kk}{{  \text{K}    }}
\nc{\Ll}{{  \text{L}    }}
\nc{\Mm}{{  \text{M}    }}
\nc{\Nn}{{  \text{N}    }}
\nc{\Oo}{{  \text{O}    }}
\nc{\Pp}{{  \text{P}    }}
\nc{\Qq}{{  \text{Q}    }}
\nc{\Rr}{{  \text{R}    }}
\nc{\Ss}{{  \text{S}    }}
\nc{\Tt}{{  \text{T}    }}
\nc{\Uu}{{  \text{U}    }}
\nc{\Vv}{{  \text{V}    }}
\nc{\Ww}{{  \text{W}    }}
\nc{\Zz}{{  \text{Z}    }}
\nc{\Xx}{{  \text{X}    }}
\nc{\Yy}{{  \text{Y}    }}

\nc{\al}{{\alpha }}
\nc{\be}{{\beta }}
\nc{\ga}{{\gamma }}
\nc{\de}{{\delta }}
\nc{\6}{{\partial }}
\nc{\del}{{\partial }}
\nc{\ep}{{\varepsilon }}
\nc{\vap}{{\epsilon }}

\nc{\ze}{{\zeta }}
\nc{\et}{{\eta }}
\rc{\th}{{\theta }}
\nc{\vth}{{\vartheta }}

\nc{\io}{{\iota }}
\nc{\ka}{{\kappa }}
\nc{\la}{{\lambda }}
\nc{\vrho}{{\varrho}}
\nc{\si}{{\sigma }}
\nc{\ups}{{\upsilon }}
\nc{\vphi}{{\varphi }}
\nc{\om}{{\omega }}

\nc{\Ga}{{\Gamma }}
\nc{\De}{{\Delta }}
\nc{\nab}{{\nabla}}
\nc{\Th}{{\Theta }}
\nc{\La}{{\Lambda }}
\nc{\Si}{{\Sigma }}
\nc{\Ups}{{\Upsilon }}
\nc{\Om}{{\Omega }}

\nc{\toc}{{\tableofcontents}}
\def\square{\hbox{\vrule\vbox{\hrule\phantom{o}\hrule}\vrule}}

\begin{document}


\nc{\ii}{{      i\in I          }}
\nc{\tww}{{     {}^{*,1}        }}
\nc{\zhc}{{ \fZ_{HC} }}
\nc{\pp}{ ^{[p]} }

\nc{\utx}{{     \DD_X           }}
\nc{\utb}{{     \DD_\BB         }}

\nc{\dx}{       \DD_X           }
\nc{\db}{       \DD_\BB         }

\nc{\z}{{       ^{\bbb 0}           }}
\nc{\st}{{    { \barr \bu}          }}

\nc{\St}{{      \SS t                   }}

\nc{\hzc}{{     \hattt{(\chi,0)} }}
\nc{\hrc}{{     \hattt{\chi,-\rho}       }}
\nc{\mr}{{      {-\rho} }}

\nc{\ob}{{      \OO_\BB         }}               
\nc{\uc}{{      u_\chi          }}
\nc{\Fp}{{      {\Bbb F}_p      }}
\nc{\frx}{{     Fr_X            }}

\nc{\Fnew}{     \mbox{Fr}_X     }
\nc{\Frx}{      \Fnew           }

\newtheorem{Thm}{Theorem}
\newtheorem{Cor}{Corollary}
\newtheorem{Lem}{Lemma}
\newtheorem{Prop}{Proposition}
\newtheorem{Sublem}{Sublemma}
\newtheorem{Fact}{Fact}
\newtheorem{Claim}{Claim}
\theoremstyle{remark}
\newtheorem{Rem}{Remark}
\newtheorem{Def}{Definition}
\newtheorem{Ex}{Example}

\newcommand{\imbed}{\hookrightarrow}
\renewcommand{\iso}{{\tii \longrightarrow}}
\newcommand{\isol}{{\tii \longleftarrow}}
\newcommand{\To}{\longrightarrow}

\newcommand{\Lotimes}{\overset{\rm L}{\otimes}}

\newcommand{\oplusl}{\bigoplus\limits}
\newcommand{\cupl}{\bigcup\limits}

\def\square{\hbox{\vrule\vbox{\hrule\phantom{o}\hrule}\vrule}}

\renewcommand{\N}{{\mathcal N}}

\renewcommand{\O}{{\mathcal O}}
\renewcommand{\F}{{\mathcal F}}
\renewcommand{\G}{{\mathcal G}}
\renewcommand{\L}{{\mathcal L}}
\renewcommand{\B}{{\mathcal B}}
\newcommand{\bO}{{\bf O}}

\newcommand{\Ntil}{{    \tilde{\mathcal N}      }}
\newcommand{\Ftil}{{    \tilde F                }}
\newcommand{\Dtil}{{    \tii \DD                 }}

\nc{\Loc}{{    \De            }}
\nc{\LLoc}{{    L\De          }}

\nc{\til}{      \tilde  }

\nc{\hattt}{{               }}
\newcommand{\Dmod}{{            mod^c(\DD)             }}
\newcommand{\Dtmodall}{{        mod^c(\Dtil)        }}        
\newcommand{\Dmodall}{{         mod^c(\Dtil)        }}
\newcommand{\Dtmod}{{           mod^c_{\hattt 0}(\tii\DD)}}

\newcommand{\Umodall}{{         mod^{fg}(U)                 }}
\newcommand{\Umodo}{{           mod^{fg}(U^0)               }}
\newcommand{\Umod}{{            mod^{fg}_{\hattt 0}(\tii U)        }}
\newcommand{\Umodallfg}{{       mod^{fg}(U)                    }}
\newcommand{\Umodofg}{{         mod^{fg}(U^0)          }}
\newcommand{\Umodoallfg}{{      mod^{fg}_{\hattt 0}(\tii U)   }}

\renewcommand{\h}{{\mathfrak h}}
\renewcommand{\t}{{\mathfrak t}}
\newcommand{\n}{{\mathfrak n}}
\newcommand{\g}{{\mathfrak g}}
\renewcommand{\bu}{{            \bullet         }}

\newcommand{\gal}{\check{\ }}

\newcommand{\epf}{\square}

\newcommand{\Zet}{{\mathbb Z}}
\newcommand{\Pn}{{\mathbb P}^n}
\newcommand{\Aone}{{\mathbb A}^1}
\newcommand{\Pone}{{\mathbb P}^1}
\newcommand{\Ql}{{  {\mathbb Q}_l       }}
\newcommand{\Qlb}{{ \bar{\mathbb Q}_l   }}
\newcommand{\Qu}{{\mathbb Q}}
\newcommand{\Ce}{{\mathbb C}}
\newcommand{\bbS}{{\mathbb S}}

\newcommand{\Db}{       \text{D}^{b}    }
\newcommand{\Dmin}{     \text{D}^-      }

\renewcommand{\sur}{\twoheadrightarrow}

\newcommand{\eps}{\epsilon}

\renewcommand{\M}{{\mathcal M}}
\renewcommand{\N}{{\mathcal N}}
\renewcommand{\B}{{\mathcal B}}
\renewcommand{\O}{{\mathcal O}}
\renewcommand{\F}{{\mathcal F}}

\newcommand{\Nt}{\tii{\mathcal N}}
\newcommand{\gtil}{{    \tii{\mathfrak g}^* }}
\newcommand{\gt}{{    \tii{\mathfrak g}}}

\newcommand{\Um}{{      \Db(U^0)                }}

\renewcommand{\b}{{\mathfrak b}}
\newcommand{\fri}{{\mathfrak i}}
\newcommand{\m}{{\mathfrak m}}
\renewcommand{\t}{{\mathfrak t}}

\nc{\diff}{{    \bbb{\bss}      }}

\nc{\Coh}{{ \CC oh  }}

\nc{\bi}{   \begin{itemize}\item        }
\rc{\i}{    \item           }
\nc{\ei}{ \end{itemize} }
\nc{\ben}{  \begin{enumerate}\item      }
\nc{\een}{  \end{enumerate}         }

\let\iso\con

\nc{\aff}{{ \tx{aff}    }}

\nc{\RGa}{{ \Rr\Gamma   }}

\nc{\uHom}{{ \underline{\rm Hom}   }}
 \nc{\uuHom}{{ \bf Hom  }}

 \nc{\PR}{{\bf pr}}
 \nc{\AS}{AS}

\nc{\Bhat}{\hatt \BB}

\nc{\Gm}{{\mathbb G}_m}
\nc{\f}[1]{ \fbox{$ $}\footnote{ \fbox{!}#1 }\fbox{$ $}     }




\title[Localization in characteristic $p$]{
Localization of modules for a semisimple
Lie algebra in prime characteristic
}

\author{
Roman Bezrukavnikov
}
\address{\small
Department of Mathematics, Office 2-284,
 Massachusetts Institute of Technology,
77  Massachusetts ave., Cambridge MA 02139.
}
\email{
bezrukav@@math.mit.edu
}

\author{    Ivan Mirkovi\'c     }
\address{\small
Department of Mathematics and Statistics,
University of Massachusetts,
\   Amherst, MA
01003, USA
}
\email{                mirkovic@@math.umass.edu        }

\author{
Dmitriy Rumynin
}
\address{\small
Mathematics Department, University of Warwick, Coventry,
\
CV4 7AL, England
}
\email{        rumynin@@maths.warwick.ac.uk         }

\thanks{
R.B. was partially supported
by NSF grant DMS-0071967 and  Clay Institute,
D.R.  by EPSRC
and
I.M.  by NSF grants.
}

\begin{abstract}
We show that, on the level of  derived categories,
representations of the Lie algebra of a semisimple algebraic group
over a field of finite characteristic
with a given (generalized) regular central character are the same
as coherent sheaves on the formal neighborhood of the
corresponding (generalized) Springer fiber.

The first step is to observe that the derived functor of global
sections provides an equivalence between the derived category of
$\DD$-modules (with no divided powers) on the flag variety and the
appropriate derived category of modules over the corresponding Lie
algebra. Thus the ``derived'' version of the Beilinson-Bernstein
localization Theorem holds
in sufficiently large positive
characteristic.
Next, one finds that for any smooth variety this algebra of
differential operators is an Azumaya algebra on
the cotangent bundle.
In the case of the flag variety it splits on  Springer
fibers, and this allows us to pass from $\DD$-modules to coherent
sheaves. The argument also generalizes to twisted $\DD$-modules. As
an application we prove Lusztig's conjecture on the number of
irreducible modules  with a fixed central character. We also give
a formula for behavior of  dimension of a module under
translation functors and reprove the Kac-Weisfeiler conjecture.

The sequel to this paper \cite{sing} treats singular
 infinitesimal characters.
\end{abstract}

\maketitle

\begin{flushright}

 {\em To Boris Weisfeiler,\ \ \ \ \ \ \ \ \ \ \ }

               \em   missing since 1985 \\
\end{flushright}


\toc

\setcounter{section}{-1}

\se{\bf Introduction } \sss*{ $\fg$-modules and $\DD$-modules } We
are interested in representations of a Lie algebra $\fg$ of a
(simply connected) semisimple algebraic group $G$ over a field
$\k$ of positive characteristic. In order to relate $\fg$-modules
and  $\DD$-modules on the flag variety $\BB$ we use the sheaf $\db$
of  {\em crystalline} differential operators (i.e. differential
operators without divided powers).

The basic observation is a version of the famous Localization Theorem
\cite{BB}, \cite{BrKa} in positive characteristic.
The center of the enveloping algebra $U(\fg)$
contains the ``Harish-Chandra part''
$\fZ_{HC}\df\ U(\fg)^G$ which is familiar from characteristic zero.
 $U(\fg)$-modules where  $\fZ_{HC}$ acts by the same character as
 on the trivial $\fg$-module $\k$ are modules over the central reduction
$U^0\df\ U(\fg)\ten_{\fZ_{HC}} \k$. Abelian categories of
$U^0$-modules  and of $\DD_\BB$-modules  are quite different.
However,  their bounded derived categories are canonically
equivalent if the characteristic $p$ of the base field $\k$ is
sufficiently large, say, $p> h$ for the Coxeter number $h$.
 More generally, one can identify the bounded
derived category of $U$-modules with a given regular (generalized)
Harish-Chandra
central character with the bounded derived category of the appropriately
twisted $\DD$-modules on $\BB$
(Theorem \ref{Localizationthm}).

\sss*{
$\DD$-modules and coherent
sheaves
}
The sheaf $\dx$ of crystalline
differential operators
on a  smooth variety   $X$ over $\k$
has a non-trivial  center,
canonically identified with
the sheaf of functions on the Frobenius twist  $T^*X\tw$ of
the cotangent bundle (Lemma \ref{io}). Moreover $\DD_X$ is
an Azumaya algebra over $T^*X\tw$ (Theorem \ref{Azumaya}).
More generally, the  sheaves of twisted differential operators are
Azumaya algebras on  twisted cotangent bundles (see
\ref{Torsors}).

When one thinks of the  algebra $U(\fg)$ as the right translation
 invariant sections of $\DD_G$, one recovers the well-known fact that
the center of $U(\fg)$ also has the ``Frobenius part''
$\fZ_{Fr}\cong\OO(\fg^*\tw)$, the functions  on the Frobenius
twist of the dual of the Lie algebra.

For  $\chi\in \fg^*$ let
$\BB_{\chi}\subset \BB$ be a connected component of the variety of
all Borel subalgebras $\fb\subset \fg$ such that
$\chi|_{[\fb,\fb]}=0$, for nilpotent $\chi$ this is the
corresponding Springer fiber.
Thus $\BB_\chi$ is naturally a subvariety of a twisted cotangent bundle of
 $\BB$.
Now, imposing the  (infinitesimal) character $\chi\in\fg^*\tw$ on $U$-modules
corresponds to considering $\DD$-modules (set-theoretically) supported on
$\BB_\chi\tw$.

Our second main observation is that
the Azumaya algebra of twisted differential operators splits on the 
formal neighborhood
of $\BB_\chi$ in the
twisted cotangent bundle. So,  the category of twisted
$\DD$-modules supported on $\BB_{\chi}\tw$ is equivalent to  the
category of coherent sheaves supported on $\BB_{\chi}\tw$ (Theorem
\ref{Achi}). Together with the localization, this provides an
algebro-geometric description of representation theory -- the
derived categories are equivalent for  $U$-modules with a
generalized $\fZ$-character
 and for coherent sheaves on the formal neighborhood
of $\BB_{\chi}\tw$ for the corresponding $\chi$.

\sss*{
Representations
}
One   representation theoretic consequence of the passage to algebraic
geometry is the
count of irreducible $U_\chi$-modules with a given regular
Harish-Chandra central character
(Theorem \ref{count}).
This was known previously when $\chi$ is regular nilpotent
in a Levi factor (\cite{FP}), and the general case was conjectured by
Lusztig (\cite{Lu1},\cite{Lu}).
In particular, we find  a
canonical isomorphism
of Grothendieck groups
of $U^0_\chi$-modules and of coherent sheaves on the Springer fiber
$
\BB_{\chi}
$.
Moreover,
 the rank of
this $K$-group  is the same as   the dimension of cohomology of
the corresponding Springer fiber in characteristic zero (Theorem
\ref{K theorem}), hence it is well understood. One of the purposes
of this paper is to provide an approach to Lusztig's elaborate
conjectural description of representation theory of $\fg$.

\sss{
}
Sections 1 and 2  deal with algebras of differential operators
$\utx$.
Equivalence $\Db(mod^{fg}(U^0))\con\ \Db(mod^c(\utb))$ and its
generalizations are proved in
$\S 3$.  In  $\S 4$ we specialize the equivalence to
objects with the $\chi$-action of the Frobenius center $\fZ_{Fr}$.
In  $\S 5$ we relate $\DD$-modules with the $\chi$-action of
$\fZ_{Fr}$ to $\OO$-modules on the Springer fiber $\BB_{\chi}$.
This leads to a dimension formula for $\fg$-modules in
terms of the corresponding coherent sheaves in $\S 6$,
here we also spell out compatibility of our functors  with translation
functors.
Finally, in $\S 7$ we calculate  the rank of
the $K$-group of the Springer fiber, and thus of the corresponding category
of $\fg$-modules.

\sss{} The origin of this study was a suggestion of James Humphreys
that the representation theory of $U^0_\chi $ should be related to
geometry of the Springer fiber $\BB_{\chi}$. This was later
supported by  the work of Lusztig \cite{Lu} and Jantzen
\cite{Ja1}, and by \cite{MR}.

\sss{} We would like to thank Vladimir Drinfeld, Michael
Finkelberg, James Humphreys, Jens Jantzen, Masaharu Kaneda, Dmitry
Kaledin,  Victor Ostrik,  Cornelius Pillen,
Simon Riche and Vadim Vologodsky 
for various
information over the years; special thanks go to Andrea Maffei
for pointing out a mistake in example \ref{exaSL3}(2) in the previous
draft of the paper.
A part of the work was accomplished
while R.B. and I.M. visited the Institute for Advanced Study
(Princeton), and the Mathematical Research Institute (Berkeley);
in addition to excellent working conditions these opportunities
for collaboration were  essential. R.B. is also grateful to the
Independent Moscow University where part of this work was done.

\sss{
Notation
}
We consider schemes over
an algebraically closed  field
$\k$
of characteristic $p>0$.
For an affine $S$-scheme $X \stackrel{q}{\ra} S$, we denote
$
q_*\OO_X
$
by $\OO_{X/S}$, or simply by $\OO_{X}$.
For a subscheme $\fY$ of  $\fX$
the  formal neighborhood
$FN_\fX(\fY)$
is an ind-scheme (a formal scheme),
the notation for the categories of modules
on $\fX$ supported on $\fY$
is introduced in
\ref{Derived categories of sheaves supported on a subscheme},
\ref{Categories} and
\ref{Categories with chi}.
The   Frobenius neighborhood
$FrN_\fX(\fY)$
is introduced in
\ref{Frobenius neighborhoods}.
The inverse image of sheaves  is denoted $f\inv$ and for
$\OO$-modules $f^*$ (both  direct images are denoted $f_*$).
We denote by $\TT_X$
and $\TT^*_X$
the
sheaves of sections of the (co)tangent bundles
$TX$ and $T^*X$.

\se{\bf
Central reductions
of the envelope $\utx $
of the tangent sheaf
}
We will
describe the center of  differential operators
(without divided powers)
as  functions on the Frobenius twist of the cotangent bundle.
Most of the material in this section is standard.

\sus{
Frobenius twist
}

\sss{ Frobenius twist of a $\k$-scheme } Let $X$ be a scheme over
an algebraically closed  field $\k$ of characteristic $p>0$. The
Frobenius map of schemes $X \ra X $ is defined as identity on
topological spaces, but the pull-back of functions is  the $p$-th
power:\ $\Frx^*(f)= f^p$ for $f\in \OO_{X\tw}= \OO_X$. The
Frobenius twist $X\tw$ of $X$ is the $\k$-scheme that coincides
with $X$ as  a scheme (i.e. $X\tw=X$ as a topological space and
$\OO_{ X\tw } =\OO_X$ as a sheaf of rings), but with a different
$\k$-structure: $a\bb{(1)}\cdot f \ \df \ a^{ 1/p } \cdot f,\
a\in\k,\ f\in\OO_{X\tw} $. It makes Frobenius map into  a map of
$\k$-schemes $X \stackrel{\mbox{\tiny Fr}_{X}}{\lra} X^{(1)} $. We
will  use the twists to keep track of using Frobenius maps. Since
$\Frx$ is a  bijection on $\k$-points,  we will often identify
$\k$-points of $X$ and $X\tw$. Also, since $\Frx$ is affine, we
may identify sheaves on $X$ with their $(\Frx)_*$-images. For
instance, if  $X$ is reduced  the $p$-th power map $ \OO_{X\tw}
\ra (\Frx)_*\OO_X$ is injective, and we think of $ \OO_{X\tw} $ as
a subsheaf $\OO_X^p\df \{f^p,\ f\in\ox\}$ of $\OO_X$.

\sss{
Frobenius neighborhoods
}
\lab{Frobenius neighborhoods}
The Frobenius neighborhood
of a subscheme $Y$ of $X$
is the subscheme $
(\Fnew)\inv Y\tw\sub X$,
we denote it $FrN_X(Y)$ or simply
$\ud X_Y$. It  contains $Y$
and
$
\OO_{\ud X_Y}
=\,
\OO_{X}
\bb{\OO_{X\tw}}
\otimes
\OO_{ Y\tw}
\
=\
\OO_{X}
\bb{\OO_{X}^p}
\otimes
\OO_{X}^p/\II_Y^p
\
=\
\OO_{X}
/
\II_Y^p\cdot \OO_X
$
for
the ideal of definition $\II_Y\subseteq\OO_X$
 of $Y$.

\sss{
Vector spaces
}
For a $\k$-vector space $V$ the
$\k$-scheme $V\tw$ has a natural
structure of a vector space over $\k$; the $\k$-linear structure is again
given by
$a\bb{(1)}\cdot
v
\
\df
\
a^{ 1/p }
v,
\
a\in \k,\ v\in V$.
We say that a map $\beta:V\ra W$
between $\k$-vector spaces
is $p$-linear
if it is additive and
$
\beta (a\cdot v)=
a^p\cdot\beta (v)
$;\
this
is the same as a linear map
$V^{(1)}\ra
W$.
The canonical isomorphism of vector spaces
$
(V^*)\tw
\con
(V\tw)^*
$
is given
by
$\alpha\ra \alpha^p$ for
$\alpha^p(v)\df\alpha(v)^p$
(here, $V^*\tw=V^*$ as a set and
$(V\tw)^*$ consists of all $p$-linear $\beta:\ V\ra \k$).
For a smooth $X$,
 canonical $\k$-isomorphisms
$T^*(X\tw)=(T^*X)\tw$
and
$(T(X))\tw\con
T(X\tw)$ are obtained from definitions.

\subsection{
The ring of  ``crystalline'' differential operators $\utx $ }
\lab{utx} Assume that $X$ is a smooth variety.
Below we will occasionally compute in local coordinates:
since $X$ is  smooth,
any point $a$  has a Zariski neighborhood $U$ with
etale coordinates\ $x_1,...,x_n$, i.e., $(x_i)$ define an etale map from
$U$ to ${\mathbb A}^n$ sending $a$ to 0. Then
$dx_i$ form a frame of $T^*X$ at $a$;
the dual frame $\del_1,...,\del_n$ of $\TT_X$
is characterized by $\del_i
(x_j)=\delta_{ij}$.

 Let  $\utx =U_{\OO_X}(\TT_X)$
denote
the   enveloping
algebra of the tangent Lie algebroid $\TT_X$; we call $\utx$ the
sheaf of crystalline differential operators. Thus $\utx$
 is generated by the algebra of functions $\ox$
and the $\ox$-module of vector fields $\TT_X$,
subject to the module and commutator relations
$f\cd \del= f\del$, \ $\del\cd f-f\cd \del= \del(f),\ \del\in\TT_X,\ f\in\ox$,
and the
Lie algebroid relations $\del'\cd \del''-\del''\cd \del'=
[\del',\del''],\ \del',\del''\in\TT_X$. In terms of a local frame
$\del_i$ of vector fields we have \ $\utx =\ \bb{I}\pl\ \OO_X\cdd\del^I$.
One readily checks that $\utx$ coincides with the object defined (in a more
general situation) in \cite{BO}, \S 4, and called there
 ``PD differential operators''.

By the definition of an enveloping algebra,
a sheaf of $\utx$ modules is just an $\OO_X$ module equipped with a flat
connection. In particular, the standard flat connection on the
structure sheaf $\OO_X$  extends to a $\utx $-action. This action
is not faithful: it provides a map from $\utx$
to the ``true'' differential operators ${\bf\D}_X\sub\
\EEnd_\k(\OO_X)$ which contain divided powers of vector fields;
the image of this map is an $\OO_X$-module of finite rank $p^{\dim X}$,
see \cite{BO} or \ref{Lagrangian} below.

For $f\in \OO_X$ the $p$-th power $f^p$ is killed by the
action of $\TT_X$, hence
for any closed
subscheme $Y\sub\ X$ we get an action of
$\utx$ on the structure sheaf $\OO_{\ud
X_Y}$ of the Frobenius neighborhood.

Being defined as an enveloping algebra of a Lie algebroid,
the sheaf of rings $\DD_X$ carries a natural ``Poincare-Birkhoff-Witt''
filtration 
 $\DD_X=\cup \DD_{X,\leq n}$, where $\DD_{X,n+1}= \DD_{X,\leq n}+
\TT_X \cdot \DD_{X,\leq n}$,
$\DD_{X,\leq 0} =\OO_X$. In the 
 following Lemma parts (a,b) 
are proved similarly to the
familiar statements in characteristic zero, while
(c) can be proved by a straightforward use of local coordinates.

\lemm
\lab{grlep}

a) We have a canonical isomorphism of the sheaves of algebras:

$gr(\DD_X)\cong \OO_{T^*X}$.

b)  $\OO_{T^*X}$ carries a Poisson algebra structure, given
by $\{ f_1,f_2\}=[\tilde f_1,\tilde f_2 ] \mod \DD_{X, \leq n_1+n_2-2}$,
$\tilde f_i\in \DD_{X,\leq n_i}$, $f_i=\tilde f_i \mod
 \DD_{X,\leq n_i-1}\in \OO_{T^*X}$, $i=1,2$.

This Poisson structure coincides with the one arising from the
standard symplectic form on $T^*X$.

c) The action of $\DD_X$ on $\OO_X$ induces an injective morphism
$\DD_{X,\leq p-1} \imbed \EEnd (\OO_X)$.

\medskip

We will use the familiar terminology, referring to the image of
$d\in \DD_{X,\leq i}$ in $ \DD_{X,\leq i}/\DD_{X,\leq i-1}\subset
\OO_{T^*X}$ as its symbol.

\subsection{
The difference $\io$ of $p\thh$ power maps on vector fields }

For
any vector field $\del\in\TT_X$, $\del^p\in \utx $ acts on
functions as another vector field which one denotes
$\del\pp\in\TT_X$. For $\del\in\TT_X$ set 
 $ \io(\del) \df \del^p-\del\pp\in \utx $. The map $\io$  lands
in the kernel of the action on $\OO_X$;
it   is injective, since it is injective on symbols.

\lemm
\lab{io1}

a) The map 
$\io:
\TT_X\tw\ra
\utx
$ is $\OO_{X\tw}$-linear, i.e., $\io(\del)+\io(\del')=\io(\del+\del')$
and 
$
\io(f\del)=\
f^p\cdd \io(\del)
,\ \del,\del'\in\TT_{X\tw},\
f\in\OO_{X\tw}$.

b) The image of  $\io$ is contained in the center of $\DD_X$.

\pf\footnote{Another proof of the lemma follows directly from
Hochschild's identity (see \cite{Ho}, Lemma 1).}
For each of the two identities in (a), 
both sides act by zero on $\OO_X$. Also, they lie in $D_{X,\leq p}$,
 and clearly coincide modulo $D_{X,\leq p-1}$.
 Thus the identities follow from Lemma \ref{grlep}(c).

b) amounts to: $[f,\io(\del)]=0$, $[\del',\io(\del)]=0$,
for $f\in \OO_X$, $\del,\del'\in \TT_X$. In both cases the left hand sides
lie in $\DD_{X,\leq p-1}$: this is obvious in the first case, and in 
the second one it
follows from the fact that the $p$-th power of an element in a Poisson
algebra in characteristic $p$ lies in the Poisson center.
The identities follow, since the left hand sides kill $\OO_X$.
 \epf

\medskip

Since $\io$ is $p$-linear,
we consider it as a linear map $\io:
\TT_X\tw\ra
\utx
$.

\lemm \lab{io} The map $ \io: \TT_X\tw\ra\ \utx $ extends to an
isomorphism of $ \fZ_X\df\ 
\OO_{T^*X\tw/X\tw}$ 
and the
center $ Z(\utx ) $. In particular, $Z(\utx ) $ contains
$
\OO_{X^{(1)} }$.

\pf
For $f\in \OO_X$ we have $f^p\in Z(\DD_X)$, because the identity
$ad(a)^p=ad(a^p)$ holds in an associative ring in
 characteristic $p$, which shows that
$[f^p,\del]=0$ for $\del\in \TT_X$. This, together with Lemma
\ref{io1}, yields a homomorphism $\fZ_X\to Z(\DD_X)$.
This homomorphism is injective, because the induced map on symbols
is the Frobenius map $\varphi\mapsto \varphi^p$, $\fZ=\OO_{T^*X\tw}
\to \OO_{T^*X}$. To see that it is surjective it suffices to see
that 
the Poisson center of the sheaf of Poisson algebras $\OO_{T^*X}$
is spanned by the $p$-th powers. Since the Poisson structure
arises from a non-degenerate two-form, a function $\varphi
\in \OO_{T^*X}$ lies in the Poisson center if and only if
$d\varphi=0$. It is a standard fact that a function $\varphi$ on
a smooth variety over a  perfect field of characteristic $p$
satisfies 
$d\varphi =0$ if and only if $\varphi=\eta^p$ for some $\eta$.
\epf

\sex If $X={\mathbb A}^n$, so $\DD_X=\k\langle x_i, \del_i\rangle$
is the Weyl algebra, then $Z(\DD_X)=\k[x_i^p,\del_i^p]$.

\sss{ The Frobenius center of enveloping algebras }
\lab{Zp}
Let  $G$ be an algebraic group over $\k$, and
$\fg$ be  its Lie algebra.
Then $\fg$ is the algebra of left invariant vector fields on $G$,
and the $p$-th power map on vector fields induces the structure of a
restricted Lie algebra on $\fg$. Considering left invariant sections
of the sheaves in Lemma \ref{io} we get an embedding
$\OO(\fg^*\tw)\aaa{\io_\fg}\inj Z(U(\fg))$;
we have  $\io_\fg(x)=x^p-x^{[p]}$ for $x\in \fg$.
 Its image is denoted
$\fZ_{Fr}$ (the ``Frobenius part'' of the center).

 From the
construction of $\fZ_{Fr}$ we see that if $G$ acts on a smooth
variety $X$ then $\fg\ra\ \Ga(X,\TT_X)$ extends to  $U(\fg)\ra\
\Ga(X,\DD_X)$ and the constant sheaf $(\fZ_{Fr})_X=
\OO(\fg^*\tw)_X$ is mapped into the center $\fZ_X= \OO_{T^*X\tw}$.
The last map comes from the moment map $T^*X\ra\ \fg^*$.

$U\fg$ is a vector bundle
of rank $p^{\dim(\fg)}$ over $\fg^*\tw$.
Any $\chi\in\fg^*$ defines a point $\chi$ of $\fg^*\tw$ and
a central reduction $U_\chi(\fg)\df\ U(\fg)\ten_{\fZ_{Fr}}\k_\chi$.

\subsection{
Central reductions } For any closed subscheme $\YY\sub T^*X $ one
can restrict $\utx $ to $\YY\tw\sub T^*X\tw$, we denote the
restriction $ \DD_{X,\YY} \ \df\ \utx \bb{ \OO_{T^*X\tw/X\tw } }
\ten \OO_{\YY\tw/X\tw } $.

\sss{
Restriction to the Frobenius neighborhood of a subscheme of $X$
}
\lab{omY}
A closed subscheme $Y\inj X$ gives a subscheme
$T^*X|Y\sub T^*X$, and the corresponding central reduction
$$
\utx
\bb{\OO_{T^*X\tw }\
}
\ten
\OO_{(T^*X|Y)^{(1)} }
\
=
\
\utx
\bb{\OO_{X\tw }\
}
\ten
\OO_{Y^{(1)} }
\
=
\
\utx
\bb{\OO_{X }\
}
\ten
\OO_{\ud X_Y }
,$$
is just the restriction of
$\utx $ to the Frobenius neighborhood of $Y$.
Alternatively, this is
the enveloping algebra
of the restriction  $\TT_X|\ud X_Y$ of the
Lie algebroid $\TT_X$.
Locally,
it is of the form
$
\
\bb{I}\pl
\
\OO_{\ud X_Y }
\del^I
$.
As a quotient of $\utx$ it is obtained by imposing
$f^p=0$ for $f\in\II_Y$.
One can say that the reason we can restrict Lie algebroid
$\TT_X$ to
the
Frobenius neighborhood
$\ud X_Y$
is that
for vector fields
(hence also for $\utx$),
the subscheme $\ud X_Y$ behaves as an open
subvariety of $X$.

Any  section $\om$ of $T^*X$ over $Y\sub X$ gives $\om(Y)\sub\
T^*X|Y$, and a  further reduction $ \DD_{X,\om(Y)} $. The
restriction to $\om(Y)\sub\ T^*X|Y$ imposes $\io(\del)=
\lb\om,\del\rb^p$, i.e., $\del^p=\del\pp+ \lb\om,\del\rb^p,\
\del\in\TT_X$. So, locally, $ \DD_{X,\om(Y)} = \
\bb{I\in\{0,1,...,p-1\}^n}\pl \ \OO_{\ud X_Y} \del^I $ \ and $ \
\del_i^p = \del_i\pp + \lb\om, \del_i \rb^p = \lb\om, \del_i \rb^p
.$

\sss{
The  ``small''
differential operators
$\DD_{X,0}$}
When $\YY$ is the zero section of $T^*X$
(i.e.,  $X=Y$ and $\om=0$),
we get the algebra
$\DD_{X,0}$
by imposing
in $\utx$ the relation $\io\del=0$,
i.e.,
$\del^p=
\del\pp
,\
\del\in\TT_X$ (in local coordinates ${\del_i}^p=0$).
The action of $\utx$ on $\ox$ factors through  $\DD_{X,0}$
since  $\del^p$ and $\del\pp$ act the same
on $\OO_X$.
Actually,
$\DD_{X,0}$ is the image of the canonical map
$\utx\ra{\bbb\D}_X$ from \ref{utx} (see \ref{Lagrangian}).

\se{\bf
Azumaya property of
$\utx $
}

\sus{
Commutative subalgebra
$\AA_X\sub\utx$
}\lab{21}
We will denote the centralizer of $\OO_X$ in $\utx $ by
$\
\AA_X
\
\df\
Z_{\utx }
(\OO_X)
$,\ and the pull-back of $T^*X\tw$ to $X$ by
\
$
\
T^{*,1}X
\
\df\
X
\tim_{
 X\tw
}
T^*X\tw
$.

\lemm
$
\AA_X= \OO_X\cdd\fZ_X
=\
\OO_{
T^{*,1}X
/ X
}
$.

\pf
The problem is local so  assume that $X$ has  coordinates
$x_i$. Then
$\utx
=
\pl\ \OO_X\del^I
$
and
$\fZ_X=\pl\ \OO_{X\tw}\del^{pI}$
(recall that
$\io(\del_i)
=
{\del_i}^p
$). So,
$
\OO_X\cdd\fZ_X
\
=\
\pl\ \OO_X\del^{pI}\
\conn\
\OO_X
\ten_{
\OO_{X\tw}
}
\fZ_X
$,
and this is the algebra
$\OO_X
\ten_{\OO_{X\tw}}
\OO_{
T^*X\tw
}
$
of functions on
$T\tww X
$.
Clearly,
$
Z_{\utx}(\OO_X)
$
contains
$\OO_X\cd\fZ_X$, and the converse
$
Z_{\utx}(\OO_X)
\sub\
\pl\ \OO_X\del^{pI}$
was already observed in the proof of
Lemma \ref{io}.
\epf

\rem\lab{p_curv} In view of the lemma, any $\utx$-module $\EE$
carries an action of $\OO_{T^{*,1} X}$; such an action is the same
as  a section $\omega$ of $\mbox{Fr}^*(\Omega^1_X)\otimes \EEnd_{\OO_X} (\EE)$.
As noted above $\EE$ can be thought of as an $\OO_X$ module with a
flat connection;
the section $\omega$ is known as {\em the $p$-curvature}
of this connection.
The section $\omega$ is parallel for the
induced flat connection on $\mbox{Fr}^*(\Omega^1_X)\otimes \EEnd_{\OO_X} (\EE)$.

\sus{
Point modules $\de^\ze$
}
\lab{Point modules}
A cotangent vector  $\ze= (b,\om)\in T^*X\tw$
(i.e.,
$b\in X\tw$ and $\om\in
T^*_a X\tw$),
defines a central reduction
$\DD_{X,\ze}= \utx\ten_{\fZ_X}\ \OO_{\ze\tw}$. 
Given a lifting $a\in T^*X$ of $b$ under the Frobenius map
(such a lifting exists since $\k$ is perfect and it is always unique),
we get  a $\utx$-module
$
\de^\xi
\df
\
\utx \ten_{\AA_X}\OO_{\xi}$, 
where we have set  $\xi=(a,\om)\in T^{*,(1)}X$. It is a central
reduction of the $\utx$-module $\de_a\df
\
\utx \ten_{\ox}\OO_a
$ of distributions at $a$, namely
$\de^\xi =
\de_a \ten_{\fZ_X}\OO_{\ze}
$.
In local coordinates
at $a$,
\ref{omY} says that
$
\DD_{X,\ze}
$ has a $\k$-basis
$
x^J
\del^I,\
{
I,J
\in\{0,1,...,p-1\}^n}
$
with
$
x_i^p=0
$
and
$
\del_i^p
=
\lb\om,\del_i\rb^p
$.

\lemm
\lab{point modules}
Central reductions of $\DD_X$ to points of $T^*X\tw$ are matrix algebras.
More precisely, in the above notations we have:
$$
\Ga(X,\DD_{X,\ze})
\
\con
\
\End_\k (\Ga(X,\de^\xi))
.$$

\pf
Let $x_1,...,x_n$ be local coordinates at $a$.
Near $a$,
$
\utx
=
\pl_{I\in\{0,...,p-1\}^n}\
\del^I\cdd\AA_X
$, hence
$\de^\xi
\cong\
\pl_{I\in\{0,...,p-1\}^n}\
\k\del^I
$.
Since $x_i(a)=0$,
$$
x_k\cdd\del^I
\
=
\
I_k\cdd\del^{I-e_k}
\
\
\
\aand\
\
\
\del_k\cdd\del^I
=
\
\Ca
\del^{I+e_k}
&
\text{
if $I_k+1<p$,
}
\\
\om(\del_i)^p\cdd\del^{I-(p-1)e_k}
&
\text{
if $I_k=p-1$.
}
\eCa
.$$
Irreducibility of $\de^\xi$ is now standard
-- $x_i$'s act on polynomials in
$\del_i$'s by derivations, so for
$0\ne
P=
\sum_{I\in\{0,...,p-1\}^n}\
\
c_I\del^I\ \in\ \de^\xi$
and a maximal $K$ with $c_K\ne 0$,
$x^K\cdd P$ is a non-zero scalar. Now multiply with $\del^I$'s to get
all of $\de^\xi$.
So
$\de^\xi$ is an irreducible
$\DD_{X,\ze}$-module.
Since
$
\dim
\DD_{X,\ze}
=p^{2\dim(X)}=
(\dim
\de^\xi
)^2
$ we are done.
\epf

Since the lifting $\xi\in T^{*,(1)}X$ of a point $\ze\in T^*X\tw$
exists and is unique, we will occasionally talk about point modules
associated to a point in $T^*X\tw$, and denote it by $\de^\ze$,
$\ze\in T^*X\tw$.

\pro
\lab{Splitting}
({\em Splitting of $\DD_X$ on $T\tww X$.})
Consider
$\utx
$
as an $\AA_X$-module
$
(\utx)_{\AA_X}
$
via the right multiplication.
The left multiplication by $\utx$ and the right multiplication by
$\AA_X$ give an isomorphism
$$
\utx
\bb{
\fZ_X
}\ten
\AA_X
\
\con
\
\EEnd_{
\AA_X
}
(
(\utx)_{\AA_X}
)
.$$

\pf
Both sides are vector bundles over $T\tww X=\Spec (\AA_X)$:\
the $\AA_X$-module
$
(\utx)_{\AA_X}
$
has a local frame
$\del^I,\
I\in
\{0,...,p-1\}^{\dim X}
$;
while
$x^J\del^I,\
J,I\in
\{0,...,p-1\}^{\dim X}
$
is a local frame for both the $\fZ_X$-module $\utx$ and
the $\AA_X$-module $\utx\ten_{\fZ_X} \AA_X$.
So,  it suffices to check that
the map is an isomorphism on fibers.
However, this is the claim of Lemma
\ref{point modules},
since
the restriction of the map to a  $\k$-point
$\ze$ of $T\tww X$
is the action of
$
(\utx
\ten_{\fZ_X}
\AA_X
)
\ten_{\AA_X}
\OO_\ze
\
=\
\utx
\ten_{\fZ_X}
\OO_\ze
\
=\
\DD_{X,\ze}
$
on
$
(\utx)_{\AA_X}
\ten_{\AA_X}
\OO_\ze
\
=\
\de^\ze
$.
\epf

\theo
\lab{Azumaya}
$\utx$ is an
Azumaya algebra
over $T^* X^{(1)}$ (non-trivial if  $\dim(X)>0$).

\pf
One of the  characterizations of Azumaya algebras
is that they are coherent as $\OO$-modules
and become matrix algebras
on a flat cover \cite{MI}.
The map $T^{\ast,1}X\ra T^\ast X\tw$
is faithfully flat, i.e., it is a flat cover, since
the Frobenius map $X\ra X\tw$ is flat for smooth $X$
(it is surjective and
on the  formal neighborhood of a point it
is given by $\k[[x_i^p]]\inj \k[[x_i]]$).
If $\dim(X)>0$, then
$\DD_X$ is non-trivial, i.e. it is
 not isomorphic to an algebra of the form $End(V)$ for a vector
bundle $V$, because  locally in the Zariski topology of $X$, 
$\DD_X$  has no zero-divisors, since $gr(\DD_X)= \OO_{T^*X}$;
while the algebra of endomorphisms of a vector bundle of rank
higher than one on an affine algebraic variety has zero divisors.
\epf

\rems (1) A related Azumaya algebra was considered in \cite{Hur}.

(2) One can give a different, somewhat shorter proof of Theorem 
\ref{Azumaya} based on the fact that a function on 
 a smooth $\k$-variety has zero
differential if and only if  it is a $p$-th power, which implies that
any Poisson ideal in $\O_{T^*X}$ is induced from $\O_{T^*X\tw}$.
This proof applies to a more general situation of the so called
Frobenius constant quantizations of  symplectic varieties in positive
characteristic, see \cite{BezruKa}, Proposition 3.8.

(3) The statement of the theorem can be
compared to the well-known fact that the algebra of differential operators
in characteristic zero is simple: in characteristic $p$ it becomes simple
after a central reduction. Another analogy is with the
 classical Stone -- von Neumann Theorem, which
asserts that $L^2(\R^n)$ is the only irreducible unitary representation of
 the Weyl algebra:  Theorem \ref{Azumaya} implies, in particular,
that the
standard quantization of functions on the Frobenius
 neighborhood of zero in $\A^{2n}_\k$ has unique irreducible representation
realized in the space of functions
 on the Frobenius neighborhood of zero in $\A^n_\k$.

\sss{
Splitting on the zero section
}
\lab{Lagrangian}
By a well known observation\footnote{The second author 
thanks  Paul Smith from whom he has learned this observation.} 
the small differential operators, i.e., the restriction $\DD_{X,0}$
of $\DD_X$ to $X\tw\sub\ T^*X\tw$, form a sheaf of matrix algebras. In
the notations above this is the observation that
the action map $(\Frx)_*\DD_{X,0}\con\ \EEnd_{\OO_{X\tw}}((\Frx)_*\ox)$
is an isomorphism by
\ref{point modules}.
Thus Azumaya algebra
$\DD_X$ splits on $X\tw$, and $(\Frx)_*\ox$ is a splitting bundle.
The corresponding equivalence between $Coh_{X\tw}$ and  $\DD_{X,0}$ modules
sends $\FF\in Coh_{X\tw}$ to the sheaf $\Frx^*\FF$ equipped with a standard
flat connection (the one for which  pull-back of a section of $\FF$
is parallel).

\rem\lab{spl_on_Z} Let  $Z\subset T^*X\tw$ be a closed subscheme, such that
the  Azumaya algebra
$\DD_X$ splits on $Z$ (see section  \ref{Splitting on  Springer fibers}
below for more examples of this situation); thus we have
 a splitting vector bundle $\EE_Z$
on $Z$ such that $\DD_X|_Z\iso End(\EE_Z)$. It is easy to see then that
$\EE_Z$ is a locally free rank one module over
$\AA_X|_Z$, thus it can be
thought of as a line bundle on the preimage $Z'$ of $Z$ in $T^*\tw X$
under the map $Fr\times id: X\times _{X\tw} T^*X\tw\to T^*X\tw$.
In the particular case when $Z$ maps isomorphically to its image
$\bar Z$ in $X$
the scheme $Z'$ is identified with the Frobenius neighborhood of $\bar Z$
in $X$. The action of $\DD_X$ equips the resulting line bundle on
$FrN(\bar Z)$ with a flat connection.
 E.g.
the above splitting on the zero-section
corresponds to the trivial line
bundle $\O_X$ with the standard flat connection.

\sus{
Torsors
}
\lab{Torsors}
A torsor  $\tii X@>\pi>> X$ for a torus $T$
defines a Lie algebroid
$\tii\TT_X\df\ \pi_*(\TT_{\tii X})^T$
with the enveloping algebra
$\tii\DD_X\df\ \pi_*(\DD_{\tii X})^T
$.
Let $\ft$ be the Lie algebra of $T$.
Locally, any trivialization of the torsor
splits the exact sequence
$0\ra \ft\ten\OO_X\ra\ \tii\TT_X\ra\TT_X\ra\ 0$
and gives
$\tii\DD_X\ \cong \DD\ten \ U\ft$.
So the map of  the constant sheaf
$U(\t)_X$ into $\tii\DD_X$,
given by the $T$-action, is a central embedding
and $\tii\DD_X$ is a deformation of $\DD_X
\cong \tii\DD_X \otimes_{S(\t)} \k_0
$ over $\ft^*$.
The center $\OO_{T^*\tii X\tw}$
of $\DD_{\tii X}$ gives  a central subalgebra
$(\pi_*\OO_{T^*{\tii X}\tw})^T
=\
\OO_{\tii T^* X\tw}
$ of $\tii\DD_X$. We combine the two into a map from
functions on
$
\tii T^*X\tw\tim_{\ft^*\tw}\ft^*
$
to $Z(\tii\DD_X)$
(the map $\ft^*\to \ft^*\tw$ is the Artin-Schreier map $\AS$;
the corresponding map on the rings of functions
$S(\ft\tw)\to S(\ft)$
is given by $\io(h)= h^p-h\pp,\ h\in\ft\tw$).
Local trivializations again show that this is an isomorphism and that
 $\tii\DD_X$ is an Azumaya
algebra on
$
\tii T^*X\tw\tim_{\ft^*\tw}\ft^*
$,
which splits on
$X\tim_{X\tw}
(\tii T^*X\tw\tim_{\ft^*\tw}\ft^*)$.

In particular, for any $\la\in\ft^*$, specialization
$\DD_X^\la\df\ \tii\DD_X\ten_{S(\ft)}\k_\la$ is an Azumaya algebra
on the twisted cotangent bundle
$T^*_{AS(\la)}X\tw
\df\
\tii T^*X\tw\tim_{\ft^*\tw}AS(\la)
$, which splits on
$
T^{*,(1)}_{AS(\la)}X
\df\
X\tim_{X\tw}T^*_{AS(\la)}X\tw
$.
For instance if $\la=d(\chi)$ is the differential of a character
$\chi$ of   $T$ then  $AS(\la)=0$; thus $T^*_{AS(\la)}X=T^*X$.
In this case  $\DD_X^\la$ is identified with the sheaf
$^{\OO_\chi}\DD_{X}\cong \OO_\chi\ten\DD_X\ten\OO_\chi\inv$ of differential
operators on sections of the line bundle $\OO_\chi$
on $X$, associated to $\tii X$ and  $\chi$.


By a  straightforward generalization of \ref{21}, \ref{Point modules},
$\tii \AA_X\dff \OO_{X\times _{X\tw}\tii T^*X\tw\times_{\ft ^*\tw}
\ft^*}$
embeds into $\tii \DD_X$.
As in \ref{Point modules}, for a
point $\zeta=(a,\omega;\lambda)$ of
 $X\times _{X\tw}\tii T^*X\tw\times_{\ft ^*\tw}
\ft^*$ we define the  point module
 $\delta^\zeta=\tii \DD_X\otimes _{\tii \AA_X}\O_\zeta$.
If $\zeta\tw=(\om,\la)$ is the corresponding point of
$\tii T^*X\tw\times_{\ft ^*\tw} \ft^*$
 then we have
$\tii \DD_X\otimes _{Z(\tii \DD_X)} \OO_{\zeta\tw}\con  End_\k(\delta^\zeta)$.

\medskip

We finish the section with a technical lemma to be used in section
\ref{Splitting on  Springer fibers}.

\lemm\lab{sdvig} Let $\nu=d(\eta)$ be an integral character.
Define a morphism $\tau_\nu$ from $\tii
T^*X\tw\times_{\ft^*\tw}\ft^*$ to itself by
$\tau_\nu(x,\la)=(x,\la+\nu)$. Then the Azumaya algebras
$\tii\DD_X$ and $\tau_\nu^*(\tii \DD_X)$ are canonically
equivalent.

\pf Recall that to establish an equivalence between two Azumaya
algebras $\CA$, $\CA'$ on a scheme $Y$ (i.e. an equivalence between their
categories of modules) one needs to provide a locally projective
module $M$ over $\CA\otimes_{\OO_Y}(\CA')^{op}$ such that
$\CA\iso End_{(\CA')^{op}}(M)$, $\CA'\iso End_{\CA}(M)$.
The sheaf $\pi_*(\DD_{\tii X})^{T,\eta}$ of sections of $\pi_*(\DD_{\tii X})$
which transform by the character $\eta$ under the action of $T$
carries the structure of such a module. \epf

\se{\bf
Localization of $\fg$-modules to $\DD$-modules on the flag variety
}

This crucial section extends the basic result of
\cite{BB}, \cite{BrKa} to positive characteristic.

\sus{
The setting
}
We define relevant triangulated categories of $\fg$-modules and
$\DD$-modules and  functors between them.

\sss{
Semisimple group $G$
}
\lab{Semisimple group G}
Let
$G$ be a semisimple simply-connected
algebraic group over $\k$.
Let $B=T\cdot N$ be a Borel subgroup
with
the unipotent radical  $N$
and a Cartan subgroup $T$.
Let $H$ be the (abstract) Cartan group of $G$ so that
$B$ gives isomorphism
$\io_\fb=(T\con B/N\cong H)$.
Let $\fg,\fb,\ft,\fn,\fh$ be the corresponding Lie algebras.
The weight lattice
$\Lambda=X^*(H)$ contains the set of roots $\De$ and of positive roots $\De^+$.
Roots in  $\De^+$ are identified with
$T$-roots
in $\fg/\fb$
via the above ``$\fb$-identification'' $\io_\fb$. Also,
$\La$ contains the root lattice $Q$ generated by
$\De$,
the dominant cone $\La^+\sub \La$ and
the semi-group $Q^+$ generated by $\Delta^+$.
Let $I\sub \De^+$  be the set of simple roots.
For a root $\al\in\De$ let $\al\mm\ch\al\in\ch\De$ be the corresponding coroot.

Similarly, $\io_\fb$ identifies
 $N_G(T)/T$
with  the Weyl group
$W\sub \Aut(H)$.
Let  $W_\aff\df
W\ltimes Q\sub
W_\aff'\df
W\ltimes \La
$
be
the affine Weyl group and the extended affine Weyl group.
We have  the standard action of $W$ on
$\Lambda$, $w:\lambda\mapsto w(\lambda)=w\cd\la$, and the
$\rho$-shift gives the dot-action
$w:\lambda \mapsto w\bu \lambda= w\bu_\rho \lambda\df\ w(\la+\rho)-\rho$
which is centered at $-\rho$,
where $\rho$ is the half sum of positive roots.
Both actions  extend
to  $W_\aff' $
so that $\mu\in \La$ acts by the $p\mu$-translation.
We will indicate
the dot-action by writing $(W,\bu)$,
this is really
the action of the $\rho$-conjugate
$ ^\rho W$
of the subgroup $W\sub W_\aff'$.

Any
weight $\nu\in\La$ defines
a line bundle $\OO_{\BB,\nu}=\OO_\nu$ on
the flag variety
$
\BB\cong G/B
$, and a standard $G$-module
$V_\nu\df\ \Hh^0(\BB,\OO_{\nu^+})$ with extremal weight $\nu$,
here
$\nu^+$ denotes the dominant $W$-conjugate of $\nu$
(notice that a dominant weight corresponds to a semi-ample
line bundle in our normalization). We will also write $\OO_\nu$
instead of $\pi^*(\OO_\nu)$ for a scheme $X$ equipped with a map
$\pi:X\to\BB$ (e.g. a subscheme of $\gtil$).

We let $\NN\subset \fg^*$ denote the nilpotent cone, i.e.
the zero set of invariant polynomials of positive degree.

\sss{ Restrictions on the characteristic $p$ }\label{restrp}
 Let $h$ be the
maximum of Coxeter numbers of simple components of $G$. If $G$ is
simple then $h= \lb \rho,\ch\al_0\rb+1$ where $\ch\al_0$ is the
highest coroot. We mostly work under the
 assumption  $p> h$, though some intermediate statements are proved
under weaker assumptions;
 a straightforward
 extension of
the main Theorem \ref{Localizationthm} with weaker assumptions on $p$ is 
recorded in 
the sequel paper \cite{sing}.
The main result is obtained for a regular
Harish-Chandra central character, and the most interesting case is
that of an integral Harish-Chandra central character; integral
regular characters exist only for $p\geq h$, hence
our choice of 
restrictions\footnote{The case $p=h$ 
is excluded because for $G=SL(p)$, $p=h$ is not very good and
$\g\not \cong \g^*$ as $G$-modules.}
on $p$.

Recall that a prime is called good if it does not coincide with
a coefficient of a simple root in the highest root \cite[\S 4]{SS},
and $p$ is very good if it is good and $G$ does not contain a factor
isomorphic to 
$SL(mp)$ \cite[3.13]{Sl}. 
 We will need a crude
observation that
 $p> h\ \imp $  very good
$\ \imp\ $ good. 

For $p$ very good $\g$  carries a non-degenerate invariant bilinear form;
also $\g$ is simple provided that $G$ is simple,
 see, e.g., \cite[6.4]{Ja}. 
We will occasionally identify $\fg$ and $\fg^*$
as $G$-modules.
This will identify the nilpotent cones $\NN$ in $\fg$ and $\fg^*$.

\sss{
The sheaf $\tii\DD$
}
\label{redu}
 Our main object is the sheaf $\DD=\DD_\B$ on  the flag variety.
Along with
$
\DD
$ we will consider its deformation $\Dtil$
 defined by the  $H$-torsor
$\tii\BB\df\ G/N @>\pi>> \BB$
as in subsection \ref{Torsors}.
Here  $G\tim H$  acts on  $\tii\BB=G/N$ by
$(g,h)\cd aN\df\ gahN
$, and this action differentiates to
a map $\fg\pl\fh@>>>\tii\TT_\BB$ which extends to
$U(\fg)\ten U(\fh)
@>>>\tii\DD_\BB$.
Then $\Dtil=\pi_* (\DD_{\tii\BB})^H$
is a deformation over $\fh^*$ of
$
\DD\cong \Dtil \otimes_{S(\fh)} \k_0
$.


The corresponding deformation
of $T^*\BB$
will be denoted
$\gtil
=
\tii T^*\BB
=\{(\fb,x)\ |\  \fb\in \BB,\ x|_{rad(\fb)}=0\}$; we have projections
$\PR_1:\gtil\to \fg^*$, $\PR_1(\fb,x)=x$ and $\PR_2: \gtil \to \fh^*$
sending $(\fb,x)$ to $x|_\fb\in (\fb/rad(\fb))^*=\fh^*$;
they yield a map
$\PR=\PR_1\times\PR_2:\gtil\to \fg^*\times_{\fh^*//W}
\fh^*$.
According to section \ref{Torsors}
the sheaf $\Dtil$ is an Azumaya algebra
on $\gtil\tw\times_{\fh^*\tw}\fh^*$ where $\fh^*$ maps to $\fh^*\tw$ by the
Artin-Schreier map.

We denote for any $B$-module $Y$  by
$Y\z$
the sheaf
of sections of the associated $G$-equivariant vector bundle on $\BB$.
For instance,
vector bundle
$\TT_\BB=[\fg/\fb]\z$ is
generated by the space $\fg$ of global sections,
so $\fg$ and  $\OO_\BB$
generate
$\DD$ as an $\OO_\BB$-algebra,
and  one finds that $\DD$ is a quotient of the smash product
$U\z=\OO_\BB\# U(\fg)$
(the semi-direct tensor product),
by the two-sided ideal
$\fb\z\cd U(\fg)\z$. So
$\DD
=
[U(\fg)/\fb U(\fg)]\z
$,
and the fiber (with respect to the left $\OO$-action)
at $\fb\in\BB$ is
$\OO_\fb\ten_\OO\DD
\cong
U(\fg)/\fb U(\fg)$.
Similarly,
$\tii\DD
=
[U(\fg)/\fn U(\fg)]\z$.

\sss{Baby Verma and point modules}\lab{baby}
Here we show that  $\tii \DD$ can be thought of as the sheaf of
endomorphisms of the ``universal baby Verma module''.

Recall the construction of the
 baby Verma module over $U(\fg)$. To define it one fixes
a Borel $\fb=\fn\oplus \ft\subset \fg$, and elements
$\chi\in \fg^*\tw
$,
$\lambda\in \ft^*$, such that $\chi|_{\fn\tw}=0$,
$\chi|_{\ft \tw}=\AS(\lambda)$ (see \ref{Torsors} for notations). For such
a triple $\zeta=(\fb, \chi; \lambda)$ one sets
 $M_\zeta=U_\chi(\fg)\otimes _{U(\fb)} \k_{\lambda}$,
where $U_\chi(\fg)$ is as in \ref{Zp}, and $ \k_{\lambda}$ is the one
dimensional $\fb$-module given by the map $\fb\to \ft \overset{\lambda}{\to}
\k$.

On the other hand, a triple
$\zeta=(\fb, \chi; \lambda)$ as above defines a point of
 $\tilde \fg^*\tw\tim_{\fh^*\tw} \fh^*$ (here we use the isomorphism
$\ft\cong \fh$ defined by $\fb$); thus
we have the corresponding point module $\delta^\zeta$ over $\tii \DD$
(see \ref{Torsors}). Pulling  back this module under the homomorphism
$U(\fg)\to \Gamma(\tii \DD)$ we get a $U(\fg)$-module (also denoted
by $\delta^\zeta$).

\spro We have $\delta^\zeta\cong M_{\fb,\chi;\lambda+2\rho}$.

\pf Let $\fn^-\subset \fg$ be a maximal unipotent subalgebra opposite to $\fb$,
and set $U_\chi(\fn^-)=U_{\chi|_{\fn^-}}(\fn^-)$.
Is suffices to check that there exists a
vector $v\in \delta^\zeta$
such that (1) the subspace $\k v$ is $\fb$-invariant, and
 $\k v\cong \k_{\lambda+2\rho}$; and (2) $\delta^\zeta$ is a free
 $U_\chi(\fn^-)$-module with generator $v$. These two statements follow
from the next lemma, which is checked by a straightforward computation
in local coordinates.

\slemm
Let $\fa$ be a Lie algebra acting\footnote{
An action of a Lie algebra $\fa$ on a variety $X$ is an action of $\fa$
on $\OO_X$ by derivations. Equivalently, it is a Lie algebra homomorphism
from $\fa$ to the algebra of vector fields on $X$.
}
on a smooth variety $X$ and
let $\tii X\to X$ be an $\fa$-equivariant torsor for a torus $T$.
Let $\zeta=(x,\chi;\lambda)$ be a point of $X\times _{X\tw}
\tii T^*X\tw\times_{ \ft^*\tw}\ft^*$, and $\delta^\zeta$ be the corresponding
point module. Let $v\in \delta^\zeta$ be the canonical generator,
$v=1\otimes 1$.

a) If  $x$ is fixed by $\fa$
then $\fa$ acts on $v$ by
$\la_x-\omega_x$, where:\
{\em (1)}
the character
$\lambda_x:\fa\to \k$
is the pairing of
$\lambda\in\ft^*$
with the
action of
$\fa$ on the fiber $\tii X_x$, and
{\em (2)}
the character $\omega_x:\fa\to \k$
is the action of  $\fa$ on the
fiber at $x$ of the canonical bundle $\om_X$.\footnote{
For a section $\Omega$ of $\om_X$ near $x$
and $\xi\in \fa$,\
$\Lie_\xi(\Omega)|_{x}=\omega_x(\xi)\cdot \Omega|_{x}$.
}

b) If, on the other hand,
the action is simply transitive at $x$
(i.e. it induces an isomorphism
$\fa\iso T_x X$),
then the map $u\mapsto u(v)$ gives an isomorphism
$U_{\chi_x}(\fa)\iso \delta^\zeta$; here $\chi_x\in \fa^*\tw$
is the pull-back of $\chi\in \tii T^*_xX$ under the action map. \ \epf

\sss{
The ``Harish-Chandra center'' of $U(\fg)$
}
\lab{HCcenter}
Let now $U=U\g$ be the enveloping algebra of $\g$.
The subalgebra
of $G$-invariants
$\fZ_{HC}\df\ (U\g)^G$
is clearly central in $U\g$.

\slemm Let the characteristic $p$ be arbitrary; the group $G$
is simply-connected, as above.

(a)
The map
$U(\fh)@>>>\Ga(\BB,\tii\DD)$
defined by the $H$-action on $\tii\BB$
gives an isomorphism
$
U(\fh)\con \Ga(\BB,\tii\DD)^G
.$

(b)
The map
$U^G
@>>>
\Ga(\BB,\tii\DD)^G\cong S(\fh)
$
gives an  isomorphism
$
U^G@>{i_{HC}}>> S(\fh)^{(W,\bu)}
$ (the ``Harish-Chandra map'').
For good $p$ this isomorphism is strictly compatible with filtrations,
where the filtration on $Z_{HC}$ is induced by the canonical
filtration on $U$, while the one the target is induced by the filtration
on $S(\fh)$ by degree.

(c)
The map
$U(\fg)\ten S(\fh)@>>>\Ga(\BB,\tii\DD)$
factors through
$
\tii U\df\ U\ten_{\fZ_{HC}} S(\fh)
$.

\pf
We borrow the arguments from \cite{Mi}.
In (a),\
$
\Ga(\BB,\tii\DD)^G
=\
\Ga(\BB,[U/\fn U]\z)^G
\cong
[U/\fn U]^B
\subb
U(\fb)/\fn U(\fb)
\cong
U(\fh)
$,
and the inclusion is an equality, as one sees
by calculating invariants for a Cartan
subgroup $T\sub B$.

For  (b), the map
$U@>>>\Ga(\BB,\tii\DD)$
restricts to a map
$
U^G
@>{i_{HC}}>>
\Ga(\BB,\tii\DD)^G
\cong
U(\fh)
$,
which fits into
$
U^G\sub
U\sur U/\fn U
\subb\ U(\fb)/\fn U(\fb)
\cong U(\fh)
$.
So,
$U^G\sub \fn U+U(\fb)$
and
$i_{HC}$ is the composition
$U^G\sub \fn U+U(\fb)\sur
[\fn U+U(\fb)]/\fn U
\cong U(\fh)$.
On the other hand a choice of a Cartan subalgebra $\ft\sub\fb$
defines an opposite Borel subalgebra
$\barr\fb$
with
$\barr\fb\cap \fb=\ft$ and $\barr\fb=\
\barr\fn\ltim \ft$.
Let us use the
$B$-identification
$\io_\fb:\fh^*\cong \ft^*
$
from \ref{Semisimple group G}
to carry over the dot-action of $W$
to
$\ft^*$
(now the shift is by
$\io_\fb(\rho)=\rho_{\barr\fn}$,
the half sum  of $T$-roots in $\barr\fn$
).
According to  \cite[9.3]{Ja}, an argument of \cite{KW}
shows that
for any simply-connected semisimple group, regardless of $p$,
the projection
$U=
(\fn U+U\barr\fn)\pl U(\ft)
@>>>
U(\ft)
$
restricts to the  Harish-Chandra
isomorphism
$
 \fZ_{HC}
@>{ \io_{ \fn,\barr\fn} } >>
S(\ft)^{ W,\bu}
$.
Therefore,
$i_{HC}=\io_\fb\ci \io_{\fn,\barr\fn}$ is an isomorphism
$
 \fZ_{HC}
\con
S(\fh)^{ W,\bu}
$.

Strict compatibility with filtrations follows from the fact that the
homomorphism $U\to \Gamma (\tii \DD)$ is strictly compatible with filtrations.
The latter follows from injectivity of the induced map on the associated 
graded algebras: $S(\g)={gr}(U)\to \Gamma(\OO_{\gt^*})\cong
{gr}(\Gamma(\tii \DD))$. Here the last isomorphism holds for good $p$,
because of vanishing of higher cohomology $H^{>0}(\BB,gr(\tii \DD))=
H^{>0}(\gt^*,\OO)$. This cohomology
vanishing for good $p$
 follows from \cite{klt}, cf. the proof of Proposition \ref{Dcohomology}
 below.
Injectivity of 
the map $\OO(\g^*)\to \Gamma (\OO_{\gt^*})$ follows from the fact that
the morphism $\gt^*\to \g^*$ is dominant. This latter fact
is a consequence of \cite[6.6]{Ja}, which claims that every
element in $\fg^*$ annihilates the radical of some Borel
subalgebra by a result of \cite{KW}.

Finally, (c)  means that
the two maps from $\fZ_{HC}$ to $\Ga(\BB,\tii\DD)$,
via $U$ and $S\fh$, are the same -- but this is the definition of the second
map.\epf


\sss{
The center of
$U(\fg)$
\cite{Ve,KW,MR1}
}
\lab{The center}
For a very good $p$ the center
$\fZ$
of  $U$
is a combination of the Harish-Chandra part
(\ref{HCcenter}) and the Frobenius part (\ref{Zp}):
$$
\fZ\conn\
\fZ_{Fr}
\ten_{\fZ_{Fr}\cap\fZ_{HC}}
\fZ_{HC}
\cong
\
\OO(\fg^*\tw\tim_{\fh^*\tw//W}\fh^*//(W,\bu))
.$$
Here, $//$ denotes the invariant theory quotient, the map
$\fg^*\tw
@>>>
\fh^*\tw//W$
is the adjoint quotient, while the map
$\fh^*//(W,\bu)
@>>>
\fh^*\tw//W
$
comes from the Artin-Schreier map
$\fh^*
@>{AS}>>
\fh^*\tw
$ defined in \ref{Torsors}.

\sss{
Derived categories of sheaves supported on a subscheme
}
\lab{Derived categories of sheaves supported on a subscheme}

 Let $\AA$ be  a coherent sheaf on a Noetherian scheme $\fX$
equipped with an associative $\OO_\fX$-algebra structure.
We denote by $mod^c(\AA)$ the abelian
category of  coherent
$\AA$-modules. We also use notations
$\Coh(\fX)$ if $\AA=\OO_\fX$
and
$mod^{fg}(\AA)$ if $\fX$ is a point.


We denote by $mod^c_{\hattt\fY}(\AA)$
the full subcategory of coherent $\AA$-modules
supported set-theoretically in
$\fY$, i.e.,
killed by some power of the ideal sheaf $\II_\fY$.
The following statement is standard.

\slemm
a) The tautological functor identifies the bounded derived
category
$
\Db(mod^c_{\hattt\fY}(\AA))
$
with a full subcategory
in
$
\Db(mod^c(\AA))
$.

b) For $\F\in \Db(mod^c(\AA))$ the following conditions are equivalent:\
i) $\F\in \Db(mod^c_{\hattt\fY}(\AA))$,\
ii) $\F$ is killed by a power of the ideal  sheaf $\II_\fY$, i.e.
the tautological arrow $\II_\fY^n\otimes _\OO \F\to \F$ is zero for some $n$;\
iii)  the cohomology sheaves of $\F$ lie in $mod^c_{\hattt\fY}(\AA)$.

 \pf
 In (a) we can replace
 $mod^c$ with $mod^{qc}$ (since $\AA$ is coherent,
$D(mod^{c}(\AA))$
 is a full subcategory of
 $D(mod^{qc}(\AA))$,
 and the same proof works for
 $D(mod^{c}_{\hattt\fY}(\AA))$
 and
 $D(mod^{qc}_{\hattt\fY}(\AA))$).
 Now it suffices to know that
 each sheaf in $mod^{qc}_{\hattt\fY}(\AA)$ embeds into an object of
 $mod^{qc}_{\hattt\fY}(\AA)$ which is injective in
 $mod^{qc}(\AA)$ (\cite[Proposition I.4.8]{Ha}).
 This follows from the
 corresponding statement for quasi-coherent  sheaves of $\OO$ modules
(see e.g. \cite{Ha}, Theorem I.7.18 and its proof),
since we
 can get a quasicoherent  injective sheaf  of $\AA$-modules from an
 injective quasicoherent
 sheaf of $\OO$-modules by coinduction.

b) Implications (i)$\Rightarrow$(ii)$\Rightarrow$(iii)
are clear by definitions, and
(iii)$\Rightarrow$(i)
is clear from (a).
\epf

\sss{
Categories of
modules with a generalized
Harish-Chandra character
}
\lab{Categories}


Let us apply
\ref{Derived categories of sheaves supported on a subscheme}
to $\tii\DD$ and $U$ (or $\tii U$),
considered as coherent sheaves over the spectra
$\tii T^*\BB\tw$ and $\fg^*\tw$ of
central subalgebras.
The interesting categories are
 $
mod^c(\DD^\la)
\sub
mod^c_{\hattt\la}(\tii\DD)
\sub
mod^c(\tii\DD)
$.
Here,
$
mod^c_{\hattt \la}(\tii\DD)
\dff
mod^c_{T^*_{AS(\la)}\BB\tw}(\tii\DD)
$
consists of those objects in $mod^c (\tii \DD)$ which are  killed by
a power of the maximal ideal $\la$
in $U\fh$.

For $\la\in\fh^*$
denote by
$U^{ \la}$  the specialization of
$U$ at the image of $\la$ in
$\fh^*//W
$
$=
\Spec(\fZ_{HC})$, i.e.,
the specialization of
$\tii U$ at $\la\in\fh^*$.
There are
analogous abelian categories
$
mod^{fg}(U^\la)\sub
mod^{fg}_{\hattt \la}(U)
\sub
mod^{fg}(U)
$,
where the category
$mod^{fg}_{\hattt \la}(U)
\dff
mod^c_{\fg^*\tw_\la}(U)
$
for
$\fg^*\tw_\la
\df
$
$
\fg^*\tw \tim_{\fh^*//W\tw} {AS(\la)}
$,
consists of $U$-modules killed by a power
of the maximal ideal in $\fZ_{HC}$.
The corresponding triangulated categories are
$
\Db(mod^{fg}(U^\la))
@>>>
$
$\Db(mod^{fg}_{\hattt \la}(U))
$
$
\sub
\Db(mod^{fg}(U))
$.

\sss{
The global section
functors on $D$-modules
}
\lab{globase}
Let $\Ga= \Ga_\OO$ be the functor of global sections on the
category
$mod^{qc}(\OO)$
of quasi-coherent sheaves on
$\B$
and let
$\Rr\Ga= \Rr\Ga_\OO$ be the derived functor
on $\Dd(mod^{qc}(\OO))$.
Recall
 from \ref{HCcenter} that
 the action of $G\tim H$ on $\tii\BB$ gives
a map
$
\tii U
@>>>
\Ga(\Dtil)$,
this gives a functor
$
mod^{qc}(\tii\DD)
@>{\Ga_{\tii\DD}}>>
mod(\tii U)$, which can be derived to
$
\Db(mod^{qc}(\tii\DD))
@>{R\Gamma_{\tii \DD}}>>
D(mod(\tii U))
$ because the category
of modules has direct limits. 
This derived
functor commutes with
 the  forgetful functors, i.e.
$
\mbox{Forg}^{\tii U}_\k
\ci
\Rr\Gamma_{\tii\DD}
=\
 \Rr\Gamma  \ci \mbox{Forg}^\Dtil_\OO
$
where
$
{\mbox{Forg}^\Dtil_\OO}:mod^{qc}(\tii\DD)
\to
mod^{qc}(\OO)
$,
$
\mbox{Forg}^{\tii U}_\k:
mod(\tii U)\to Vect _\k
$
are
the  forgetful functors.
This is true since
the category
      $mod^{qc}(\tii\DD)$ has enough objects acyclic for
    the functor of global sections $\Rr\Ga$ (derived in quasi-coherent
    $\OO$-modules) -- namely, if $U_i@>{j_i}>>\ \BB,\ \ii$,
    is an affine open cover then for any object $\F$ in
    $mod^{qc}(\tii\DD)$ one has
    $\FF\inj\ \pl_\ii\ (j_i)_*(j_i)^*(\F)$.
Since
$\Ga$ has finite homological dimension, $\Rr\Gamma_{\tii\DD}$
actually  lands in the bounded derived category.

\slemm \lab{coherence} The (derived)  functor of global sections
preserves coherence, i.e. it sends the full subcategory  $\Db( mod^c(\tii\DD))
\subset \Db(mod^{qc}(\tii \DD))$ into the full subcategory
$ \Db(mod^{fg}(\tii U))\subset
\Db(mod(\tii U))$.

\pf First notice that since $\tii U$ is noetherian, $
\Db(mod^{fg}(\tii U)) $ is indeed identified with
$\Db_{fg}(mod(\tii U))$, the full subcategory in $\Db(mod (\tii
U))$ consisting of complexes with finitely generated cohomology.

The map $\tii U@>>>\Ga\tii\DD$ is compatible
with natural filtrations and it produces
a proper map $\mu$ from
$\Spec(Gr(\tii\DD))
=
G\tim_B\ \fn^\perp$  to the affine
variety $
\Spec(Gr(\tii U))\cong\fg^*\tim_{\fh^*//W}\fh^*$
(here, $gr(\fZ_{HC})\cong \OO(\fh^*)^W$
by Lemma \ref{HCcenter}(b)).
Any coherent $\tii\DD$-module $M$ has a coherent filtration,
i.e., a lift to a filtered $\tii\DD$-module
$M_\bu$
such that $gr(M_\bu)$ is coherent for $Gr(\tii\DD)$.
Now, each $R^i\mu_*(gr(M_\bu))$
is a  coherent sheaf on
$\Spec(Gr(\tii U))$, i.e,
$\Hh^*(\BB,gr(M_\bu))$
is a finitely generated module over
$Gr(\tii U)$.
The filtration on $M$ leads to a spectral sequence
$\Hh^*(\BB,gr(M))\ \imp\ gr(\Hh^*(\BB,M))$,
so
$gr(\Hh^*(\BB,M))$
is
a subquotient of
$\Hh^*(\BB,gr(M))$, and therefore  it is also finitely generated.
Observe that
the induced filtration on
$
\Hh^*(\BB,M)
$
makes it into a filtered module for
$
\Hh^*(\BB,\DD)
$
with its induced filtration.
Since
$\tii U
@>>>
\Hh^0(\BB,\DD)
$
is a map of filtered rings,
$
\Hh^*(\BB,M)
$
is also  a filtered module for
$\tii U$.
Now, since
$
gr (\Hh^*(\BB,M))
$
is a finitely generated module for
$
gr(\tii U)
$,
we find that
$
\Hh^*(\BB,M)
$ is finitely generated  for
$
\tii U
$.
This shows that
$R\Gamma_{\tii\DD}$
maps
$
\Db(mod^c(\tii\DD))
$
to 
$
\Db_{fg}(mod(\tii U))\cong \Db(mod^{fg}(\tii U))
$.
\epf

It follows
 from \ref{HCcenter} that
the  canonical map $\tii U\to \DD^\la$
factors
for any
$\la
\in\fh^*$
to
$U^\la\to \DD^\la
$.
So, as above,
we  get functors
$
mod^c_{\hattt \la}(\tii\DD)
@>{\Gamma_{\tii\DD,\la}}>>
mod^{fg}_{\hattt \la}(\tii U)$,
$
mod^c(\DD^\la)
@>{\Gamma_{\DD^\la}}>>
$
$
mod^{fg}(U^\la)
$.
The derived functors
$
\Db(mod^{c}_{\hattt \la}( \tii\DD))
@>{\Rr\Gamma_{\tii\DD,\la} }>>
\Db(mod^{fg}_{\hattt \la}( U))$,\
$\Db(mod^{c} (\DD^\la))
@>{\Rr\Gamma_{\DD^\la}}>>
\Db(mod^{fg}( U^\la))$
are defined
and compatible  with the forgetful functors.


\sus{Theorem}
\lab{Localizationthm}
({\bf The main result.})
{\em
Suppose\footnote{The restriction on $p$ is discussed in \ref{restrp} 
above.}  
that  $p> h$.
For any
regular
$\la\in\fh^*$  the global section functors
 provides equivalences
 of triangulated categories
\begin{equation}
\label{equivD0}
\Rr\Gamma_{\DD^\la}:
\Db(mod^c(\DD^\la))
\con
\Db(mod^{fg}(U^\la));
\end{equation}
\begin{equation}
\label{equivD}
\Rr\Gamma_{\tii\DD,\la}:
\Db(mod^c_{\hattt \la}(\tii\DD))
\con
\Db(mod^{fg}_{\hattt \la}( U)).
\end{equation}
}

\begin{Rem}
\lab{Remark 1}
In the characteristic zero case Beilinson-Bernstein
(\cite{BB}, see also \cite{Mi}),
proved that for a dominant $\lambda$
the functor of global sections provides an equivalence between
the abelian categories
$mod^c(\DD^\lambda)
@>>>
mod^{fg}(U^\lambda)$. The
analogue for crystalline differential operators
in characteristic
$p$ is evidently false: for any line bundle $\L$ on $\B$ the line bundle
$\L^{\otimes p}$ carries a natural structure of a $\DD$-module
(\ref{Lagrangian}); however
$\Rr^i\Gamma (\L^{\otimes p})$ may certainly be nonzero for $i>0$.
Heuristically, the analogue of characteristic zero results
about dominant weights is not available in characteristic $p$, because
a weight can not be dominant (positive) modulo $p$.


However, for a generic $\la\in\fh^*$
it is very easy to see that
global sections give an equivalence of abelian categories
$mod^c(\DD^\la)
@>>>
mod^{fg}(U^\la)$.
If $\io(\la)$ is regular, 
the twisted cotangent bundle $T^*_{\io(\la)}\BB$
is affine, so $\DD^\la$-modules are equivalent to modules for
$\Ga(\BB,\DD^\la)$,
and
$\Ga(\BB,\DD^\la)=U^\la$
is proved in
\ref{Dcohomology}.
\end{Rem}


\begin{Rem}
Quasicoherent and ``unbounded'' versions of the equivalence, say $
D^?(mod^{qc}(\DD^\la)) @>{\Rr\Gamma_{\DD^\la}}>> D^?(mod(U^\la))
$, $?=+, \, -$ or $b$, follow formally from the coherent versions
since
$R\Ga_{\DD^\la}$ and its adjoint (see \ref{Localization}) commute
with homotopy direct limits.
For  completions to formal neighborhoods see
\ref{Equivalences on formal neighborhoods}.
\end{Rem}

\sss{The strategy of the proof of Theorem \ref{Localizationthm} }
We concentrate on the second statement, the first one follows (or
can be proved in a similar way). First we observe that the functor
of global sections
 $\Rr\Ga_{\tii \DD,\la}:\
\Db(mod^c_\lambda(\tii\DD ))\ra\ \Db(mod^{fg}_\la(U))$ has
left adjoint --  the localization functor 
$\LL^{\hatt\la}$ . A straightforward modification of a known
characteristic zero argument shows that the composition of the two
adjoint functors in one order is isomorphic to identity. The
theorem then follows from a certain abstract property of the
category  
$\Db(mod^c_{ \la}(\tii\DD))$
which we call the {\it (relative) Calabi-Yau} property (because
the derived category of coherent sheaves on a Calabi-Yau manifold
provides a typical example of such a category). This property of
 $\Db(mod^c_{ \la}(\tii\DD))$ will be
derived from triviality of the canonical class of 
$\gtil$.


\begin{Rem}
One can give another proof of Theorem \ref{Localizationthm} with a
stronger restriction on characteristic $p$, which  is closer to the
original proof by Beilinson and Bernstein \cite{BB} of the characteristic zero
statement. (A similar proof appears in an earlier preprint version
of this paper.) Namely, for  fixed weights $\lambda$, $\mu$ and
large $p$ one can use the Casimir element in $Z_{HC}$ to show that
the sheaf $\OO_\mu\otimes \MM$ is a direct summand in the sheaf of
$\g$ modules $V_\mu\otimes \MM$ for a $\DD_\la$-module $\MM$ (where
$\la$ is assumed to be integral and regular). Choosing $p$, such
that this statement holds for a finite set of weights $\mu$, such
that $\OO_\mu$ generates $\Db(Coh(\BB))$, we deduce from Proposition
\ref{Dcohomology} that the functor $R\Gamma$ is fully faithful.
Since the adjoint functor $\LL$ is easily seen to be fully faithful
as well (see Corollary \ref{compid}), we get the result.
\end{Rem}

\sus{ Localization functors} \lab{Localization}

\sss{Localization for categories with
generalized Harish-Chandra
character}
We start with the localization functor $Loc$ from (finitely generated)
$U$-modules to $\tii\DD$ modules, $Loc(M)= \tii \DD \otimes
_{U} M$.
Since $U$ has finite homological dimension
it has a left derived functor
$
\Db(mod^{fg}(U))
@>{\LL}>>
\Db(mod^c(\tii\DD))
$.
Fix $\lambda\in \h^*$, for any $M\in
\Db(mod_\lambda^{fg}(U))$ we have a canonical decomposition
$\LL(M)= \oplusl_{\mu \in W\bu \lambda} \LL^{\lambda\to\mu}(M)$
with $\LL^{\lambda\to\mu}(M)\in
\Db(mod^c_\mu (\tii \DD))$.
Localization with the generalized character $\la$
is the functor
$
\LL^{\hatt \lambda}
\df
\LL^{\lambda\to\lambda}:
\Db(mod^{fg}_\la(U))
@>>>
\Db(mod^c_\la(\tii\DD))
$.

\lemm The functor $\LL$ is left adjoint to $\Rr\Ga$, and
$\LL^{\hatt \la} $ is left adjoint to
$\Rr\Gamma_{\tii\DD,\la}$.

\pf It is easy to check that the functors between abelian
categories $\Gamma: mod^{qc}(\tii \DD)\to mod(U)$, $Loc:mod(U)\to
mod^{qc} (\tii \DD)$ form an adjoint pair. Since $mod^{qc}(\tii
\DD)$ (respectively, $mod(U)$) has enough injective (respectively,
projective) objects, and the functors $\Gamma$, $Loc$ have bounded
homological dimension
 it follows that their derived functors form
an adjoint pair.
Lemma \ref{coherence} asserts
that $R\Gamma$ sends $\Db(mod^c(\tii \DD))$ into
$\Db(mod^{fg}(U))$; and it is immediate to check that $\LL$ sends
$\Db(mod^{fg}(U))$ to $\Db(mod^c(\tii \DD))$. This yields the
first statement. The second one follows from the first one. \epf

\sss{Localization for categories with
a
fixed Harish-Chandra
character
}
We now turn to the categories appearing in
\eqref{equivD0}. The functor $Loc$ from the previous subsection
restricts to a functor $Loc^\lambda: mod^{fg}(U^\lambda) \to
mod^c(\DD^\lambda),\
Loc^\lambda(M)=\DD^\la\ten_{U^\la}M$. It has a left derived functor
$\LL^\la:\Dmin(mod^{fg}(U^\la)\to \Dmin(mod^c(\DD^\la)),\
\LL^\lambda(M)=\DD^\la\Lten_{U^\la}M
$. Notice that
the algebra $U^\lambda$ may a'priori have infinite homological
dimension\footnote{For regular $\la$ the finiteness of homological
dimension will eventually follow from the equivalence
\ref{Localizationthm}.}, so $\LL^\la$ need not  preserve  the
bounded derived categories. The next lemma shows that
it does for regular $\la$.
\lemm
\lab{commut}
a) $\LL^\la$ is left adjoint to the functor
$
\Dmin (mod^c(\DD^\la))
@>{R\Gamma_{\DD^\la}}>>
\Dmin(mod^{fg}(U^\la))$.

b) For regular $\la$ the localizations at $\la$ and the generalized
character $\la$ are compatible, i.e., for the obvious functors
$\Dmin(mod^{fg}(U^\la))
@>{i}>>
\Dmin(mod_\la^{fg}(U))$ and
$
\Dmin(mod^c(\DD^\la))
@>{\io}>>
\Dmin(mod_\la^c(\tii\DD))$, there is a  canonical isomorphism
$$\iota \circ \LL^\la\cong \LL^{\hatt \la} \circ i,$$
and this isomorphism is compatible with the adjunction arrows in the
obvious sense.

\pf a) is standard.
%
 To check (b) observe that if  $\la$ is regular
    for the dot-action of $W$, then
the projection $\h^* \to \h^*/(W,\bu)$ is etale at $\la$; thus we
have $\OO(\h^*)^{\hatt \la}\Lotimes _{\OO(\h^*/(W,\bu))} \k_\la =
\k$, where $\OO(\h^*)^{\hatt \la}$ is the completion of $\OO(\h^*)$
at the maximal ideal of $\la$. It follows that $\tii \DD^{\hatt \la}
\Lotimes _{U} U^\la= \DD^\la$, where $\tii \DD ^{\hatt \la}= \tii
\DD \otimes _{\OO(\h^*)} \OO(\h^*)^{\hatt \la}$. It is easy to see
from the definition that $\LL^{\hatt\la}(M)\cong \DD^{\hatt\la}\Lten
_{U} M$ canonically, thus we obtain the desired isomorphism of
functors. Compatibility of this isomorphism with  adjunction follows
from the definitions. \epf



\cor The functor $\LL^\la$ sends the bounded derived category $
\Db (mod^c(\DD^\la))$ to $ \Db(mod^{fg}(U^\la))$ provided $\la$ is
regular. \epf

\sus{ Cohomology of $\tii \DD$ } The computation in this section
will be used to check that $R\Gamma_{\tii \DD,\la} \circ \LL^{\hatt
\la}\cong id$ for regular $\la$.

\pro
\lab{Dcohomology} Assume that $p$ is very good. 
Then we have $\tii U\con\Rr\Gamma(\tii\DD)$ and also
$U^\la\con\Rr\Gamma(\DD^\la) $ for $\la\in\fh^*$.

\pf The sheaves of algebras $\DD^\la$, $\tii \DD$ carry filtrations
by the order of a differential operator, the associated graded sheaves
are, respectively, $\OO_{\Nt}$ and $\OO_{\gt^*}$. 
Cohomology vanishing for $\DD$, $\tii \DD$
follows from cohomology vanishing of the 
associated graded sheaves. For $\OO_{T^*\BB}$ this is Theorem 2 of \cite{klt},
which only requires $p$ to be good for $\g$.
The case of $\gt^*$ is a formal consequence. To see this consider a 
 two-step $B$-invariant filtration
on $(\fg/\fn)^*$ with associated graded $\fh^*\oplus (\fg/\fb)^*$.
It induces a filtration on $\gt^*$ considered as a vector bundle on $\BB$. 
The associated graded of the corresponding filtration on
$\OO_{\gt^*}$ (considered as a sheaf on $\BB$) is $S(\fh)\otimes 
\OO_{\Nt}$. Cohomology vanishing of the last sheaf follows from the one
for $\OO_{\Nt}$, and implies one for $\OO_{\gt^*}$.

Furthermore, higher cohomology vanishing for the associated
graded sheaves
 $\OO_{\Nt}=gr(\DD^\la)$, $\OO_{\gt^*}=gr(\tii \DD)$ 
implies that the natural maps
$gr(\Gamma(\DD^\la))\to \Gamma(\O_{\Nt})$, $gr(\Gamma(\tii \DD))
\to \Gamma (\gt^*)$ are isomorphisms.

We will show that the maps $U^\la\to \Gamma (\DD^\la)$,
$\tii U\to \Gamma(\tii \DD)$ are isomorphisms by showing that the 
induced maps on the associated graded algebras are. Here
the filtration on $U^\la$ is induced by the canonical
filtration on $U$, and the one  on $\tii \DD$
is induced by the canonical filtration on $U$ and
the degree filtration on $S(\h)$.

The associated graded rings of $U^\la$, $\tii U$ 
are quotients of, respectively,
$S(\g)$ and $S(\fg)\otimes S(\fh)$. Moreover, in view
of Lemma \ref{HCcenter}(b), they are quotients of, respectively, 
$S(\g)\otimes_{S(\g)^G}\k$ and $S(\g)\otimes _{S(\g)^G}S(\h)$.
It remains to show that the maps
$S(\g)\otimes_{S(\g)^G}\k\to \Gamma (\O_{\Nt})$,
$S(\g)\otimes _{S(\g)^G}S(\h)\to \Gamma (\OO_{\gt^*})$
are isomorphisms. Here the maps are readily seen to be induced by 
the canonical morphisms $\Nt\to \g^*$
and $\gt^*\to \g^*\times _{\h^*/W} \h^*$.

Since $p$ is very good, we have a $G$-equivariant isomorphism
 $\g\cong \g^*$, see \ref{restrp},
 thus it suffices
to show that the global functions on the nilpotent variety
$\NN\subset \g$ map isomorphically to the ring of global functions
on $\Nt\cong \n\times^B G$. Moreover, the \' etale slice theorem of
\cite{BaRi}
shows that for very good $p$
 there exists a $G$-equivariant isomorphism between
$\N$ and the subscheme 
 ${\mathcal U}\subset G$
defined by the $G$-invariant polynomials on $G$ vanishing at the unit
element, cf. \cite[9.3]{BaRi}.
Thus it
reduces the task to showing that the ring of regular
 functions on  ${\mathcal U}$
maps isomorphically to the ring of global functions on
$N \times ^B G$.  This follows once we know that
$\mathcal U$ is reduced and normal and the Springer map 
$N \times ^B G\to \mathcal U$
is birational.
These facts can be found in \cite{St} for all $p$:
 $\mathcal U$ is reduced and normal
by 3.8, Theorem 7, it is irreducible by 3.8, Theorem 1, while the 
Springer map is a resolution of singularities by 3.9, Theorem 1.

Finally, surjectivity of the map $S(\fg)\otimes_{S(\fh)^W}
 S(\fh)\to \Gamma(\O(\gt^*))$
follows from surjectivity established in the previous
paragraph 
by  graded Nakayama lemma;
notice that higher cohomology vanishing for  $\O_{\gt^*}$
implies that $\Gamma(\O_{\Nt})=\Gamma(\O_{\gt^*})\otimes_{S(\fh)} \k$.
Injectivity of this map is clear  from the fact that
$S(\h)$ is free over $S(\h)^W$ for very good $p$
by \cite{De}, cf. also \cite[9.6]{Ja}. Hence 
 $S(\fg)\otimes_{S(\fh)^W}
 S(\fh)$ is free over $S(\g)$, while the map $\gt^*\to \g^*\times_{\h^*/W}\h^*$
is an isomorphism over the open set of regular semisimple elements in 
$\g^*$ for any $p$. \epf

%

\cor
\lab{compid}
a)
The composition $\Rr\Ga_{\tii\DD} \circ \LL:
\Db(mod^{fg}(U)) \to \Db(mod^{fg}(\tii U))$ is isomorphic to the
functor $M\mapsto M\otimes _{\fZ_{HC}} S(\fh)$.

b)
For regular $\la$ the
adjunction map
$id\to \Rr\Ga_{\tii\DD,\la} \circ
\LL^{\hatt\la}$ is an isomorphism
on ${\Db(mod^{fg}_\la(U))}$.

c)
For any $\la$,
the adjunction map is an isomorphism
$id@>>>\Rr\Ga_{\DD^\la}\ci \LL^\la$
on $\Dmin(mod^{fg}(U^\la))$.

\pf For any $U$-module $M$ the action of $U$ on
$\Gamma_{\tii\DD} (\LL(M))$
extends to the action of
$\Gamma (\tii \DD)= \tii U$.
So the adjunction map
$M\to \Gamma_{\tii\DD} (\LL(M))$
extends to
$S(\fh)\otimes _{\fZ_{HC}} M
=\tii U \ten_U M
\to \Gamma_{\tii\DD} \circ \LL (M)$.
Proposition \ref{Dcohomology} implies that if $M$ is a free module
then this map is an isomorphism, while higher derived functors
$R^i\Ga_{\tii\DD} (\LL(M))$, $i>0$, vanish. This yields statement (a).
(c) is proved in the same way using the second claim in
Proposition \ref{Dcohomology}.

To deduce (b) observe that for regular
$\la$ and  $M\in \Db(mod_\la^{fg}(U))$, we have canonically
$M\otimes _{\fZ_{HC}} S(\fh)\cong \pl_W M$.
The adjunction morphism viewed as
$M\to \pl_W\ M$,
equals $\sum _W id_M$ (when $M$ is the
restriction of $U$ to the $n\thh$ infinitesimal neighborhood of $\la$
this follows by restricting
$\tii U\con R\Ga(\tii\DD)$).
Now the claim follows
since $\Rr\Ga_{\tii \DD,\la}
(\LL^{\hatt\la}(M))
$ is one of the summands.
\epf

\sus{Calabi-Yau categories}
\lab{CY}
 We recall some generalities about Serre functors in triangulated
categories; we refer to the original paper\footnote{We
slightly generalize the definition of \cite{BK}, cf. also \cite{BezruKa}.}
 \cite{BK} for details.

Let $\O$ be a finite type commutative algebra over a field; and
let $D$ be an $\O$-linear triangulated category.
A structure of an $\OO$-triangulated category on $D$
is a functor $\RHom_{D/\OO}: D^{op}\times D\to \Db(mod^{fg}(\OO))$,
together with a functorial isomorphism
$Hom_D(X,Y)\cong
\Hh^0(\RHom_{D/\OO}(X,Y))$.

For any quasi-projective variety $Y$, the triangulated category
$\Db(\Coh(Y))$ is equipped with a
canonical anti-auto-equivalence, namely the Grothendieck-Serre
duality $\D_Y=\RHHom_\OO(-,K_Y)$ for the dualizing complex
$K_Y=(Y\to pt)^!\k$.

By an {\it $\OO$-Serre functor } on $D$
we will mean an
auto-equivalence $S:D\to D$
together with
a natural (functorial) isomorphism
$\RHom_{D/\OO}(X,Y)\cong
\D_\O(\RHom_{D/\OO}(Y,SX))$
for all $X,Y\in D$.
If a
Serre functor exists, it is unique up to a unique isomorphism.
An  $\OO$-triangulated category will be called {\it Calabi-Yau} if
for some $n\in \Z$
the shift functor $X\mapsto X[n]$
admits a structure of an $\OO$-Serre functor.

For example, if $X$ is a smooth variety over $k$ equipped with a
projective morphism $\pi:X\to Spec(\O)$ then $D=\Db(\Coh_X)$
is $\OO$-triangulated by
$\RHom_{D/\OO}(\F,\G) \df\ R\pi_* \RHHom(\F,\G)$.
The functor $\F\mapsto \F\otimes \om_X[\dim X]$ is naturally a Serre
functor with respect to $\O$; this is true because
Grothendieck-Serre duality commutes with proper direct images,
and the dualizing complex for a smooth $X$ is $K_X\con\om_X[\dim(X)]$, so
$$
\D_\OO (R\pi_*\RHHom(\F,\G))
\cong
R\pi_*(\D_X \RHHom(\F,\G))
\cong
R\pi_*\RHHom (\G, \F\otimes \om_X[\dim X])
.$$
We will need the following generalization of this fact. Its proof
is straightforward and left to the reader.\footnote{Details of 
the proof can also be found in the sequel 
paper \cite{sing}.}

\lemm
\lab{AzuSerre}
Let
$\CA$ be an Azumaya algebra on
a smooth variety $X$ over $k$, equipped
with a projective morphism $\pi:X\to Spec(\O)$. Then
$\Db(mod^c(\CA))$ is
naturally $\O$-triangulated and the functor  $\F\mapsto
\F\otimes \om_X[\dim X]$ is naturally a Serre functor with respect
to $\OO$.
In particular, if $X$ is a Calabi-Yau manifold
(i.e., $\om_X\cong \OO_X$) then
the
$\O$-triangulated category
$\Db(mod^c(\CA))$ is Calabi-Yau.

Application of the above notions to our situation is based on the
following lemma.
A similar argument was used e.g.
in \cite{MuMu}, Theorem 2.3.

\lemm\lab{abstract} Let
$D$ be a Calabi-Yau $\OO$-triangulated category
for  some
commutative finitely generated algebra $\OO$.
Then a sufficient condition for
a triangulated functor $L:C\to D$
to be  an equivalence is given by

i) $L$ has a right adjoint functor $R$ and the adjunction
morphism $id\to R\circ L$ is an isomorphism, and

ii) $D$ is indecomposable, 
i.e. $D$ can not be written as
$D=D_1\oplus D_2$ for nonzero triangulated categories $D_1,D_2$;
and $C\ne 0$.

\pf
Consider any full subcategory $\CC\sub D$ invariant under the shift functor.
The right orthogonal is the full subcategory
$\CC^\perp=\{y\in D;\ Hom_D(c,y)=0\ \forall c\in \CC\}$.
If $S$ an $\OO$-Serre functor for $D$
then
$S^{-1}:\CC^\perp\to
{^\perp \CC}$ (the left orthogonal of $\CC$), since  for
$y\in \CC^\perp$ and $c\in \CC$ one has
$
\Hh^n\RHom_{D/\OO}(c,y)
=\Hom_D(c,y[n])=
\Hom_D(c[-n],y)=0,\ n\in\Z$,
 hence
$\RHom_{D/\OO}(c,y)=0$,
and then
$D_\OO\RHom_{D/\OO}(S\inv y,c)
=
\RHom_{D/\OO}(c,y)=0$.
In particular, if $D$ is Calabi-Yau relative to
$\O$, then $^\perp\CC=\CC^\perp$.

Now,
condition (i) implies that $L$ is a full embedding, so we will
regard it as the inclusion of a full subcategory $C$ into $D$.
Moreover, for $d\in D$, any cone $y$ of the map
$LR(d)\to d$
is in $C^\perp$.
Therefore,  $y\in{^\perp C}$, and then $d\cong LR(d)\pl y$.
This yields a
decomposition $D=C\oplus C^\perp$. So, condition (ii)
implies that $C^\perp=0$ and $L$ is an equivalence. \epf

 Another useful simple fact is

\lemm\lab{indecompo} (cf. \cite{MuMu}, Lemma 4.2)
 Let $X$ be a connected scheme quasiprojective over a field $\k$,
and let $\CA$ be an Azumaya algebra on $X$. Then the category
$\Db(mod^c(\CA))$ is indecomposable. Moreover, if $Y\subset X$ is
a connected closed subset then $\Db(mod^c_Y(X,\CA))$ is
indecomposable.

\pf  Assume that $\Db(mod^c(\CA))=D_1\oplus D_2$ is a decomposition
invariant under the shift functor. Let $P$ be an indecomposable
summand of the free $\CA$-module. Let $L$ be a very ample line bundle
on $X$ such that $0\ne H^0(L\otimes \HHom_\CA (P,P))
=
\Hom_\CA (P,P\ten L)$. For any
$n\in \Zet$ the $\CA$-module $P\otimes L^{\otimes n}$ is
indecomposable, hence belongs either to $D_1$ or to $D_2$.
Moreover, all these modules belong to the same summand, because
$Hom_\AA(P\otimes L^{\otimes n}, P\otimes L^{\otimes m})\ne 0$ for
$n\leq m$.  If $\F$ is an object of the other summand, then we
have
$Ext^\bu_\AA (P\otimes L^{\otimes n}, \F)=0$ for all $n$.
However,
since $\AA$ is  Azumaya algebra,
$P\ne 0$ is a  locally projective $\AA$-module and $X$ is connected,
$\FF\ne 0$ would imply
$\RHHom_\AA(P,\F)\ne 0$
(this claim reduces to the case when $\AA$ is a matrix algebra and then to
$\AA=\OO_X$).
So $\FF=0$
(otherwise
$H^*(X,\RHHom_\AA (P, \F)\ten L^{\otimes -n})$
could not be
zero for all $n$),
and this proves the first statement.
The second claim follows: for any closed subscheme
$Y'\subset X$ whose topological space equals $Y$,
the image of
$\Db(mod^c(Y',\CA|_{Y'}))$ under the push-forward functor lies in
one summand of any  decomposition
$\Db(mod^c_Y(X,\CA))=D_1\oplus D_2$.
\epf

\sus{Proof of Theorem \ref{Localizationthm}} The canonical line
bundle on $\gtil$ is trivial; hence the same is true for
$\gtil\tw\tim_{\fh^*\tw}\fh^* $, the spectrum of the center of $\tii
\DD$ (see \ref{The center}). Thus Lemma \ref{AzuSerre} shows
that $\Db(mod^c(\tii \DD))$ is  Calabi-Yau with respect to
$\O(\g^*)$.

It follows from the definitions that a full triangulated
subcategory in a Calabi-Yau category with respect to some algebra
$\O$ is again a Calabi-Yau category with respect to $\O$.
Therefore, \eqref{equivD} follows from Corollary \ref{compid}(b) and
Lemmas \ref{abstract}, \ref{indecompo}.

To deduce \eqref{equivD0} from \eqref{equivD} we use Lemma
\ref{commut}(b). It says that the functors $i$, $\iota$ send the
adjunction arrows into adjunction arrows; since $i$, $\iota$ kill
no objects, and the adjunction arrows in
$\Db(mod^c_\la(\tii\DD))$, $\Db(mod_\la^{fg} (U))$ are
isomorphisms, we conclude that the adjunction arrows in
$\Db(mod^c(\DD^\la))$, $\Db(mod^{fg} (U^\la))$ are isomorphisms,
which implies \eqref{equivD0}. \epf

 \se{\bf Localization  with a
generalized Frobenius character }

\sus{
Localization on  (generalized) Springer fibers
}
\lab{Localization_Springer}

The map
$U\to \tii\DD$ restricts to  a map of central algebras
$\OO(\fg^*\tw)\to \OO_{\tii\fg^*\tw}$.
So, the commutative part of the localization mechanism
is the resolution
$\tii\fg^*\tw\to \fg^*\tw$.
Therefore, the specialization of the algebra $U$ to
$\chi\in\fg^*\tw$ will correspond to the
restriction of $\tii\DD$ to the corresponding Springer fiber.

{F}rom here on we keep in
mind
that the Weyl group always acts by the dot
action and we   write $\fX//W$ instead of $\fX//(W,\bu)$
for the invariant theory quotients.

\sss{
Categories with  a generalized character $\chi$
of the Frobenius center
}
\lab{Categories with chi}

Recall that the center
$\fZ=\ \OO(\fg^*\tw\tim_{\fh^*\tw//W}\fh^*//W)$
of  $U$ is generated by subalgebras
$\fZ_{Fr}=\ \OO(\fg^*\tw)$
and
$\fZ_{HC}=\ \OO(\fh^*//W)$
which the map $U(\fg)\ra\ \Ga\tii\DD$ sends   to
the central subalgebras
$
\OO(\tii T^*\BB\tw)
$ and
$S\fh$ of $\tii\DD$
(\ref{The center}).

For $\la\in\fh^*,\ \chi\in\fg^*$,  the notation
$U^\la,U_\chi,U^\la_\chi$ denotes
the specializations of $U$ to the characters $\la,\chi,(\la,\chi)$
of $\fZ_{HC},\fZ_{Fr},\fZ$.
Similarly,  the sheaf of algebras
$\tii\DD$
has specializations $\DD^\la\df\ \tii\DD^\la,
\tii\DD_{\chi},
\DD^{\la}_{\chi}
$.
As in \ref{Derived categories of sheaves supported on a subscheme},
we denote the full subcategories with
a generalized character
$\ze\in\{\la,\chi,(\la,\chi)\}$
of
$\fZ_{HC},\fZ_{Fr}$ or $\fZ$,
by
$
mod^{c}_{\ze}(-)
\sub\
mod^{c}(-)
$,
and one has
$\Db(mod^{c}_{\ze}(-))
\sub
\Db(mod^{c}(-))$.
For later use we notice that
$mod_\chi^{fg}(U)$ can be viewed as the
category
$
mod^{fl}(U^\la_{\hatt\chi})
$
of finite length modules for
the completion
$U^\la_{\hatt\chi}$ of $U_\la$ at $\chi$.

According to
\ref{The center}
the specialization $\fZ^\la$ of  the center $\fZ$ of $U$
is the space of functions on
$
\fg^*\tw_\la
\df
(\fg^*\tw\tim_{\fh^*\tw//W}\fh^*//W)\tim_{\fh^*//W}\la
=\
\fg^*\tw\tim_{\fh^*\tw//W}AS(\la)
$.
For instance, any integral $\la$ is killed by the Artin-Schreier map,
so
$
\fg^*\tw_\la
=
\NN\tw
$
and
$U^\la$ is an $\OO(\NN\tw)$-algebra.

\sss{
(Generalized) Springer fibers
}
\lab{Springer fibers}

Fix  $(\chi,\nu)\in {\fg^*\tw}\times_{{\fh^*\tw}//W} \fh^*$, and define
$\BB_\chi,\BB_{\chi,\nu}\subset \gtil$ by
 $\BB_\chi=\PR_1^{-1}(\chi)$, $\BB_{\chi,\nu}=\PR^{-1}(\chi,\nu)$
(notations of \ref{redu});
we equip $\BB_\chi$, $\BB_{\chi,\nu}$
 with the {\em reduced}\footnote{``Reduced'' will only be used in
        lemma \ref{Springer Slodowy}c. It is irrelevant in \S4 and \S5
        since we only use formal neighborhoods of the fiber.}
 subscheme structure.  When $\chi$ is nilpotent (so $\nu=0$
and $\BB_{\chi,\nu}=\BB_{\chi}$)
 it is called a Springer fiber; otherwise
we call it a generalized Springer fiber.

One can show  that  $\BB_{\chi,\nu}$ is connected; in fact
it is a Springer fiber for the centralizer of $\chi_{ss}$ where
$\chi=\chi_{ss}+\chi_{nil}$ is the Jordan decomposition.
 Thus $\BB_{\chi,\nu}$ is a
connected component of $\BB_\chi$.
Via the projection
$ \gtil @>\pi>> \BB$ the (generalized) Springer fiber can be identified
with a subscheme
$\pi(\BB_{\chi,\nu})$ of $\BB$,
and
$\BB_{\chi,\nu}$ is a section
of
$\gtil$ over
$\pi(\BB_{\chi,\nu})$.

\lemm
\lab{small equivalences}
If $\la\in \fh^*$ is regular and
$(\chi,AS(\la))\in {\fg^*\tw}\times_{{\fh^*\tw}//W}\fh^*\tw$,
the  equivalences in Theorem
\ref{Localizationthm}
restrict
to
$$
\Db(mod^c_{\chi}(\DD^\la))
\cong\
\Db(mod^{fg}_{\chi}(U^{\la}))
,\ \ \
\Db(mod^c_{ {\la,\chi}}(\tii \DD))
\cong\
\Db(mod^{fg}_{ {\la,\chi}}( U))
.$$

\pf
$\OO(\fg^*\tw)$ acts on the categories $mod^c(\tii\DD)$,
$mod^{fg}(U)$
etc,
and on their derived categories.
The  equivalences in Theorem
\ref{Localizationthm}
are equivariant under $\OO(\fg^*\tw)$
and therefore they restrict to the full subcategories of
objects on which the $p$-center acts by the
generalized character $\chi$
(cf. Lemma \ref{Derived categories of sheaves supported on a subscheme}).

\cor
\lab{KComparison}
If $\la$ is regular and
$(\chi,AS(\la))\in{\fg^*\tw}\times_{{\fh^*\tw}//W}{\fh^*\tw}$,
localization gives a canonical isomorphism\
$K(U^\la_\chi)\cong\ K(\DD_\chi^\la)$.

\pf
By Lemma
\ref{small equivalences},
localization gives
isomorphism
$K(\Db(mod^{fg}_{\hattt \chi}(U^\la)))
\con
K(\Db(mod^c_{\hattt \chi}(\DD^\la)))
$. This simplifies to the desired isomorphism since
$$
K(U^\la_\chi)
\df\
K(mod^{fg}(U^\la_\chi))
\con\
K(mod^{fg}_{\hattt \chi}(U^\la))
\cong\
K(\Db(mod^{fg}_{\hattt \chi}(U^{\la})))
,$$
the first isomorphism
is the fact that
the  subcategory
$mod^{fg}(U^\la_\chi)
$
generates
$
mod^{fg}_\chi(U^\la)
$
under extensions,
 and the second
is the equality of $K$-groups of a triangulated category
(with a bounded $t$-structure), and of its heart.
Similarly,
$$
K(\DD_{\chi}^{\la})
\df\
K(mod^c(\DD_\chi^{\la}))
\con\
K(mod^c_{\hattt \chi}(\DD^{\la}))
=\
K(\Db(mod^c_{\hattt \chi}(\DD^{\la}))).
\ \ \ \ \epf
$$

\se{\bf
Splitting of the Azumaya algebra of crystalline differential operators
on  (generalized) Springer fibers
}
\lab{Splitting on  Springer fibers}

\sus{
$\DD$-modules and coherent sheaves
}
Since $\tii\DD$ is an Azumaya algebra  over
$
\tii T^*\BB\tw\tim_{\fh^*\tw}\fh^*
$, for $\la\in\fh^*$
we will view $\DD^\la$ as an Azumaya algebra  over
$T^*_\nu\BB\tw$ where
\ud{$\nu=\AS(\la)$}
(see \ref{Torsors}).
The aim of this section is the following

\theo
\lab{Achi}
a) For any $\la\in\fh^*$, Azumaya algebra
$\tii\DD
$ splits on
the formal neighborhood
in
$
\tii T^*\BB\tw\tim_{\fh^*\tw}\fh^*
$
of
$
\BB_{\chi}\tw\times_{\fh^*\tw}\la
\cong
\BB_{\chi,\nu}\tw
$,
i.e., there
is a vector bundle $\MM^\la_\chi$ on this formal neighborhood,
such that
the restriction of $\tii\DD$ to the neighborhood
is isomorphic to
$\EEnd_\OO(\MM_\chi^\la)$.

b) The functor $\F\mapsto  \M^\la_\chi\otimes_\O \F$ provides
equivalences
$$
\Coh_{\BB_{\chi,\nu}\tw}(\tii T^*\BB\tw\tim_{\fh^*\tw}\fh^*)
\con\
mod^c_{\chi,\la}(\tii\DD)
,$$
$$
\Coh_{\BB_{\chi,\nu}\tw}(T^*_\nu\BB\tw)
\con\
mod^c_\chi(\DD^\la)
.$$

\pf
(b) follows from (a).
Lemma \ref{sdvig} shows that to check statement (a) for particular $(\chi,\la)$
it suffices to check it for $(\chi,\la+d\eta)$ for some character $\eta:H\to
\Gm$.

Let us say that $\la\in \fh^*$ is {\it unramified} if for any
coroot $\alpha$
we have either $\langle \alpha, \la+\rho\rangle=0$, or
 $\langle \alpha, \la\rangle \not \in \Fp$. We claim that
 for any  $\la\in \fh^*$ one
can find a character $\eta:H\to \Gm$ such that $\la+d\eta$ is unramified.
For this it suffices to show the existence of $\mu \in \h^*(\Fp )$, such that
$\langle \alpha, \la+\rho \rangle = \langle \alpha, \mu\rangle$ for any coroot $\alpha$,
such that $\langle \alpha, \la\rangle\in \Fp$. These conditions constitute a system
of linear equations over $\Fp$, which have a solution over the bigger field $\k$. By
standard linear algebra
they also have a solution over $\Fp$.

Thus it suffices to check (a) when $\la$ is unramified.
The next proposition shows that for unramified $\la$ the
restriction of
$\tii\DD$ to the formal neighborhood of
${\BB_{\chi}\tw}\times_{\fh^*\tw}\la$
is isomorphic to the
pull-back of an Azumaya algebra on
the formal neighborhood
$\hatt\chi\tw=FN_{\fg^*}(\chi)\tw$ of $\chi$ in $\fg^*\tw$. The latter splits
by \cite[IV.1.7]{MI}  (vanishing of the Brauer group of a complete local ring
with separably closed residue field). \epf

\sus{
Unramified Harish-Chandra characters
}
\lab{Steinberg}
Let $\fh^*_{unr} \subset \fh^*$ be the open set of all unramified weights.
Let $\fZ_{unr}$ be the algebra of functions on $\fg^*\tw\times _{\fh^*\tw//W}
\fh^*_{unr}
\sub
Spec(\fZ)$ (see \ref{The center}).

\pro\lab{AU}
a) $U\otimes_{\fZ} \fZ_{unr}$ is an Azumaya algebra over $\fZ_{unr}$.

b) The action map $U\otimes_{\fZ}\O(\gtil\tw\times_{\fh^*\tw} \fh^*)\to
\tii \DD$ induces an isomorphism
 $$
U\otimes_{\fZ}\O(\gtil\tw\times_{\fh^*\tw} \fh^*_{unr})
\
\iso
\
\tii \DD|_{\gtil\tw\times_{\fh^*\tw} \fh^*_{unr}}
.$$


\pf (a) is proved in \cite{BG}, Corollary 3.11; moreover, it is shown
in {\it loc. cit.} that for $\fz\in \fZ_{unr}$ and a baby Verma module
$M$ with central character $\fz$ we have an isomorphism
$U(\fg)\otimes_{\fZ}\k_\fz \iso End_\k(M)$. This implies (b)
in view of  Proposition \ref{baby}. \epf

\rems\lab{Verma} 1) Consider the restriction of $\MM_\chi^0$
to the reduced subscheme $\BB_\chi\tw$.
 In view of Remark \ref{spl_on_Z} it defines (and is defined by)
a line bundle with a flat
 connection on the Frobenius neighborhood of $\BB_\chi$ in $\BB$. The
 requirement that  the sheaf on $T^*X\tw$ arising from the bundle with
connection lives on $\BB_\chi\tw$ is equivalent to the equality between the
$p$-curvature of the connection and  the section of
$\Omega^1_\BB|_{\BB_\chi}$ defined by $\chi$
(cf. Remark \ref{p_curv}).\footnote{As pointed out in Remark \ref{p_curv}
the $p$-curvature of a $\DD_X$-module $\EE$ is a parallel
section of $Fr^*(\Omega^1)
\otimes \EEnd(\EE)$. If $\EE$ is a line bundle we get a parallel section of
$Fr^*(\Omega^1)$, i.e. a section of $\Omega^1$; for a line bundle with
a flat connection on $FrN_X(Y)$ its $p$-curvature is a section of
 $\Omega^1_X|_Y$.}

For some particular cases such a line bundle with a flat
connection was constructed in \cite{MR}.
Notice that already in the case $G=SL(3)$, and $\chi$
subregular this line bundle is non-trivial for any choice of the
splitting bundle $\MM_\chi^\lambda$ (see, however,  equality
\eqref{MM prime} in the proof of Lemma \ref{class of E} below).

2) The choice of  a character $\eta\in \Lambda$ such
that $\lambda +d\eta$ is unramified, provides a particular splitting
line bundle $\MM_\chi^\lambda=\MM_\chi^\lambda(\eta)$ in Theorem
\ref{Achi}(a): apply the equivalence of Lemma  \ref{sdvig} to the
trivial (equivalently, lifted from $\hatt{\nu}\tw$) splitting vector
bundle on the formal neighborhood of
 $\BB_\chi\tw\times _{\fh^*\tw}(\lambda+d\eta)$. It is easy to see
then that
$\MM_\chi^\lambda(\eta+p\zeta)=\MM_\chi^\lambda(\eta)\otimes \OO_{-\zeta}$.

3) One can show that the Azumaya algebra  $U\otimes_{\fZ} \fZ_{unr}$
splits on some closed subvarieties of $Spec( \fZ_{unr})$. E.g.
 the Verma module
$M^\fb(-\rho)\df\ ind_{U\fb}^{U\fg}\ \k_{-\rho}$  is easily seen
to be a splitting
module
on
$\fn\times \{-\rho\}
$.

\sus{
$\fg$-modules and coherent sheaves
}
By putting together known equivalences
(Theorem \ref{small equivalences} and Theorem \ref{Achi}(b)),
we get

\theo
\lab{better}
If $\la\in\fh^*$ is regular
and
$(\chi,\la)\in \fg^*\tw\times_{\fh^*\tw//W}\fh^*$ with $(\chi,W\bullet \la)\in Spec(\fZ)$, then
there are 
equivalences  (set $\nu=\AS(\la)$)
$$
\Db(mod^{fg}_{\hattt\chi}(U^\la))
\cong\
\Db(mod^c_{\hattt\chi}(\DD^\la))
\cong\
\Db(\CC oh_{\BB_{\chi,\nu}\tw}(T^*_\nu\BB\tw))\
;$$
$$
\Db(mod^{fg}_{\hattt (\la,\chi)}(U))
\cong\
\Db(mod^{c}_{\hattt (\la,\chi)}(\tii \DD))
\cong
\Db(\CC oh_{\BB_{\chi,\nu}\tw}(\tii T^*\BB\tw\tim_{\fh^*\tw}\fh^*))\
.$$

\rem The equivalences depend on the choice of the splitting bundle
$\MM_\chi^\lambda$ in Theorem \ref{Achi}(a); thus on the choice of
 $\eta\in \Lambda$
such that $\lambda+d\eta$ is unramified (see Remark \ref{Verma}(2)).
Replacing $\eta$ by $\eta+p\zeta$ we get another equivalence, which is
the composition of the first one with twist by $\OO_\zeta$.

\exs \label{exaSL3} Let us list some  objects
in $mod_\chi^{fg}(U^\lambda)$ whose image in the derived category of
coherent sheaves can be computed explicitly. We leave the proofs as an
exercise to the reader.

0) A baby Verma module $M_{\b,\chi;\lambda+2\rho}$ corresponds to
a sky-scraper sheaf, see section \ref{baby}.

\medskip

For simplicity of notation, in the next two examples we set $\lambda=0$,
and  normalize the equivalences by setting $\eta=(p-1)\rho$; notice
that for $\chi=0$ this choice gives the splitting on the zero section $\BB_0$
from \ref{Lagrangian}.

1) Let $G$ be simple and simply-laced, and $\chi$ a subregular nilpotent.
 Recall that the irreducible components of the
(reduced) Springer fiber are indexed by the simple roots of $G$, each
component is a projective line.
The images
of irreducible objects in $mod_\chi(U^0)$ are as follows:
$\O_{\Pone_\alpha}(-1)[1]$; and $\O_{\pi^{-1}(\chi)}$. Here
$\Pone_\alpha$ runs over the set of irreducible components of $\B_\chi\tw$,
$\pi:T^*\B\tw\to \N\tw$ is the projection, and $\pi^{-1}$ stands for
 the
scheme-theoretic preimage. Notice that the same objects appear in the
geometric theory of McKay correspondence, \cite{KV}.

2) $G=SL(3)$, $\chi=0$. See the appendix for a description of this example.

%


\sus{
Equivalences on formal neighborhoods
}
\lab{Equivalences on formal neighborhoods}
We will extend Theorem
\ref{better} to the formal neighborhood
of $\chi$.\footnote{
The same argument gives
extension  to the formal neighborhood of $\la$.
}
For  $\la,\chi,\nu$ as in
\ref{better},
denote by $\hatt\chi$ and $\hatt{\BB_{\chi,\nu}}$
the formal neighborhoods of $\chi$ in
$\PR_1(T^*_\nu\BB)$ and $\BB_{\chi,\nu}$ in $T^*_\nu\BB$.

\theo
\lab{even better}
There are canonical equivalences
$
\Db_{fg}( U^\la_{\hatt\chi}         )
\cong\
\Db_c(  \DD^\la_{\hatt\chi}         )
\cong\
\Db_c( \OO_{\hatt{\BB_{\chi,\nu}\tw}}   )\
$.

\pf
Our main reference for sheaves on a formal scheme $\fX$
is \cite{AJL}. We consider
the full subcategory
$D_c^b(\OO_\fX)$ of the derived category
$D(\OO_\fX)$ of the abelian category of {\em all}
$\OO_\fX$-modules by requiring that cohomology sheaves are coherent
(and almost all vanish).
Denote by    $U^\la_{\hatt\chi},\DD^\la_{\hatt\chi}$
the restrictions of the coherent $\OO$-algebras
$U^\la,\DD^\la$ to $\hatt\chi,\hatt{\BB_{\chi,\nu}}$.
Now, (coherent) $\DD^\la_{\hatt\chi}$-modules
are
(coherent)
$\OO_{ \hatt{\BB_{\chi,\nu}} }$-modules with extra structure,
and this allows us to lift the direct image
functor
$
R\mu_*:\Db_c( \OO_{\hatt{\BB_{\chi,\nu}\tw}}    )\
\to
\Db_c( \OO_{\hatt\chi}  )\
$
to
$
R\mu_*:
\Db_c( \DD^\la_{\hatt\chi}  )\
\to
\Db_c( U^\la_{\hatt\chi}    )\
$ (as in  \ref{globase}).
The proof that this is an equivalence follows the proof of Theorem
\ref{Localizationthm}. First, $R\mu_*
(\DD^\la_{\hatt\chi}    )
\cong
U^\la_{\hatt\chi}
$ follows from
\ref{Dcohomology} by the formal base change for proper maps
(\cite{EGA}, Theorem 4.1.5.). Then, for the Calabi-Yau trick
(\ref{CY}) one uses the Grothendieck duality for formal schemes
(\cite{AJL}, Theorem 8.4, Proposition 2.5.11.c and 2.4.2.2).
The second equivalence follows from Theorem
\ref{Achi} above. \epf

\sss{}
In the remainder of the section, for simplicity,
{\em $\la$ is integral regular and $\chi\in \N$.}

\cor \lab{count} For $p> h$ there is a natural isomorphism of
Grothendieck groups $ K(U^\la_\chi)\cong K(\BB_\chi\tw) $. In
particular, the number of irreducible $U^\la_\chi$-modules is the
rank of $K(\BB_{\chi})$. (This rank is calculated below in Theorem
\ref{K theorem}.)

\pf It is well known that  for a closed embedding
$\iota:\fX\imbed \fY$ of Noetherian schemes we have an isomorphism
 $K(\fX)\iso K(\Coh_\fX(\fY))$ induced by the functor $\iota_*$.
In particular, $K(\BB_\chi\tw)
\cong
K(\CC oh_{\BB_{\chi}\tw}(T^*\BB\tw))
\cong
K(\CC oh_{\BB_{\chi}\tw}(\tii T^*\BB\tw\tim_{\fh^*\tw}\fh^*))$. \epf

\rems
(a)
In the case when $\chi$ is regular nilpotent
in a Levi factor the corollary is a fundamental
observation  of Friedlander and Parshall
(\cite{FP}). The general case was conjectured
by Lusztig (\cite{Lu1},\cite{Lu}).

(b)
Theorem \ref{Achi} provides a natural  isomorphism of
K-groups. However, if one is  only interested in
the number of irreducible modules 
(i.e., the size of the K-group),
one does not
need the splitting (i.e., section \ref{Splitting on  Springer
fibers}). Indeed, one can show that for any Noetherian scheme $X$,
and an Azumaya algebra $\AA$ over $X$ of rank $d^2$, the forgetful
functor from the category of $\AA$-modules to the category of
coherent sheaves induces an isomorphism
$K(\AA-mod)\otimes_{\Zet}\Zet[\frac{1}{d}]\iso
K(\Coh(X))\otimes_{\Zet}\Zet[\frac{1}{d}]$.

\sus{Equivariance}
Let $H$ be a group. An $H$-category\footnote{The term  ``a weak $H$-category''
would be more appropriate here, since we do not
fix isomorphisms between $[gh]$ and $[g]\circ [h]$; we use the shorter
expression, since the more rigid structure does not appear in this
paper.}
is a category $\mathcal C$ with
functors
$[g]: {\mathcal C} \ra {\mathcal C}, \ g\in H$,
such that $[e_H]$ is isomorphic to the
identity functor,
and $[gh]$ to $[g]\circ [h]$ for  $g,h\in H$.
If ${\mathcal C}$ is  abelian or triangulated $H$-category
we ask that the functors
$[g]$ preserve the additional structure,
and
then $K({\mathcal C})$ is an $H$-module.
An $H$-functor is a functor
${\mathcal F}: {\mathcal C} \ra {\mathcal C}^\prime$
between $H$-categories such that
$[g]\circ {\mathcal F} \cong  {\mathcal F} \circ [g]$
for $ g\in H$.
If it induces a map of $K$-groups
$K({\mathcal F}): K({\mathcal C}) \ra ({\mathcal C}^\prime)$,
then  this  is a map  of $H$-modules.

The actions of the group
$G(\k)$ on $U$ and $\BB$
make all categories
in the localization
theorem \ref{Localizationthm}
into $G(\k)$-categories, while the categories
appearing in Theorem \ref{Achi}(b) (for $\nu=0$)
are $G_\chi(\k)$ categories.
 Standard methods show that
the action of
$G_\chi(\k)$ on these K-groups factors through $A_\chi=\pi_0(G_\chi)$.

\pro
\lab{equivariance}
The isomorphism $K(U^\la_\chi)\cong K(\BB_\chi\tw)$
in Corollary
\ref{count}
is an isomorphism of
$A_\chi$-modules.

\pf The functors
$\Rr\Gamma_{\DD^\la}$
and
$\Rr\Gamma_{\tii\DD,\la}$
are clearly $G(\k)$-functors.
Thus it suffices to check that
the Morita equivalences in Theorem \ref{Achi}
are $G_\chi(\k)$-functors.

We will use  a general observation that if a group $H$
acts on a split Azumaya algebra $A$ with a center $Z$
and a splitting module $E$ is $H$-invariant
(in the sense that $^gE\cong E$ for any $g\in H$),
then the Morita equivalence defined by $E$ is an $H$-functor.
Indeed, for $g\in H$ a choice of
an $A$-isomorphism $\psi_g:\ ^gE\con E$ gives for each $A$-module
$M$ a $Z$-isomorphism
$$
^g(E\otimes_AM)
\stackrel{Id}{\longrightarrow}
\;^gE\otimes_A(^gM)
\stackrel{\psi_g \otimes Id}{\longrightarrow}
E\otimes_A(^gM).
$$
Thus we have to check that  the splitting bundle $\MM_\chi^\lambda$
 of Theorem \ref{Achi}
is  $G_\chi(\k)$ invariant.
The equivalence between the Azumaya algebras $\DD^\la$ and $\DD^{\la+d\eta}$
 from Lemma \ref{sdvig} is clearly $G(\k)$, and hence $G_\chi(\k)$
equivariant. Then our Azumaya algebra is  $G_\chi(\k)$ equivariantly
identified with  the pull-back of an Azumaya algebra on
 $\hatt\chi \tw$ (see the proof
of Theorem \ref{Achi}), and  $\MM_\chi^\la$ is the pull-back of a splitting
bundle from $\hatt \chi\tw$; thus it is enough to see that the latter
is $G_\chi(\k)$ invariant. This is obvious, since any two vector bundles
(and also any two modules over a given Azumaya algebra) on $\hatt \chi \tw$
 of a given rank are isomorphic. \epf

\rems
(1) Proposition
\ref{equivariance}
can be used to sort out how many simple
modules in a regular block are twists of each other,
a question raised by Jantzen (\cite{Ja3}).
For instance, if $G$ is of type $G_2$
and
$p> 6$,
we find that 3 out of 5 simple modules in a regular
block are twists of each other.

(2) We expect that Proposition \ref{equivariance} can be strengthened:
the splitting bundle $\MM_\chi^\lambda$ can be shown to carry
 a natural $G_\chi(\k)$
equivariant structure, thus
the equivalences of Theorem \ref{Achi}(b) can be enhanced
to equivalences of strong $G_\chi(\k)$ categories
(the isomorphisms $[gh]\cong [g]\circ [h]$
are fixed and satisfy natural compatibilities).
We neither prove nor use this fact here.

\se{\bf Translation functors and dimension of $U_\chi$-modules }
In this section we spell out compatibility between the
localization functor and translation functors, and use our results
to derive some rough information about the dimension of
$U_\chi$-modules for $\chi\in\NN$. So we consider only integral
elements of $\fh^*$ and these we view as differentials of elements of $\La$.
Similar methods can be applied to computation of
the characters of the maximal torus in the centralizer of $\chi$
acting on an irreducible
$U_\chi$-module. We keep the assumption
$p> h$.

\sus{ Translation functors } For $\la\in\La$, $\DD^\la\df\
^{\OO_\la}\DD$ is canonically isomorphic to $\DD^{d\la}$ for the
differential $d\la$ and we also denote  $U^\la\df
U^{d\la}$ etc. We denote by $M@>>>[M]_\la$ the  projection of the
category of finitely generated
$\fg$-modules with a locally finite action of $\zhc$
to its direct summand $ mod^{fg}_\la(U) \df\ mod^{fg}_{d\la}(U) $.
For $\la,\mu\in\La$ the translation functor $T_\la^\mu:
mod^{fg}_\la(U) @>>> mod^{fg}_\mu(U) $ is defined by $
T_\la^\mu(M)\df\ [V_{\mu-\la}\ten M]_\mu$
where $V_{\mu-\la}$ is the standard $G$-module with an extremal
weight $\mu-\la$ as defined in 3.1.1.

Notice that the translation functor is well-defined. First,
$V_{\mu-\la}\ten M$ is finitely generated by \cite[Proposition
3.3]{Ko}. To  show that the action of $\zhc$ on $V_{\mu-\la}\ten
M$ is locally finite we can assume that $M$ is annihilated by a
maximal ideal  $I_\eta$ of $\zhc$.
By \cite[Theorem 1]{MR1},
for a very good $p$ there is a
ring homomorphism $\Upsilon : {\mathfrak Z}_{\mathbb
Z}\rightarrow \zhc ={\mathfrak Z}_{\mathbb Z}\otimes_{\mathbb
Z}\k$ where ${\mathfrak Z}_{\mathbb Z}$ is the center of
$U({\mathfrak g}_{\mathbb Z})$.
By
\cite[Theorem 5.1]{Ko}, for each $x\in\mbox{im}(\Upsilon)$,
on
$V_{\mu-\la}\ten M$
\begin{equation}
\lab{kostant} \prod_{\nu} (x-\eta(x)-\nu(x))=0,
\end{equation}
where $\nu$ run over the weights of $V_{\mu-\la}$. Thus $\zhc$ is
spanned by elements satisfying equation (\ref{kostant}). It
follows that the action of $\zhc$ on $V_{\mu-\la}\ten M$ is
locally finite.


We review some standard ideas.
For $\la,\mu,\eta\in\La$ we denote by
$\WW_\eta$ the weights of $V_\eta$ and $\WW_\la^\mu \df\
(\la+\WW_{\mu-\la}) \cap\ W_\aff'\bu\mu $.
Since we assume  $p> h$,\ \
$\WW_\la^\mu =\ (\la+\WW_{\mu-\la}) \cap\ W_\aff\bu\mu $.

\sss{} \lab{Translation of cohomology} For $\MM\in
\Db(mod^c_\la(\tii\DD))$,
the sheaf of $\fg$-modules
$ V_{\eta}\ten
\MM=(V_{\eta}\ten\OO)\ten_\OO \MM $ is an extension of terms
$V_{\eta}(\nu)\ten(\OO_\nu\ten_\OO\MM)$ where $\nu$
runs over the set of weights $\WW_{\eta}$  and $V_{\eta}(\nu)$
is the corresponding weight space.
Since $\OO_\nu\ten_\OO\MM\in \Db(mod^c_\la(\tii\DD))$ we get
the local finiteness of the
$\zhc$-action on the sheaf 
$V_{\eta}\ten \MM$.
Therefore,  translation functors  commute in a sense with taking
the cohomology of $\DD$-modules
$$
T^\mu_\la(\Rr\Ga_{\tii\DD,\la}\MM) =\ [V_{\mu-\la}\ten
\Rr\Ga_{\tii\DD,\la}\MM]_\mu =\ [\Rr\Ga_{\tii\DD}(V_{\mu-\la}\ten
\MM)]_\mu \cong\ \Rr\Ga_{\tii\DD,\mu}([V_{\mu-\la}\ten \MM)]_\mu)
.$$ Moreover, $ [V_{\mu-\la}\ten_\OO\MM]_\mu $ is a successive
extension of terms $V_{\mu-\la}(\nu)\ten(\OO_\nu\ten_\OO\MM)$ for
weights $\nu\in\WW_\la^\mu-\la\sub \WW_{\la-\mu}$. There are two
simple special cases:

\lemm \lab{Up and down} Let $\la,\mu$ lie in the same closed
alcove $\AA$.

(a) (``Down''.) If $\mu$ is in the closure of the facet of $\la$
then
$$
T^\mu_\la(\Rr\Ga_{\tii\DD,\la}\MM) \ \cong\ \Rr\Ga_{\tii\DD,\mu}
(\OO_{\mu-\la}\ten_\OO\MM) .$$

(b) (``Up''.) Let  $\la$ lie in the single wall $H$ of $\AA$ and
$\mu$ be regular. If $s_H(\mu)<\mu$ for the reflection $s_H$ in
the $H$-wall, then
$$
\Rr\Ga_{\tii\DD,s_H(\mu)} (\OO_{\la-\mu}\ten_\OO\MM) @>>>
T^\mu_\la(\Rr\Ga_{\tii\DD,\la}\MM)
@>>> \Rr\Ga_{\tii\DD,\mu} (\OO_{\mu-\la}\ten_\OO\MM) .$$

\pf This follows from \ref{Translation of cohomology} and the
following combinatorial observation from \cite[Lemmas 7.7 and
7.8]{Ja0}:\ \ if $\la,\mu\in\La$ lie in the same alcove then
$$
\WW_\la^\mu =\ (\la+\WW_{\mu-\la}) \cap\ W_\aff\bu\mu \ =\
(W_\aff)_\la\bu\mu \ \sub\ \la+W\cd(\mu-\la) .$$ Indeed, the
assumption  in (a) implies that
$(W_\aff)_\mu\sub (W_\aff)_\la$, hence
$\WW_\la^\mu=\{\mu\}$, while in (b) we assume
$(W_\aff)_\la=\{1,s_H\}$, hence $\WW_\la^\mu=\{\mu,s_H(\mu)\}$, and
$s_H(\mu)$ appears earlier in the filtration since $s_H(\mu)<\mu$.
\epf

\sus{Dimension} We set $R=\prod\limits_{\alpha} \langle
\rho,\ch\alpha\rangle$ where $\alpha$ runs over the set of
positive roots of $G$.

\theo
\lab{dimdim}
 Fix $\chi\in \N$ and a regular $\lambda\in \Lambda$.
For any  module
$M\in mod^{fg}_{(\lambda,\chi)}(U)$
there exists a polynomial $\bd_M\in \frac{1}{R}\Zet
[\Lambda^*]$ of degree $\le \dim(\BB_\chi)$,
such that for any $\mu\in \Lambda$ in the closure of the
alcove of $\lambda$, we have
$$\dim(T_\lambda^\mu(M))=\bd_M(\mu).$$
Moreover,
$\bd_M(\mu)=p^{\dim \BB}\bd_M^0(\frac{\mu+\rho}{p})$
for
another
polynomial $\bd_M^0\in\frac{1}{R}\Zet  [\Lambda^*]$,
such that
$\bd_M^0(\mu)\in \Zet$ for $\mu \in \Lambda$.

\rems
\lab{lifts to regular or dimsum}
(0) The theorem is suggested by the experimental data kindly
provided by J.~Humphreys and V.~Ostrik.

(1) The proof of the theorem gives an explicit description of
$\bd_M$ in terms of the corresponding coherent sheaf $\FF_M$
on $\BB_\chi\tw$.

(2) For $\mu$ and $\la$ as above, any module
$N\in mod^{fg}_{(\mu,\chi)}(U)$ is of the form
$T^\mu_\la M$ for some $M\in mod^{fg}_{(\lambda,\chi)}(U)
$.\footnote{Also, exactness of $T^\mu_\la$ implies that if
$N$ is irreducible we can choose
$M$ to be irreducible.}
Indeed, according to Lemma
\ref{Up and down}.a and
Proposition \ref{compid}.c,
$T^\mu_\la\RGa(\OO_{\la - \mu}\ten\LL^\mu N)
=N$.
Since $T^\mu_\la$ is exact we can choose $M$ as the
zero cohomology
of $\RGa(\OO_{\la-\mu}\ten\LL^\mu N)$.

\cor
\lab{Kac Weisfeiler}
The dimension of any  $N\in mod^{fg}_\chi(U)$
is divisible by $p^{\codim_\BB \BB_\chi}$.

\pf
To apply the theorem observe that $\dim(N)<\yy$, so
we may assume that $\fZ_{HC}$ acts by a generalized
eigencharacter. Since $\chi\in\NN$ eigencharacter is
necessarily integral,
so it lifts to some
$\mu\in\La$.
We choose a regular $\la$ so that $\mu$ is
in the closure of the $\la$-facet,
and $M\in mod^{fg}_{(\lambda,\chi)}(U)
$
as in the remark
\ref{lifts to regular or dimsum}(2).
Then Theorem
\ref{dimdim}
says that
$\dim(N)
=p^{\dim \BB}\ \cd\ \bd_M^0(\fra{\mu+\rho}{p})
$.
For $\de=\deg(d_M^0)=\deg(d_M)\le \dim(\BB_\chi)$,\
rational number
$\dim(N)/p^{\dim(\BB)-\de}=\
p^\de\cd\bd_M^0(\fra{\mu+\rho}{p})
$
is an integer since the denominator divides
both $R$ and a power of $p$, but
$R$ is prime to $p$ for $p> h$
(for any root $\al$,\
$\lb\rho,\ch \al\rb<h$).
\epf

\rem
The statement of the corollary was conjectured by Kac and
Weisfeiler  \cite{KW},
 and proved by Premet  \cite{Pr} under less restrictive assumptions on
$p$.
We still found it worthwhile to point out how this
famous fact is related to our methods.

Our basic observation is

\lemm
\lab{class of E}
Let $\MM_\chi^\lambda$
 be the splitting vector bundle for the restriction of the Azumaya
algebra  $\DD^{ \la}$
 to $\BB_\chi\tw$, that was constructed in the proof of Theorem \ref{Achi}.
We have an equality in
$K^0(\BB_\chi\tw)$:
\begin{equation}\lab{inKBchi}
[\MM_\chi^\lambda]=[(Fr_{\BB})_*\O_{p\rho+\lambda}|_{\BB_\chi\tw}].
\end{equation}

\pf
 Since   $\DD^{ \la}$ contains the algebra
of functions on
$\BB\tim_{\BB\tw}
T^*\BB\tw
$,
any $\DD^\la$-module $\F$ can be viewed as a
quasi-coherent sheaf $\FF'$ on $\BB\times _{\BB\tw} T^*\BB\tw$.
If $\F$ is a splitting bundle of the restriction
$\DD^\la\big |_{Z\tw}$ for a  closed subscheme
$Z\subset
T^*\BB$, then $\F'$ is a line bundle on
$\BB\times _{\BB\tw}Z\tw$.
It remains to show that the equality
\begin{equation}
\lab{MM prime}
[(\MM_\chi^\lambda)']=[\O_{p\rho+\lambda}|_{FrN(\BB_\chi)}]
\end{equation}
holds in $K(FrN(\BB_\chi))$.
The construction in the
proof of Theorem \ref{Achi} shows that
$(\MM_\chi^\lambda)'=\O_\lambda\otimes
(\MM_\chi^0)'$, thus it suffices to check
\eqref{MM prime} for one $\lambda$. We will do it for
$\lambda=-\rho$
by constructing a line bundle $\LL$ on
$FrN(\BB_\chi)\times \Aone$
such that the restriction of $\LL$
at $1\in\A^1$ is
$(\MM_\chi^{-\rho})'$, and at $0$ it is
$\OO_{(p-1)\rho}|_{FrN(\BB_\chi)}$;
existence of such a line bundle implies the desired statement by
rational invariance of $K^0$.

Let $\tii\fn\subset T^*\BB$ be the preimage of $\fn\subset \N$ under
the Springer map. Remark \ref{Verma}(3) together with Proposition \ref{AU}(b)
 show that there exists
a splitting bundle
$\tii \MM$ for $\DD^{-\rho}\big|\tii\fn\tw$ whose restriction to
 $\BB_\chi\tw$ is $\MM$;
we thus get a line bundle
 $\tii \MM'$ on
$\BB\times_{\BB\tw}\tii\fn\tw$.
Its  restriction  to the
zero section $\BB\subset \BB\times _{\BB\tw}T^*\BB\tw$  is a line
bundle on $\BB$ whose direct image under Frobenius is isomorphic to
$\OO_\BB^{ p^{\dim \BB}}$. It is easy to see that the only such line
bundle is $\OO_{(p-1)\rho}$. Thus we can let $\LL$
be the
pull-back of $\tii\MM'$ under the map $FrN(\BB_\chi)\times \Aone\to
\BB\times _{\BB\tw}\tii\fn\tw$, $(x,t)\mapsto (x,(Fr(x), t\chi))$.
\epf

We also recall the standard numerics of line bundles on the
flag variety.

\lemm
\lab{nume}
 For any $\F\in \Db(\Coh(\BB))$ there exists
a polynomial $\bd_\F\in \frac{1}{R}\Zet[\Lambda^*]$
such that for $\lambda \in \Lambda$ the
 Euler characteristic of
$R\Gamma(\F\otimes \O_\lambda)$ equals $\bd(\lambda)$.  Moreover,
we have
\begin{equation}\lab{deg of d}
\deg(\bd_\F)\leq \dim \supp (\F);
\end{equation}
\begin{equation}\lab{Fr on dim}
\bd_{Fr^*(\F)}(\mu) = p^{\dim \BB} \bd_\F(\frac{\mu+(1-p)\rho}{p}).
\end{equation}

{\it Proof.} The existence of $\bd_\F$
is well-known,
for line bundles it is given by the Weyl dimension formula, and the general case follows since
the classes of line bundles generate $K(\BB)$.
The degree estimate
follows from Grothendieck-Riemann-Roch once we recall
 that $ch_i(\F)=0$ for $i< \codim\, \supp (\F)$ because
the restriction map $\Hh^{2i}(\BB)\to \Hh^{2i}(\BB-\supp(\F))$ is an isomorphism
for such $i$.
To prove the polynomial identity \eqref{Fr on dim} it suffices
to  check it for $\mu=p\nu-\rho$, $\nu \in \Lambda$. Then it follows
from the well-known
 isomorphism
$Fr_*(\O_{-\rho}) \cong \O_{-\rho}^{\ \pl\ p^{\dim (\BB)}}
$
which
implies that
$$
Fr_*( Fr^*(\F)\otimes \O_{p\nu-\rho})\cong Fr_*(Fr^*(\F\otimes \O_\nu)
\otimes \O_{-\rho}) \cong \F\otimes \O_\nu \otimes Fr_*(\O_{-\rho})
$$
is isomorphic to the sum of $p^{\dim \BB}$ copies of
$
 \FF\otimes \OO_{\nu-\rho}
$.
\epf

\sss{
Proof of the theorem
}
Let
$\F_M\in
\Coh(
\tii\BB_{\hatt\chi}\tw)
$
be the image of $M$ under the equivalence of Theorem \ref{Achi}, i.e.,
$\LL^\la M\cong \MM_\la\ten \FF_M$;
and let $[\FF_M]\in K(\Coh(
\tii\BB_{\hatt\chi}\tw))
=K(\BB_\chi\tw)$ be its class.
According to
Lemma \ref{Up and down}(a)
$$
T_\lambda^\mu (M)
=
\RGa(\OO_{\mu-\la}\ten \LL^\la M)
=
\RGa(\OO_{\mu-\la}\ten \MM_\la \ten \FF_M)
=
\RGa(\MM_\mu
\ten
\FF_M)
.$$
Let
$\int
$
stand for Euler characteristic of
$\RGa$, so that
$
\dim (T_\lambda^\mu (M))
=
\int_{\BB_\chi\tw}\
[\MM_\mu]
\bb{}\cdot
[\FF_M]
$,
where
the multiplication means  the action of $K^0$ on $K$.
Now Lemma \ref{class of E}
lets us rewrite this as
(denote
by $f^*,f_*$ the standard functoriality of
Grothendieck groups and
$\BB_\chi\tw\aa{i}\inj\BB\tw$),
$$
\int_{\BB_\chi\tw}\
i^*[(Fr_\BB)_*\OO_{p\rho+\mu}]
\bb{}\cdot
[\FF_M]
=
\int_{\BB\tw}\
[(Fr_{\BB})_*\OO_{p\rho+\mu}]
\cdot
i_*[\FF_M]
=
\int_\BB
\ \O_{p\rho+\mu}\cdot Fr_\BB^*(i_*[\FF_M])
,$$

So, Lemma \ref{nume}
shows that
$$
\dim(T^\mu_\la M)
=
\bd_{Fr_\BB^*(i_*\FF_M)}(p\rho+\mu)
= p^{\dim \BB}\cd \bd_{\FF_M}(\frac{\mu+\rho}{p}).
$$
Taking into account \eqref{deg of d}, \eqref{Fr on dim}
 we see that the polynomial
$\bd_M^0=\bd_{i_*\FF_M} $ satisfies the conditions
of the theorem. $\square$

\se{\bf
K-theory of Springer fibers
}

In this section we prove Theorem \ref{K theorem}.

\sus{
Bala-Carter classification of nilpotent orbits
\cite{Sp}
}
\lab{Bala-Carter classification of nilpotent orbits}
Let $G_\Z$ (with the Lie algebra $\fg_\Z$)
be the split reductive group scheme over $\Z$ that gives $G$
by extension of scalars: $(G_\Z)_\k=G$.
Fix a split Cartan subgroup
$T_\Z\sub G_\Z$ and a   Bala-Carter datum, i.e.,
a pair
$(L,\la)$ where
$L$ is Levi factor of $G_\Z$ that contains $T_\Z$,
and
$\la$ is a cocharacter of
$T_\Z\cap L'$ (for the derived subgroup
$L'$ of $L$),
such that the
$\la$-weight spaces
$(\fl')^0
$
and $(\fl')^2
$
(in
the Lie algebra $\fl'$ of $L'$),
have the same rank.
To such datum one  associates
for any closed field $k$ of good characteristic
a nilpotent orbit
in $\fg_k$ which we will denote
$\al_k$. It is characterized by: $\al_k$ is dense in $(\fl'_k)^2$.
This gives a bijection between $W$-orbits of Bala-Carter data
and
nilpotent orbits in $\fg_k$. In particular
the classification of
nilpotent orbits over a closed field is uniform
for all   good characteristics (including zero).
This is used  in the formulation of:

\theo
\lab{K theorem}
For
$p> h$ the
Grothendieck group of $\Coh(\BB_{\chi})$ has no torsion and
its rank  coincides with the dimension of
the cohomology of the corresponding   Springer
fiber over a field of characteristic zero.

\sss{
}

The absence of torsion is clear from
Corollary \ref{count}.
The rank will be found   from
known favorable properties of K-theory and
cohomology of Springer fibers
using the Riemann-Roch Theorem.
We start with recalling some standard basic facts about the
K-groups.

\sss{
Specialization in K-theory
}
\lab{specsec}
Let $X$ be a
Noetherian scheme,   flat over a discrete valuation ring $\OO$.
Let   $\eta=Spec(k_\eta), s=Spec(k_s)$ be respectively the
generic and the special point of $Spec(\OO)$
and denote $X_s@>{i_s}>> X@<{i_\eta}<<X_\eta$.
The specialization map
$sp:K(X_\eta) \to K(X_s)$ is defined by
$sp(a)\df\
(i_s)^*(\tii a)$ for $a\in K(X_s)$ and any extension $\tii a\in K(X)$ of
$a$ (i.e. $(i_\eta)^*\tii a=a$).
To see that this makes sense we use
the excision    sequence
$$
K(X_s) \overset{ (i_{s})_*}{\To} K(X)
\overset{(i_\eta)^*}\To K(X_\eta)\to 0
$$
and observe that
$(i_{s})^{ *}(i_{s})_*=0$ on $K(X_s)$ since
the flatness of $X$ gives
exact triangle $\F[1]\to (i_s)^*(i_{s})_* (\F)\to \F$ for $\F\in
\Db(\Coh_{X_s})$.

\sss{
A lift to the formal
neighborhood of $p$
}
\lab{lift}
Assume now that
$\OO$ is the ring of integers
in a  finite extension $\KK=k_\eta$ of $\Q_p$,
with an embedding of the residue field $k_s$ into $\k$.

Let
$G_\OO$
be the group scheme
$(G_\Z)_\OO$ over $\OO$
(extension of scalars), so that $(G_\OO)_\k=G$,
and similarly for the Lie algebras.
By a result of Spaltenstein
\cite{Sp},
one can choose $x_\OO\in\fg_\OO$
so that
(1)
its images
in $\fg_\KK$ and in $\fg_{k_s}$
lie in nilpotent orbits
$\al_\KK$ and $\al_{k_s}$,
(2)
the $\OO$-submodule
$[x_\OO,\fg_\OO]\sub\ \fg_\OO$ has a complementary submodule
 $Z_\OO$,
(3)
for the Bala-Carter  cocharacter
$ G_{m,\Z}@>{\la}>>G_\Z$ (see
\ref{Bala-Carter classification of nilpotent orbits}),
$x_\OO$ has weight $2$
and the sum of all positive weight spaces
$\fg_\OO^{{}^{> 0}}$, lies in
$[x_\OO,\fg_\OO]$.
We denote by
$\BB_{\chi}^\OO$ the
Springer fiber
at $x_\OO$  (i.e., the $\OO$-version of $\BB_\chi$ from
\ref{Springer fibers}),
so it is defined   as the
reduced part of the inverse of $x_\OO$ under the moment map.

\lemm
\lab{Springer Slodowy}
(a)
$Z_\OO$ can be chosen $G_m$-invariant and with weights $\le 0$.

(b) Now
$S_\OO=x_\OO+Z_\OO$
is a slice
to the orbit $\al$ in the sense that: (i)
the conjugation $G_\OO\tim_\OO S_\OO@>>>\fg_\OO$ is smooth,
(ii) the $G_m$-action on $\fg$ by  $c\bu y\df\ c^{-2}\ \cd\ ^{\la(c)}y$,
contracts $S_\OO$ to $x_\OO$.

(c) The  Springer fiber
$X=\BB_{\chi}^\OO$ of $x_\OO$
 is flat\footnote{Though one expects that
 the scheme theoretic fiber is also flat,
this version is good enough for the specialization machinery.} over $\OO$
and
the Slodowy scheme $\tii S_\OO$ (defined as
the preimage of $S_\OO$  under the
Springer map), is  smooth over $\OO$.

\pf
(a) is elementary: if $M\sub A\sub B$ and $M$ has a complement $C$
in $B$ then it has a complement $A\cap C$ in $A$.
Now $[x_\OO,\fg_\OO]$ is $G_m$-invariant
and each weight space
$[x_\OO,\fg_\OO]^n$
has a complement in
$[x_\OO,\fg_\OO]$, hence in
$\fg_\OO$, and then  also a complement
$Z_\OO^n$ in $\fg_\OO^n$. So,
$Z_\OO=
\pl_n\
Z_\OO^n$
is a $G_m$-invariant complement.
Claim
(b$_{\text{ii}}$)
is clear. The smoothness in
(b$_{\text{i}}$)  is valid
on a neighborhood of
$G_\OO\tim_\OO x_\OO$
by (2)
(the image of the differential at a point in
$G_\OO\tim_\OO S_\OO$ is $[x_\OO,\fg_\OO]+Z_\OO$).
Then  the general case follows from the contraction in
(b$_{\text{ii}}$).

In (c), the smoothness of $\tii S_\OO$
follows from
(b$_{\text{i}}$)
by  a formal base change argument
(\cite{Sl}, section 5.3).
Finally, to see that $\BB_{\chi}^\OO$
is flat we use the
cocharacter $\la$ to define a parabolic subgroup $P_\OO\sub G_\OO$ such
that its Lie algebra is
$\fg_\OO^{{}^{\ge 0}}$.
Let
$\BB_{x_\OO}$
be the
scheme theoretic Springer fiber at  $x_\OO$, i.e., the
scheme theoretic inverse of $x_\OO$ under the moment map.
Following  Proposition 3.2 in
\cite{DLP}
we will see that the intersection of
$\BB_{x_\OO}$  with each $P_\OO$ orbit in the flag variety
$\BB_\OO$ is smooth over $\OO$.

Each $w\in W$ defines a Borel subalgebra
$^w\fb_\OO$ of $\fg_\OO$. We view it also as an
$\OO$-point $p^w_\OO$ of the flag variety
$\BB_\OO$ over $\OO$, and use it to generate  a $P_\OO$-orbit
$\fO_w\sub \BB_\OO$.
Consider the maps
$$
\fO_w
@<{\psi_w}<<
P_\OO
@>{\phi}>>
\fg_\OO^{\ge 2}
,
$$
where $\phi$ is given by
$
P_\OO\cong P_\OO\tim_\OO {x_\OO}
@>{}>>
\fg_\OO^{\ge 2}
,\  (g,y)\mm\ ^{g\inv}y
$,\
and
$\psi_w$ by
$
P_\OO\cong P_\OO\tim_\OO {p^w_\OO}
@>{}>>
\fg_\OO^{\ge 2}
,\  (g,p)\mm\ {g}p
$.
Here, $\psi_w$ is smooth as the quotient map
of a group scheme by a smooth group subscheme,
and $\phi$ is smooth since
property (3)
 implies that
 $\fg_\OO^{{}^{\ge 2}}\sub\
[x_\OO,\fg_\OO]^{\ge 2}=\
[x_\OO,\fg_\OO^{{}^{\ge 0}}]
=\
[x_\OO,\fp_\OO]$.
Now,
$
\BB_{x_\OO}
\cap \fO_w$
is  smooth over $\OO$ since the
scheme theoretic inverses
$\psi_w\inv
(\BB_{x_\OO}\cap \fO_w)$
and
$
\phi\inv
(\fg_\OO^{\ge 2}\cap\ ^w\fb_\OO)
$
coincide.

Now we see that any
$p$-torsion function $f$  on an open affine piece $U$
of $\BB_{x_\OO}$ has to be nilpotent
(so the functions on the reduced scheme
$\BB^\OO_\chi$
have  no $p$-torsion and
$\BB^\OO_\chi$ is flat over $\OO$).
The restriction of $f$
to each stratum is zero (strata are
smooth, in particular flat). However any closed point of
$U
$
lies in the restriction  $U_s$ to the special point,
hence in one of the strata.
Since $f$ vanishes at closed points of $U$ it is nilpotent.
\epf

\sss{}
We will use the rational $K$-groups $K(\fX)_\Q\df\ K(\fX)\ten_\Z\Q$
where $\fX$ is a Springer fiber $\BB_{\chi}^A$ over $A$ which could be
$\C,\OO,\eta,s,\k$ etc.
The main claim in this section is

\pro
\lab{spe}
Assume that
$\pl_i\ \Hh_{et}^{2i}(\BB_{\chi}^{\barr\KK},\Ql(-i))$ is a trivial
 $Gal(\overline{\KK}/\KK)$ module\footnote{
A finite extension $\KK/\Qu_p$
satisfying this assumption exists by Lemma \ref{Tate}.}.

(a)
The specialization
 $sp:
K(\BB_{\chi}^\eta)_\Qu
\con
K(\BB_{\chi}^s)_\Qu
$
identifies the $K$-groups over generic and
special points.

(b)
The base change map
identifies the $K$-groups  over the
special point  and  over $\k$.
Also, for any embedding
$
\KK
\imbed \Ce
$
the corresponding base
change maps identifies
$K$-groups over
the generic point and over  $\C$:
\begin{equation}
\label{etaCe}
K(\BB_{\chi}^\eta)_\Qu\iso K(\BB_{\chi}^\Ce)_ \Qu,
\end{equation}
\begin{equation}\label{stok}
K(\BB_{\chi}^s)_ \Qu\iso K(\BB_{\chi}^\k)_ \Qu.
\end{equation}

\sss{
Proposition \ref{spe} implies
Theorem \ref{K theorem}
}
\lab{implies}
In the  chain of isomorphisms
$$
K(\BB_{\chi}^\k)_ \Qu
\conn
K(\BB_{\chi}^s)_ \Qu
\bb{sp}\conn
K(\BB_{\chi}^\eta)_\Qu
\con
K(\BB_{\chi}^\Ce)_ \Qu
\bb{\tau}\con
A_\bu(\BB_{\chi}^\Ce)_\Qu
\cong
\Hh^*(\BB_{\chi}^\Ce,\Q)
,$$
the first three are
provided by the proposition.
It is shown in \cite{DLP} that the Chow group
$A_\bu(\BB_{\chi}^\Ce)$ is a free
abelian group of finite rank equal to
$\dim \Hh^*(\BB_{\chi}^\Ce,\Q)$.
Finally,
by \cite{Fu}, Corollary 18.3.2, the
``modified Chern character'' $\tau_{\BB_{\chi}^\C}$ provides
the fourth isomorphism.
\epf


\sus{
Base change from $\KK$ to $\C$
}
The remainder is devoted to the proof of
Proposition \ref{spe}.
We  need two  standard auxiliary lemmas on Galois action.

\lemm
\label{KL}
Let $L/K$ be a  field extension. Let $X$ be a scheme of finite
type over $K$. Then the base change map
$bc=bc_K^L:K(X)_\Qu \to
K(X_L)_ \Qu$  is injective. If
$L/K$ is a composition of a purely transcendental and a normal
algebraic extension (e.g. if $L$ is algebraically closed) then the
image of $bc$ is the space of invariants $K(X_L)_\Qu^{Gal(L/K)}$.

\pf
If $L/K$ is a finite normal extension, then the
direct image (restriction of scalars) functor induces a map
$res:K(X_L)\to K(X)$, such that $res \circ bc = deg(L/K)\cdot id$,
and $bc\circ res(x) = n \cdot \sum_{\gamma \in Gal(L/K)}
\gamma(x)$, where $n$ is the inseparability degree of the
extension $L/K$. This implies our claim in this case; injectivity
of $bc$ for any finite extension follows.

If $L=K(\alpha)$ where $\alpha$ is transcendental over $K$, then
$K(X)\iso K(X_L)$; this follows from the
excision
sequence
$\pl _{t\in \Aone}
K(X\tim t)\to K(X\times \Aone) \to
K(X_{K(\alpha)})\to 0
$
(where $t$ runs over the  closed
points in $\Aone_K$),
since the first map is zero and
$
K(X\times \Aone)
\cong
K(X)
$.

If $L$ is finitely generated over $K$, so that there exists a
purely transcendental subextension $K\subset K'\subset L$ with
$|L/K|<\infty$, then injectivity follows by comparing the previous
two special cases; if $L/K'$ is normal we also get the description
of the image of $bc$.

Finally, the general case follows from the case of a finitely
generated extension by passing to the limit.
 \epf

\lemm
\label{Tate}
For all $i$
the Galois group $Gal(\overline{\KK}/\KK)$ acts on the
$l$-adic cohomology $\Hh_{et}^{2i}(\BB_{\chi}^{\barr\KK},\Ql(-i))$
through  a finite quotient.


\pf
The cycle map
$c_{\Ql}:{A_i(\BB_{\chi}^{\barr\KK})}_\Ql \to
\Hh^{2i}_{et}(\BB_{\chi}^{\barr\KK}, \Ql(- i))^*
$,
defined by
$\langle c_{\Ql}([Z]),h\rangle =\int h|_Z$ for an $i$-dimensional
 cycle $Z$
(here $\int:
\Hh^{2i}_{et}
(Z,\Ql(-i))\to \Ql$ is the canonical map),
is compatible with the $Gal(\Kbar/\KK)$ action.
It is an isomorphism since
$\bar \KK\cong \Ce$ and the results of \cite{DLP} show
that the cycle map
$c:A_i(\BB_{\chi}^\Ce)
\to
\Hh_{2i}(\BB_{\chi}^\Ce,\Zet)$ is an isomorphism.

In order to factor
the action of $Gal(\Kbar/\KK)$ on
$A_*(\BB_{\chi}^{\barr\KK})$
through $Gal(\KK'/\KK)$
we choose  a finite
set of cycles $Z_i$
whose classes form a basis in $A_*(\BB_{\chi}^\Ce)_\Qu$,
and then a finite
 subextension $\KK'\subset \Kbar$ such that all $Z_i$ are defined
over $\KK'$.
\epf

\sss{
Proof of
\eqref{etaCe}
}\label{pretaCe}
Lemma \ref{KL} says that
$
K(\BB_{\chi}^\KK)_\Qu
=
K(\BB_{\chi}^{\barr\KK})_\Qu^{\ Gal(\Kbar/\KK)}
$
so it suffices to see that the Galois action on
$K(\BB_{\chi}^{\barr\KK})_\Qu
$ is trivial.
However, \ref{implies} and the proof of \ref{Tate}
provide
$
Gal(\Kbar/\KK)
$-equivariant isomorphisms
$
K(\BB_{\chi}^{\barr\KK})_\Qu
\
\bb{ \tau   
}\con
\
A_\bu(\BB_{\chi}^{\barr\KK})_\Qu
\
\bb{c_{\Ql}}\con
\
\Hh_{et}^\bu(\BB_{\chi}^\Kbar,\Ql(-i))^*$.
\epf

\nc{\Eul}{{\EE ul}}

\sus{
The specialization map in \ref{spe}(a) is injective
}
\lab{Specialization is injective}
For this we will use the pairing of $K$-groups of
the Springer fiber and of the Slodowy variety.
Let $\fX$ be a proper variety over a field
$k$, and $i:\fX\imbed \fY$ be a closed embedding, where
$\fY$ is smooth
over $k$. We have a bilinear pairing $\Eul=\Eul_k: K(\fY)\times
K(\fX)\to \Z$, where $\Eul([\F],[\G])$ is the Euler characteristic
of $Ext^\bu(\F,i_*\G)$.

Let us now return to the situation of  \ref{specsec}, and assume
that $X$ is proper over $\OO$, and that $i:X\imbed Y$ is a closed
embedding, where $Y$ is smooth over $\OO$.
For $a\in K(Y^\eta)$, $b\in K(X^\eta)$ we have
$$\Eul_s(sp(a),
sp(b))=\Eul_\eta(a,b)
$$
since  $
(\Ll i_s^*)\RHHom(\FF,\GG)
\cong
\RHHom(\Ll i_s^*\FF,\Ll i_s^*\GG)
$
for $\FF\in \Db(\Coh(Y)),\GG\in \Db(\Coh(Y))$.
In particular, if the pairing $\Eul_\eta$
is non-degenerate in the second variable,
specialization $sp:K(X^\eta)\to
K(X^s)$ is injective.

Since the Slodowy scheme
$\tii S_\OO$ is smooth (in particular flat) over $\OO$
(Lemma \ref{Springer Slodowy}),
we can apply these considerations to
$X=\BB_{\chi}^\OO$, and $Y=\tii S_\OO$.
It is proved in  \cite{Lu} II, Theorem 2.5,
that the pairing
$(\Eul_\Ce)_\Qu: K(Y^\Ce)_ \Qu \times
K(X^\Ce)_ \Qu \to \Qu$ is non-degenerate.
Since
$K(X^\eta)_\Q\iso K(X^\Ce)_\Q$
is proved in
\ref{Specialization is injective}
and the same argument shows that
 $K(Y^\eta)_\Q\iso K(Y^\Ce)_\Q$;
the pairing $\Eul_\eta$ is
also non-degenerate  and then
$sp$ is injective.
\epf

\begin{Rem} The proof of
Lemma \ref{nondegp} below can be adapted  to give a proof  that
$\Eul_\k$ is non-degenerate if  $\k$ has large positive characteristic.
One can then deduce that the same holds
for $\k=\C$. This would give  an alternative proof
 of the result from
\cite{Lu} II mentioned above.
\end{Rem}

\sus{
Upper bound on the
$K$-group
}
Here we use
another
Euler pairing
to prove  that
\begin{equation}\label{ineq}
\dim_\Q K(\BB_{\chi}^\k)_ \Qu\ \le\ \dim_\Q \Hh^\bu (\BB_{\chi}^\Ce,\Q).
\end{equation}
Besides $K(\fX)=K(\Coh(\fX))$ one can consider
$K^0(\fX)$, the Grothendieck group of vector bundles (equivalently,
of complexes of finite homological dimension) on $\fX$.
When $\fX$ is proper over a field we have
another Euler
pairing $Eul_\fX: K^0(\fX)\times
K(\fX) \to \Z$ by $Eul_\fX([\FF],[\GG])=[\RHom(\FF,\GG)]$.

\lemm
\label{nondegp}
The Euler
pairing $Eul_\fX$  for $\fX={\BB_{\chi}^\k}$
is non-degenerate in the second factor,
i.e. it yields an injective map $K(\fX)\imbed \Hom(K^0(\fX),\Zet)$.

\pf
Let $\BB_\chi\aa{\io}\inj\Bhat_\chi$ be the formal neighborhood of $\BB_\chi$
in $T^*\BB$.
For any vector bundle $V$ on $\Bhat_\chi$ and
$\GG\in \Db({\BB_\chi})$, one has
$\RHom^\bu(V,\io_*\GG)
\cong
\RHom^\bu(\io^*V,\GG)$.
So it suffices
to show  that the Euler pairing
$\Eul:
K(\Bhat_\chi)\times K(\BB_\chi)
\to \Zet,\
\Eul([V],[\GG])=\
[\RHom^\bu(V,\io_*\GG)]$, is a perfect pairing.

Let us interpret this pairing  using  localization.
The first of the
isomorphisms
(see  \ref{Categories with chi} for notations)
$$
K(\BB_{\chi})
\cong
K(mod^{fl}(U_{\hatt\chi}^{ 0}))
\aand
K(\Coh(\Bhat_\chi))
\cong
K(mod^{fg}(U_{\hatt\chi}^{ 0}))
,$$
comes from Theorem \ref{better} (notice that
$mod^{fl}(U_{\hatt\chi}^{ 0})= mod_{\chi}(U^0)$, see
 \ref{Categories with chi}), and the second one
from  Theorem
\ref{even better} (notice that $K^0(\Bhat_\chi)\con K(\Bhat_\chi)$
because $T^*\BB$ is smooth).
The above Euler pairing now becomes
the
Euler pairing
 $$
K(mod^{fg}(U_{\hatt \chi}^0))\times  K(mod^{fl}(U_{\hatt \chi}^0))\to \Zet
.$$
However,  the completion  $U_{\hatt \chi}^0$
of $U^0$ at $\chi$  is a complete multi-local algebra
of finite homological dimension: this follows from
 finiteness of homological dimension of
$U^0$, which is clear from Theorem \ref{Localizationthm}. Thus
the latter
 pairing is perfect, because the classes of irreducible and
of indecomposable projective modules provide dual bases in
$ K(mod^{fl}(U_{\hatt \chi}^0)) $ and $K(mod^{fg}(U_{\hatt \chi}^0))$
respectively.
\epf

\lemm
\label{less}
 If $\fX$ is a projective variety over a field, such that
the pairing $Eul_\fX$ is non-degenerate in the second factor $K(\fX)$,
then the
following composition
of the modified Chern character $\tau$ and the $l$-adic cycle map
$c_{\Ql}$,
is injective
$$
K(\fX)_\Ql @>{\tau}>> A_\bu(\fX)_ \Ql
@>{c_{\Ql}}>>
\oplusl_i (\Hh^{2i}_{et}(\fX,\Ql(-i)))^*
.$$

\pf
The pairing $Eul_\fX$ factors through the
modified Chern character by the
Riemann-Roch-Grothendieck Theorem
\cite{Fu} 18.3, and then through the cycle map
by \cite{Fu} Proposition 19.1.2, and the text after Lemma
19.1.2
(this reference uses the cycle map for complex varieties and ordinary
Borel-Moore homology, however the proofs adjust to the
$l$-adic cycle map).
\epf

\lemm
\label{dimindepe}
$
\dim_\Qlb\ \Hh^*_{et}(\BB_{\chi}^\k,\Qlb)
=
\dim_\Q\ \Hh^*(\BB_{\chi}^\Ce,\Q)
$.

{\em Proof.\footnote{This argument was  explained to us by
Michael Finkelberg.}}
Since the
decomposition of the Springer sheaf into irreducible perverse sheaves
is independent of $p$,
the calculation
of  the cohomology of
Springer fibers (i.e., the stalks of the Springer sheaf),
reduces to the calculation of
stalks of intersection cohomology sheaves
of irreducible local systems
on nilpotent orbits. However, Lusztig proved that
the latter one is
independent of $p$
for good $p$ (\cite{Lu2} section 24, in particular
Theorem 24.8 and Subsection 24.10).
\epf

\sss{
Proof of the upper bound
\eqref{ineq}
}
Lemmas
\ref{nondegp} and \ref{less} give the embedding
$
K(\BB_{\chi}^\k)_\Ql
@>{c_\Ql\ci \tau}>>
\oplusl_i
\Hh^{2i}_{et}(\BB_{\chi}^\k,\Ql(-i))^*
$.
Together with
Lemma
\ref{dimindepe}
this gives
$
\dim_\Q K(\BB_{\chi}^\k)_\Q
\le\
\dim_\Ql
\Hh^{*}_{et}(\BB_{\chi}^\k,\Ql(-i))
=
dim_\Q\Hh^*(\BB_{\chi}^\Ce,\Q)
$.
\epf

\sss{
End of the proof of Proposition \ref{spe}
}
We compare the K-groups via
$$
K(\BB_{\chi}^\C)_\Q
\cong
K(\BB_{\chi}^{\barr\KK})_\Q
\bb{\cong}{@<{bc_\KK^{\barr\KK}}<<}
K(\BB_{\chi}^\eta)_\Q
\aa{sp}\inj
K(\BB_{\chi}^s)_\Q
\aa{bc_{k_s}^\k}\inj
K(\BB_{\chi}^\k)_\Q
.$$
The first two isomorphisms are a particular case of
\eqref{etaCe} proved in \ref{pretaCe};
specialization is injective by
\ref{Specialization is injective},
and the base change
$
bc_{k_s}^\k
$
is injective by
Lemma
\ref{KL}.
Actually, all maps have to be isomorphisms since
\eqref{ineq}
says that
$\dim_\Q K(\BB_{\chi}^\k)_\Q$ is bounded above by
$\dim_\Q \Hh^\bu (\BB_{\chi}^\Ce,\Q)
=
\dim_\Q
K(\BB_{\chi}^\C)_\Q
$.
\epf

\end{document}


\maketitle

We compute for $\frakg=\fraksl(3)$
the coherent sheaves corresponding to
the irreducible $U_{\hat{0}}^{0}$-modules and their projective
covers under the equivalence
$
\calD^{b}\Coh_{\calB^{(1)}}(\wcalN^{(1)}) \overset{\Upsilon}{\to}
\calD^{b}(U_{\hat{0}}^{0}\text{-}\Mod^{fg})
$.
We use the normalization of this equivalence from
the subsection 5.3.3, so that
for every $\calF
\in \text{Coh}(\calB^{(1)})$ we have
$\Upsilon(i_{*}\calF)=R\Gamma(\calB,\Fr_{\calB}^{*}\calF)$.

\section{Notations}
We keep the notations of the article, with $G=\SL(3,\bk)$, and
denote $\alpha_1,\alpha_2$ the simple roots of $G$
and $\omega_{1},\omega_{2}$ the fundamental weights. Let
$s_j$ be the reflection $s_{\alpha_{j}} \in W$. We denote by
$\calB \overset{i}\to\wcalN \overset{p}\to\calB$
the inclusion of the zero section
and the natural projection. There are two
natural maps $\pi_{j} : \calB \to \mathbb{P}^{2}$ mapping a flag
$0 \subset V_{1} \subset
V_{2} \subset \bk^{3}$ to $V_{j},\ j=1,2$.  For $n \in \mathbb{Z}$ we have
isomorphisms:
$\pi_{i}^{*}\calO_{\bbP^{2}}(n) \cong \calO_{\calB}(n
\omega_{i}),\ i=1,2$,
and
$\Fr_\calB^*\calO_{\calB^{(1)}}(\lambda)\cong
\calO_{\calB}(p\lambda)$
for
$
\lambda\in \Lambda$.
We will study
irreducible $G$-modules
$L(\lambda)$ of highest weight $\lambda$ for
reduced dominant weights $\lambda$ in $W_{\aff}' \bullet 0$.
Recall the
exact sequence

 $(*) \ 0 \to \Omega_{\bbP^{2}}^{1} \to
\calO_{\bbP^{2}}(-1)^{\oplus 3} \to \calO_{\bbP^{2}} \to 0.$

For simplicity, in what follows we will omit the Frobenius
twist ${}^{(1)}$ (except in the proof of theorem \ref{thmirred},
where we have to be more careful); it should appear on
(almost) every variety we consider.

\section{Irreducible modules}

\begin{theorem}\label{thmirred}

The irreducible $U_{\hat{0}}^{0}$-modules and the corresponding
coherent sheaves are :
\medskip

\begin{center} \begin{tabular}{|c|c||c|c|}
  \hline
  $L(0)=\bk$ & $i_{*}\calO_{\calB}$ & $L((p-2)\omega_{1} + \omega_{2})$ & $i_{*}\pi_{1}^{*}(\Omega_{\bbP^{2}}^{1}(1))[1]$ \\
  $L((p-3)\omega_{2})$ & $i_{*}\calO_{\calB}(-\omega_{1})[2]$ & $L(\omega_{1} + (p-2)\omega_{2})$ & $i_{*}\pi_{2}^{*}(\Omega_{\bbP^{2}}^{1}(1))[1]$ \\
  $L((p-3)\omega_{1})$ & $i_{*}\calO_{\calB}(-\omega_{2})[2]$ &  $L((p-2)\rho)$ & $\calL$ \\
  \hline
\end{tabular} \end{center}
\medskip where
$\calL$
is the cone of the
 only (up to a constant) nonzero morphism
$i_{*}\calO_{\calB}\to
i_{*}\calO_{\calB}(-\rho)[3]$.

\end{theorem}

\emph{Proof.}\
We have
$\Upsilon(i_{*}\calO_{\calB})=R\Gamma(\calB,\calO_\calB)=\bk$.
Also, $\Upsilon (i_*\calO_{\calB}(-\omega_j))=R\Gamma
(\calO_{{\bbP}^2}(-p))$, which
gives the claim for $L((p-3)\omega_{j})$,
$j=1,2$.

Similarly
$\Upsilon(i_{*}\pi_{1}^{*}(\Omega_{(\bbP^{2})^{(1)}}^{1}(1))[1])=
R\Gamma(\calB,\Fr_{\calB}^{*}\pi_{1}^{*}
(\Omega_{(\bbP^{2})^{(1)}}^{1}(1)))[1]$.
Using the exact sequence $(*)$ we obtain a distinguished triangle
$R\Gamma(\calB,\calO_{\calB})^{\oplus 3} \to
R\Gamma(\calB,\calO_{\calB}(p\omega_{1})) \to
\Upsilon(i_{*}\pi_{1}^{*}(\Omega_{(\bbP^{2})^{(1)}}^{1}(1))[1]).$
Here the second arrow is the inclusion of
$G$-modules $L(\omega_{1})^{(1)} \hookrightarrow
H^{0}(p\omega_{1})$. Hence
$\Upsilon(i_{*}\pi_{1}^{*}(\Omega_{(\bbP^{2})^{(1)}}^{1}(1))[1])
\cong L((p-2)\omega_{1} + \omega_{2}).$ The claim for
$L(\omega_{1} + (p-2)\omega_{2})$ follows by applying the
outer automorphism of $\fraksl(3)$.

Finally, the last irreducible module $L((p-2)\rho) $
is a quotient of the Weyl module
$[H^{0}((p-2)\rho)]^{*}$, moreover,
we have a short exact sequence
$
0 \to \bk \to [H^{0}((p-2)\rho)]^{*} \to
L((p-2)\rho) \to 0
.$ Applying $\Upsilon^{-1}$, we get  distinguished triangle
$i_{*}\calO_{\calB}\to
i_{*}\calO_{\calB}(-\rho)[3]\to \calL$, where we used
that $\Upsilon(i_{*}\calO_{\calB^{(1)}}(-\rho))
=R\Gamma(\calB, \calO_{\calB}(-p\rho))
=[H^{0}((p-2)\rho)]^{*}[-3]$
by Serre duality. Since $\Hom (\bk,  [H^{0}((p-2)\rho)]^{*})$
is one dimensional, we see that the first arrow in this
triangle is the unique (up to a constant) map between the two objects.
\qed

\emph{Remark.}
We have just shown, using  equivalence $\Upsilon$, that
 $\Ext_{\wcalN}^{3}(i_{*}\calO_{\calB},i_{*}\calO_{\calB}(-\rho))$
is one dimensional. One can
compute this Ext group more directly:
using the Koszul resolution of
$\calO_{\calB}$ over ${\rm S}(\calT_{\calB})$ one can identify it with
$H^{3}(-\rho) \oplus H^{2}(\Omega^{1}_{\calB} (-\rho)) \oplus
H^{1}(\Omega^{2}_{\calB}(-\rho)) \oplus
H^{0}(\Omega^{3}_{\calB}(-\rho))$. One can show
that $H^{3}(-\rho)$, $H^{0}(\Omega^{3}_{\calB}(-\rho))$
and $H^{1}(\Omega^{2}_{\calB}(-\rho))$ are $0$, while
 $H^{2}(\Omega^{1}_{\calB}(-\rho))\cong \bk$ (by Serre
duality the last claim is equivalent to $H^1(\calT_\calB(-\rho))
\cong \bk$, which is checked below).

\section{Projective covers}

\begin{theorem}

The coherent sheaves corresponding to the projective covers of the
irreducible modules are : \medskip

\begin{center} \begin{tabular}{|c|c||c|c|}
  \hline
  $i_{*}\calO_{\calB}$ & $\calP$ & $i_{*}\pi_{1}^{*}(\Omega_{\bbP^{2}}^{1}(1))[1]$ & $\calO_{\wcalN}(\omega_{1})$ \\
  $i_{*}\calO_{\calB}(-\omega_{1})[2]$ & $p^{*}((\pi_{2}^{*}\Omega_{\bbP^{2}}^{1})(\omega_{1} + 2\omega_{2}))$
  & $i_{*}\pi_{2}^{*}(\Omega_{\bbP^{2}}^{1}(1))[1]$ & $\calO_{\wcalN}(\omega_{2})$ \\
  $i_{*}\calO_{\calB}(-\omega_{2})[2]$ & $p^{*}((\pi_{1}^{*}\Omega_{\bbP^{2}}^{1})(2\omega_{1} + \omega_{2}))$
  & $\calL$ & $\calO_{\wcalN}(\rho)$ \\
  \hline
\end{tabular} \end{center}
\medskip where
$\calP$ is the nontrivial extension of
$
\calO_{\wcalN}(\rho)$
by
$\calO_{\wcalN}$
given by a non-zero element in the one dimensional space
$H^{1}(\calT_{\calB}(-\rho))\subset H^1(\calO_{\wcalN}(-\rho))$.

\end{theorem}

\emph{Remark.} In fact, the sheaves corresponding to the projective
covers are vector bundles on the  formal completion of
$\wcalN$ at $\calB$. The objects displayed in the above
table are vector bundles on $\wcalN$. The former are obtained
from the latter by pull-back to the formal completion.
\medskip

\emph{Proof.}\  We only have to check that for each
$\calP_{i}$ in the
list and each irreducible $\calL_{j}$, we have
$\Ext^{*}_{\wcalN}(\calP_{i},\calL_{j})=\bk^{\delta_{ij}}$. Let us begin
with $\calO_{\wcalN}(\rho)$. We have
$\Ext_{\wcalN}^{*}(\calO_{\wcalN}(\rho),i_{*}\calO_{\calB}) \cong
\Ext^{*}_{\calB}(\calO_{\calB}(\rho),\calO_{\calB})\cong
H^{*}(\calB,\calO_{\calB}(-\rho))=0$ by adjunction. Similarly for
$i_{*}\calO_{\calB}(-\omega_{j})[2]$ ($j=1,2$). The sequence $(*)$
gives
$
\Ext_{\wcalN}^{*}
(\calO_{\wcalN}(\rho),i_{*}\pi_{j}^{*}(\Omega_{\bbP^{2}}^{1}(1))[1])
=
\Ext_{\calB}^{*}(\calO_{\calB}(\rho),\pi_{j}^{*}(\Omega_{\bbP^{2}}^{1}(1))[1])
=0$
($j=1,2$). Using the distinguished triangle from the definition of
 $\calL$ we get
$\Ext_{\wcalN}^{*}(\calO_{\wcalN}(\rho),\calL)=\bk$. The cases of
$\calO_{\wcalN}(\omega_{j})$ ($j=1,2$) are similar.

Now let us consider
$p^{*}((\pi_{1}^{*}\Omega_{\bbP^{2}}^{1})(2\omega_{1} +
\omega_{2}))$. The exact sequence $(*)$ and Borel-Weil-Bott Theorem
\cite{Jagps} give the
result for the first 5 irreducible modules. For $\calL$, we have
$\Ext^{*}_{\calB}((\pi_{1}^{*}\Omega_{\bbP^{2}}^{1})(2\omega_{1} +
\omega_{2}),\calO_{\calB})=0$, and in computing
$\Ext^{*}_{\calB}((\pi_{1}^{*}\Omega_{\bbP^{2}}^{1})(2\omega_{1} +
\omega_{2}),\calO_{\calB}(-\rho)[3])$, two non-zero modules appear
in degree $0$ : $[H^{3}(-2\rho)]^{\oplus 3}$ and $H^{0}(\omega_{1})$. The
map between these two modules is an isomorphism as in the proof of
Theorem \ref{thmirred}, hence
$\Ext^{*}_{\wcalN}(p^{*}((\pi_{1}^{*}\Omega_{\bbP^{2}}^{1})(2\omega_{1}
+ \omega_{2})),\calL)=0$.

We claim that $H^{1}(\calT_{\calB}(-\rho)) \cong \bk$,
this follows by the Borel-Weil-Bott Theorem
from the exact sequence $0\to \calO_{\calB}
(\alpha_1)\to \calT_{\calB} \to \pi_2^*(\calT_{{\bbP}^2})\to 0$,
and vanishing of $R\Gamma( \pi_2^*(\calT_{{\bbP}^2})(-\rho))$
(see, e.g., \cite{D}).
Thus we have the line
$H^{1}(\calT_{\calB}(-\rho))
\subset
H^{1}({\rm S}(\calT_{\calB})(-\rho))
=
\Ext_{\wcalN}^{1}(\calO_{\wcalN}(\rho),\calO_{\wcalN})
$, which  defines a triangle
$\calO_{\wcalN} \to \calP \to \calO_{\wcalN}(\rho)$.
Standard calculations give
the result for $\calP$ and the first three irreducible modules. The triangle
defining $\calP$ gives
$
\Ext_{\wcalN}^{*}(\calP,i_{*}\pi_{1}^{*}(\Omega_{\bbP^{2}}^{1}(1))[1])=
H^{*}(\pi_{1}^{*}(\Omega_{\bbP^{2}}^{1}(1)))[1]
.$
Using $(*)$, we have an exact sequence $0 \to
H^{0}(\pi_{1}^{*}(\Omega_{\bbP^{2}}^{1}(1))) \to \bk^{3} \to
H^{0}(\omega_{1}) \to H^{1}(\pi_{1}^{*}(\Omega_{\bbP^{2}}^{1}(1))) \to 0$
with invertible middle arrow
(the other cohomology modules vanish).

Finally, let us show that $\Ext^{*}_{\wcalN}(\calP,\calL)=0$.
We have  $ R\Hom_{\wcalN}(\calP,i_{*}\calO_{\calB}) \cong
R\Gamma(\calO_{\calB})\cong \bk, \;
R\Hom_{\wcalN}(\calP,i_{*}\calO_{\calB}(-\rho)[3]) \cong
R\Gamma(\calO_{\calB}(-2\rho)[3])\cong \bk$, thus we only need
to check that for nonzero morphisms $b : i_{*}\calO_{\calB} \to
i_{*}\calO_{\calB}(-\rho)[3]$, $\phi:\calP\to i_{*}\calO_\calB$
we have $b\circ \phi \ne 0$. It is clear from Remark after Theorem
\ref{thmirred} that $b=i_*(\beta)\circ \delta$, where
$\delta:i_*\calO_\calB\to i_*\calT_\calB[1]$ is the class of the extension
$0\to  i_*\calT_\calB\to \calO_{\wcalN}/\calJ_\calB^2\to i_*\calO_\calB\to 0$,
and
$\beta: \calT_\calB [1]\to \calO_\calB(-\rho)[3]$ is a non-zero morphism;
here $\calJ_\calB$ is the ideal sheaf on the zero section in $\wcalN$.

We claim that $\delta\circ \phi=i_*(\gamma) \circ \psi$, where
$\psi:\calP\twoheadrightarrow i_*\calO_\calB(\rho)$ and $\gamma:
\calO_\calB(\rho)\to \calT_\calB[1]$
 are nonzero morphisms. This follows from the definition
of $\calP$, which implies that $\calP$ has a quotient,
which is an extension of $i_*\calO_\calB\oplus i_*\calO_\calB(\rho)$
by $i_*\calT_\calB$, such that the corresponding class
in $\Ext^1(i_*\calO_\calB, i_*(\calT_\calB))$ equals $\delta$,
while the corresponding class in  $\Ext^1(i_*\calO_\calB(\rho),
i_*(\calT_\calB))$ is non-trivial and is an image under $i_*$
of an extension of coherent sheaves on $\calB$.

It remains to show that the composition $i_*\beta\circ i_*\gamma
\circ \psi$ is nonzero. The composition $\beta\circ \gamma \in
\Ext^3(\calO_\calB(\rho),\calO_\calB(-\rho))=H^3(\calB,
\calO(-2\rho))= \bk$ is nonzero, because it coincides with the
Serre duality pairing of nonzero elements $\beta$, $\gamma$ in the
two dual one-dimensional spaces $H^1(\calT_\calB(-\rho))$,
$H^2(\calT_\calB^*(-\rho))$. Consequently, the composition
$i_*(\beta\circ \gamma) \circ \psi$ is also nonzero, since under
the isomorphism $\Hom (\calP, i_*\calO_\calB(-\rho)[3])\cong \Hom
(i^*\calP, \calO_\calB(-\rho)[3])\cong \Hom (\calO_\calB\oplus
\calO_\calB(\rho), \calO_\calB(-\rho)[3])$ it corresponds to the
composition of
 $\beta\circ \gamma$ and projection
to the second summand.
 \qed